\documentclass[10pt]{articleHJ}
\usepackage[english]{babel}
\usepackage{amsmath, amsthm}
\usepackage{mathtools}

\usepackage{amssymb}

\usepackage[text={6.0in,8.6in},centering,letterpaper]{geometry}

  \usepackage{paralist}
  \usepackage{graphics} 
  \usepackage{epsfig} 
\usepackage{graphicx}  
\usepackage[colorlinks=true]{hyperref}
\hypersetup{urlcolor=blue, citecolor=red}

\newlength{\defbaselineskip}
\newcommand{\setlinespacing}[1]%
           {\setlength{\baselineskip}{#1 \defbaselineskip}}

\newcommand{\singlespacing}{\setlength{\baselineskip}{\defbaselineskip}}

\newcommand{\Z}{\ensuremath{\mathbb{Z}}}

\newcommand{\R}{\ensuremath{\mathbb{R}}}
\newcommand{\C}{\ensuremath{\mathbb{C}}}

\renewcommand{\epsilon}{\ensuremath{\varepsilon}}

\renewcommand{\Re}{\mathop\mathrm{Re}\nolimits}
\renewcommand{\Im}{\mathop\mathrm{Im}\nolimits}
\renewcommand{\ker}{\ensuremath{\mathrm{ker}}}
\newcommand{\Range}{\ensuremath{\mathrm{Range}}}

\newcommand{\codim}{\ensuremath{\mathrm{codim}}}
\renewcommand{\span}{\ensuremath{\mathrm{span}}}

\DeclarePairedDelimiter\abs{\lvert}{\rvert}

\DeclarePairedDelimiter{\ip}\langle\rangle
\DeclarePairedDelimiter{\nrm}\lVert\rVert

\theoremstyle{plain}
\newtheorem{theorem}{Theorem}[section]
\newtheorem{proposition}[theorem]{Proposition}
\newtheorem{lemma}[theorem]{Lemma}
\newtheorem{corollary}[theorem]{Corollary}

\theoremstyle{definition}

\newtheorem{remark}[theorem]{Remark}

\newtheorem*{assumption*}{\assumptionnumber}
\providecommand{\assumptionnumber}{}
\makeatletter
\newenvironment{assumption}[2]
 {%
  \renewcommand{\assumptionnumber}{Assumption ($#1$)}%
  \begin{assumption*}%
  \protected@edef\@currentlabel{$#1$}%
 }
 {%
  \end{assumption*}
 }
\makeatother
\newcommand{\asref}[2]{\ref{$#1$}}

\numberwithin{equation}{section}

\begin{document}

\begin{frontmatter}

\title{Nonlinear stability of pulse solutions for the discrete FitzHugh-Nagumo equation with infinite-range interactions}
\journal{...}
\author[LDA]{W. M. Schouten-Straatman\corauthref{coraut}},
\corauth[coraut]{Corresponding author. }
\author[LDB]{H. J. Hupkes}
\address[LDA]{
  Mathematisch Instituut - Universiteit Leiden \\
  P.O. Box 9512; 2300 RA Leiden; The Netherlands \\ Email:  {\normalfont{\texttt{w.m.schouten@math.leidenuniv.nl}}}
}
\address[LDB]{
  Mathematisch Instituut - Universiteit Leiden \\
  P.O. Box 9512; 2300 RA Leiden; The Netherlands \\ Email:  {\normalfont{\texttt{hhupkes@math.leidenuniv.nl}}}
}

\date{\today}
\journal{Discrete  and Continuous Dynamical Systems A}
\begin{abstract}
\singlespacing

We establish the existence and nonlinear stability of travelling pulse solutions for the discrete FitzHugh-Nagumo equation with
infinite-range interactions close to the continuum limit. For the verification of the spectral properties, we need
to study a functional differential equation of mixed type (MFDE) with unbounded shifts. We avoid the use of exponential
dichotomies and phase spaces, by building on a technique developed by Bates, Chen and Chmaj for the discrete Nagumo equation.
This allows us to transfer several crucial Fredholm properties from the PDE setting to our discrete setting.

\end{abstract}



\end{frontmatter}

\section{Introduction}\label{introduction}
The FitzHugh-Nagumo partial differential equation (PDE) is given by
\begin{equation}\label{firstfirstequation}\begin{array}{lcl}u_t&=&u_{xx}+g(u;r_0)-w\\[0.2cm]
w_t&=&\rho(u-\gamma w),\end{array}\end{equation}
where $g(\cdot;r_0)$ is the cubic bistable nonlinearity given by
\begin{equation}\begin{array}{lcl}g(u;r_0)&=&u(1-u)(u-r_0)\end{array}\end{equation}
and $\rho,\gamma$ are positive constants.
This PDE is commonly used as a simplification of the Hodgkin-Huxley equations, which describe the propagation of signals through nerve fibres.
The spatially homogeneous version of this equation was first stated by FitzHugh in 1961 \cite{FITZHUGH19662} in order to describe the potential felt
at a single point along a nerve axon as a signal travels by. A few years later \cite{FITZHUGH1966}, the diffusion term in (\ref{firstfirstequation}) was
added to describe the dynamics of the full nerve axon instead of just a single point. As early as 1968 \cite{FITZHUGH1968}, FitzHugh released a computer
animation based on numerical simulations of (\ref{firstfirstequation}). This video clip clearly shows that (\ref{firstfirstequation}) admits
isolated pulse solutions resembling the spike signals that were measured by Hodgkin and Huxley in the nerve fibres of giant squids \cite{HODHUX1952}.\par

As a consequence of this rich behaviour and the relative simplicity of its structure, (\ref{firstfirstequation}) has served as a prototype for several similar systems. For example, memory devices have been designed using a planar version of (\ref{firstfirstequation}), which supports stable stationary, radially symmetric spot patterns \cite{KAMIN2006}. In addition, gas discharges have been described using a three-component FitzHugh-Nagumo system \cite{ORG1998,SCHENK1997}, for which it is possible to  find stable travelling spots \cite{HEIJSTER2014}.\par

Mathematically, it turned out to be a major challenge to control the interplay between the excitation and recovery dynamics and rigorously construct the travelling pulses visualized by FitzHugh in \cite{FITZHUGH1968}. Such pulse solutions have the form
\begin{equation}\begin{array}{lcl}(u,w)(x,t)&=&(\overline{u}_0,\overline{w}_0)(x+c_0 t),\end{array}\end{equation}
in which $c_0$ is the wavespeed and the wave profile $(\overline{u}_0,\overline{w}_0)$ satisfies the limits
\begin{equation}\begin{array}{lcl}\lim\limits_{|\xi|\rightarrow\infty}(\overline{u}_0,\overline{w}_0)(\xi)&=&0.\end{array}\end{equation}
Plugging this Ansatz into (\ref{firstfirstequation}) and writing $\xi=x+c_0t$, we see that the profiles are homoclinic solutions to the travelling wave ordinary differential equation (ODE)
\begin{equation}\label{firstodequation}\begin{array}{lcl}c_0\overline{u}_0'(\xi)&=&\overline{u}_0''(\xi)+g(\overline{u}_0(\xi);r_0)-\overline{w}_0(\xi)\\[0.2cm]
c_0\overline{w}_0'(\xi)&=&\rho\big[\overline{u}_0(\xi)-\gamma \overline{w}_0(\xi)\big].\end{array}\end{equation}
The analysis of this equation in the singular limit $\rho\downarrow 0$ led to the birth of many techniques in
geometric singular perturbation theory, see for example \cite{JONES1995} for an interesting overview.
Indeed, the early works \cite{CARP1977,HAST1976,JONES1984,JONESKOPLAN1991} used geometric techniques such as the Conley index, exchange lemmas and differential forms to construct pulses and analyze their stability. A more analytic approach was later developed in \cite{KRUSANSZM1997}, where Lin's method was used in the $r_0\approx\frac{1}{2}$ regime to connect a branch of so-called slow-pulse solutions to (\ref{firstodequation}) to a branch of fast-pulse solutions. This equation is still under active investigation, see for example \cite{CARTER2016,CARTER2015}, where the birth of oscillating tails for the pulse solutions is described as the unstable root $r_0$ of the nonlinearity $g$ moves towards the stable root at zero.\par

Many physical, chemical and biological systems have an inherent discrete structure that strongly influences their dynamical behaviour. In such settings \textit{lattice differential equations} (LDEs), i.e. differential equations where the spatial variable can only take values on a lattice such as $\Z^n$, are the natural replacements for PDEs, see for example \cite{BatesInfRange,HJHSTBFHN,MPA}. Although mathematically it is a relatively young field of interest, LDEs have already appeared frequently in the more applied literature. For example, they have been used to describe phase transitions in Ising models \cite{BatesInfRange}, crystal growth in materials \cite{CAHN} and phase mixing in martensitic structures \cite{VAIN2009}.\par

To illustrate these points, let us return to the nerve axon described above and reconsider the propagation of electrical signals through nerve fibres. It is well known that electrical signals can only travel at adequate speeds if the nerve fibre is insulated by a myeline coating. This coating admits regularly spaced gaps at the so-called nodes of Ranvier \cite{RANVIER1878}. Through a process called saltatory conduction, it turns out that excitations of nerves effectively jump from one node to the next \cite{LILLIE1925}. Exploiting this fact, it is possible \cite{EVVPD18} to model this jumping process with the discrete FitzHugh-Nagumo LDE
\begin{equation}\label{ditishetalgemenealgemeneprobleem}\begin{array}{lcl}\dot{u}_j&=&\frac{1}{h^2}(u_{j+1}+u_{j-1}-2u_j)+g(u_j;r_0)-w_j\\[0.2cm]
\dot{w}_j&=&\rho[u_j-\gamma w_j].\end{array}\end{equation}
The variable $u_j$ now represents the potential at the $j^{\text{th}}$ node, while the variable $w_j$ denotes a recovery component. The nonlinearity $g$ describes the ionic interactions. Note that this equation arises directly from the FitzHugh-Nagumo PDE upon taking the nearest-neighbour discretisation of the Laplacian on a grid with spacing $h>0$. \par

Inspired by the procedure for partial differential equations, one can substitute a travelling pulse Ansatz
\begin{equation}
\begin{array}{lcl}
\label{eq:int:trvPulseAnsatz}
(u_j,w_j)(t)&=&(\overline{u}_h,\overline{w}_h)(hj+c_ht)\end{array}
\end{equation}
into (\ref{ditishetalgemenealgemeneprobleem}). Instead of an ODE, we obtain the system
\begin{equation}\begin{array}{lcl}c_h\overline{u}_h'(\xi)&=&\frac{1}{h^2}[\overline{u}_h(\xi +h )+\overline{u}_h(\xi -h)-2\overline{u}_h(\xi)]+g(\overline{u}_h(\xi);r_0)-\overline{w}_h(\xi)\\[0.2cm]
c_h\overline{w}_h'(\xi)&=&\rho[\overline{u}_h(\xi)-\gamma \overline{w}_h(\xi)]\end{array}\end{equation}
in which $\xi=hj+c_ht$. Such equations are called functional differential equations of mixed type (MFDEs), since they contain both advanced (positive) and retarded (negative) shifts.\par
In \cite{HJHSTBFHN,HJHFZHNGM}, Hupkes and Sandstede studied (\ref{ditishetalgemenealgemeneprobleem}) and showed that for small values of $\rho$ and $r_0$ sufficiently far from $\frac{1}{2}$, there exists a locally unique travelling pulse solution of this system and that it is asymptotically stable with an asymptotic phase shift. No restrictions were required on the discretisation distance $h$, but the results relied heavily on the existence of exponential dichotomies for MFDEs. As a consequence, the techniques developed in \cite{HJHFZHNGM,HJHSTBFHN} can only be used if the discretisation involves finitely many neighbours. Such discretisation schemes are said to have finite range.\par

Recently, an active interest has arisen in non-local equations that feature infinite-range interactions. For example, Ising models have been used to describe the infinite-range interactions between magnetic spins arranged on a grid \cite{BatesInfRange}. In addition, many physical systems, such as amorphous semiconductors \cite{GU1996} and liquid crystals \cite{Ciuchi}, feature non-standard diffusion processes, which are generated by fractional Laplacians. Such operators are intrinsically non-local and hence often require infinite-range discretisation schemes \cite{Ciaurri}.\par

Our primary interest here, however, comes from so-called neural field models, which aim to describe the dynamics of large networks of neurons. These neurons interact with each other by exchanging signals across long distances through their interconnecting nerve axons \cite{BRESS2011,BRESS2014,PINTO2001,SNEYD2005}. It is of course a major challenge to find effective equations to describe such complex interactions. One model that has been proposed \cite[Eq. (3.31)]{BRESS2011} features a FitzHugh-Nagumo type system with infinite-range interactions.\par

Motivated by the above, we consider a class of infinite-range FitzHugh-Nagumo LDEs that includes the prototype
\begin{equation}\begin{array}{lcl}\label{infiniterangeversion}\dot{u}_j&=&\frac{\kappa}{h^2}\sum\limits_{k\in\Z_{>0}}e^{-k^2}[u_{j+k}+u_{j-k}-2u_j]+g(u_j;r_0)-w_j\\[0.2cm]
\dot{w}_j&=&\rho[u_j-\gamma w_j],\end{array}\end{equation}
in which $\kappa>0$ is a normalisation constant. In \cite{Faye2015}, Faye and Scheel studied equations such as (\ref{infiniterangeversion}) for discretisations with infinite-range interactions featuring exponential decay in the coupling strength. They circumvented the need to use a state space as in \cite{HJHFZHNGM},
which enabled them to construct pulses to (\ref{infiniterangeversion}) for arbitrary discretisation distance $h$. Very recently \cite{Faye2016}, they developed a center manifold approach that allows bifurcation results to be obtained for neural field equations.\par

In this paper, we also construct pulse solutions to equations such as (\ref{infiniterangeversion}), but under weaker assumptions on the decay rate of the couplings. Moreover, we will establish the nonlinear stability of these pulse solutions, provided the coupling strength decays exponentially. However, both results do require the discretisation distance $h$ to be very small.\par

In particular, we will be working in the continuum limit. The pulses we construct can be seen as perturbations of the travelling pulse solution of the FitzHugh-Nagumo PDE. However, we will see that the travelling wave equations are highly singular perturbations of (\ref{firstodequation}), which poses a significant mathematical challenge. On the other hand, we do not need to use exponential dichotomies directly in our non-local setting as in \cite{HJHSTBFHN}. Instead we are able to exploit the detailed knowledge that has been obtained using these techniques for the pulses in the PDE setting.\par

Our approach to tackle the difficulties arising from this singular perturbation is strongly inspired by the work of Bates, Chen and Chmaj. Indeed, in their excellent paper \cite{BatesInfRange}, they study a class of systems that includes the infinite-range discrete Nagumo equation
\begin{equation}\begin{array}{lcl}\dot{u}_j&=&\frac{\kappa}{h^2}\sum\limits_{k\in\Z_{>0}}e^{-k^2}[u_{j+k}+u_{j-k}-2u_j]+g(u_j;r_0),\end{array}\end{equation}
in which $\kappa>0$ is a normalisation constant. This equation can be seen as a discretisation of the Nagumo PDE
\begin{equation}\label{nagumopde}\begin{array}{lcl}u_t&=&u_{xx}+g(u;r_0).\end{array}\end{equation} The authors show that, under some natural assumptions, these systems admit travelling front solutions for $h$ small enough.\par

In the remainder of this introduction we outline their approach and discuss our modifications, which significantly broaden the application range of these methods. We discuss these modifications for the prototype (\ref{infiniterangeversion}), but naturally they can be applied to a broad class of systems.\par

\subsection*{Transfer of Fredholm properties: Scalar case.}
An important role in \cite{BatesInfRange} is reserved for the operator $\mathcal{L}_{h;\overline{u}_{0:\text{sc}};c_{0:\text{sc}}}$ given by
\begin{equation}\label{scalaroperator}\begin{array}{lcl}\mathcal{L}_{h;\overline{u}_{0:\text{sc}};c_{0:\text{sc}}}v(\xi)&=&c_{0:\text{sc}}v'(\xi)-\frac{\kappa}{h^2}\sum\limits_{k\in\Z_{>0}}e^{-k^2}\Big[v(\xi+hk)+v(\xi-hk)-2v(\xi)\Big]\\[0.2cm]
&&\qquad-g_u(\overline{u}_{0:\text{sc}}(\xi);r_0)v(\xi),\end{array}\end{equation}
where $\overline{u}_{0:\text{sc}}$ is the wave solution of the scalar Nagumo PDE (\ref{nagumopde}) with wavespeed $c_{0:\text{sc}}$. This operator arises as the linearisation of the scalar Nagumo MFDE
\begin{equation}\begin{array}{lcl}c_{0:\text{sc}}u'(\xi)=\frac{\kappa}{h^2}\sum\limits_{k\in\Z_{>0}}e^{-k^2}\Big[v(\xi+hk)+v(\xi-hk)-2v(\xi)\Big]+g_u(\overline{u}_{0:\text{sc}}(\xi);r_0)v(\xi),\end{array}\end{equation}
around the wave solution $\overline{u}_{0:\text{sc}}$ of the scalar Nagumo PDE (\ref{nagumopde}). This operator should be compared to
\begin{equation}\begin{array}{lcl}\mathcal{L}_{0;\overline{u}_{0:\text{sc}};c_{0:\text{sc}}}v(\xi)&=&c_{0:\text{sc}}v'(\xi)-v''(\xi)-g_u(\overline{u}_{0:\text{sc}}(\xi);r_0)v(\xi),\end{array}\end{equation}
the linearisation of the scalar Nagumo PDE around its wave solution.\par

The key contribution in \cite{BatesInfRange} is that the authors fix a constant $\delta>0$ and use the invertibility of $\mathcal{L}_{0;\overline{u}_{0:\text{sc}};c_{0:\text{sc}}}+\delta$ to show that also $\mathcal{L}_{h;\overline{u}_{0:\text{sc}};c_{0:\text{sc}}}+\delta$ is invertible. In particular, they consider weakly-converging sequences $\{v_n\}$ and $\{w_n\}$ with $(\mathcal{L}_{h;\overline{u}_{0:\text{sc}};c_{0:\text{sc}}}+\delta)v_n=w_n$ and try to find a uniform (in $\delta$ and $h$) upper bound for the $L^2$-norm of $v_n'$ in terms of the $L^2$-norm of $w_n$. Such a bound is required to rule out the limitless transfer of energy into oscillatory modes, a key complication when taking weak limits. To obtain this bound, the authors exploit the bistable structure of the nonlinearity $g$ to control the behaviour at $\pm\infty$. This allows the local $L^2$-norm of $v_n$ on a compact set to be uniformly bounded away from zero. Since the operator $\mathcal{L}_{h;\overline{u}_{0:\text{sc}};c_{0:\text{sc}}}+\delta$ is not self-adjoint, this procedure must be repeated for the adjoint operator.\par

\subsection*{Transfer of Fredholm properties: System case.}
Plugging the travelling pulse Ansatz
\begin{equation}\begin{array}{lcl}(u,w)_j(t)&=&(\overline{u}_h,\overline{w}_h)(hj+c_ht)\end{array}\end{equation}
into (\ref{infiniterangeversion}) and writing $\xi=hj+c_ht$, we see that the profiles are homoclinic solutions to the equation
\begin{equation}\label{secondodequation}\begin{array}{lcl}c_h\overline{u}_h'(\xi)&=&\frac{\kappa}{h^2}\sum\limits_{k>0}e^{-k^2}\Big[\overline{u}_h(\xi+kh)+\overline{u}_h(\xi-kh)-2\overline{u}_h(\xi)\Big]\\&&+g(\overline{u}_h(\xi);r_0)-\overline{w}_h(\xi)\\[0.2cm]
c_h\overline{w}_h'(\xi)&=&\rho\Big(\overline{u}_h(\xi)-\gamma \overline{w}_h(\xi)\Big).\end{array}\end{equation}
We start by considering the linearised operator $\mathcal{K}_{h;\overline{u}_0;c_0}$ of the system (\ref{secondodequation}) around the pulse solution $(\overline{u}_0,\overline{w}_0)$ of the FitzHugh-Nagumo PDE with wavespeed $c_0$. This operator is given by
\begin{equation}\label{eerstelinearisation}\begin{array}{lcl}\mathcal{K}_{h;\overline{u}_0;c_0}\left(\begin{array}{l}v\\ w\end{array}\right)(\xi)&=&\left(\begin{array}{l}\mathcal{L}_{h;\overline{u}_0;c_0}v(\xi)+w(\xi)\\
c_0 w'(\xi)-\rho v(\xi)+\rho\gamma w(\xi)\end{array}\right),\end{array}\end{equation}
where $\mathcal{L}_{h;\overline{u}_0;c_0}$ is given by equation (\ref{scalaroperator}), but with $\overline{u}_{0:\text{sc}}$
replaced by $\overline{u}_0$ and $c_{0:\text{sc}}$ by $c_0$. \par

In {\S}\ref{singularoperator}
we use a Fredholm alternative as described above to establish the
invertibility of $\mathcal{K}_{h;\overline{u}_0;c_0}+\delta$ for fixed $\delta>0$. However, the transition from a scalar
equation to a system is far from trivial.
Indeed, when transferring the Fredholm properties there are multiple cross terms that need to be controlled.
We are aided here by the relative simplicity of the terms in the equation that involve $\overline{w}$.
In particular, three of the four matrix-elements
of the linearisation (\ref{eerstelinearisation}) have constant coefficients.
We emphasize that it is not sufficient to merely assume that the limiting state $(0,0)$ is a stable
equilibrium of (\ref{infiniterangeversion}). In \cite{SCHPerExtensions} we explore
a number of structural conditions that allow these types of arguments to be extended to general multi-component systems.

\subsection*{Construction of pulses.}
Using these results for $\mathcal{K}_{h;\overline{u}_0;c_0}$ a standard fixed point argument
can be used to show that the system (\ref{infiniterangeversion}) has a locally unique travelling pulse
solution $(\overline{U}_h(t))_j=(\overline{u}_h,\overline{w}_h)(hj+c_ht)$ for $h$ small enough,
which converges to a travelling pulse solution of the FitzHugh-Nagumo PDE as $h\downarrow 0$.
Indeed, one can mimic the approach developed in \cite[{\S}4]{BatesInfRange},
which in turn closely follows the lines of a standard proof of the implicit function theorem. \par

\subsection*{Spectral stability.}
The natural next step is to study the linear operator $\mathcal{K}_{h;\overline{u}_h;c_h}$ that arises after linearising the system (\ref{infiniterangeversion}) around its new-found pulse solution. This operator is given by
\begin{equation}\begin{array}{lcl}\mathcal{K}_{h;\overline{u}_h;c_h}\left(\begin{array}{l}v\\ w\end{array}\right)(\xi)&=&\left(\begin{array}{l}\mathcal{L}_{h;\overline{u}_h;c_h}v(\xi)+w(\xi)\\
c_0 w'(\xi)-\rho v(\xi)+\rho\gamma w(\xi)\end{array}\right),\end{array}\end{equation}
where $\mathcal{L}_{h;\overline{u}_h}$ is given by equation (\ref{scalaroperator}), but with $\overline{u}_{0:sc}$ replaced by $\overline{u}_h$ and $c_{0:\text{sc}}$ by $c_h$. The procedure above can be repeated to show that for fixed $\delta>0$, it also holds that $\mathcal{K}_{h;\overline{u}_h;c_h}+\delta$ is invertible for $h$ small enough. However, to understand the spectral stability of the pulse, we need to consider the eigenvalue problem
\begin{equation}\label{ditishetprobleem2}\begin{array}{lcl}\mathcal{K}_{h;\overline{u}_h;c_h}v+\lambda v&=&0\end{array}\end{equation}
for fixed values of $h$ and $\lambda$ ranging throughout a half-plane. Switching between these two points of view turns out to
be a delicate task.\par

We start in \S \ref{waveproperties} by showing that $\mathcal{K}_{h;\overline{u}_h;c_h}$ and its adjoint $\mathcal{K}_{h;\overline{u}_h;c_h}^*$ are Fredholm operators with one-dimensional kernels. This is achieved by explicitly constructing a kernel element for $\mathcal{K}_{h;\overline{u}_h;c_h}^*$ that converges to a kernel element of the adjoint of the operator corresponding to the linearised PDE. An abstract perturbation argument then yields the result.\par

In particular, we see that $\lambda=0$ is a simple eigenvalue of $\mathcal{K}_{h;\overline{u}_h;c_h}$. In \S \ref{spectralanalysis} we establish that in a suitable half-plane, the spectrum of this operator consists precisely of the points $\{k2\pi i c_h\frac{1}{h}:k\in\Z\}$, which are all simple eigenvalues. We do this by first showing that the spectrum is invariant under the operation $\lambda\mapsto \lambda+\frac{2\pi i c_h}{h}$, which allows us to restrict ourselves to values of $\lambda$ with imaginary part in between $-\frac{\pi c_h}{h}$ and $\frac{\pi c_h}{h}$. Note that the period of the spectrum is dependent on $h$ and grows to infinity as $h\downarrow 0$. This is not too surprising, since the spectrum of the linearisation of the PDE around its pulse solution is not periodic. However, this means that we cannot restrict ourselves to a fixed compact subset of the complex plane for all values of $h$ at the same time. In fact, it takes quite some effort to keep the part of the spectrum with large imaginary part under control.\par

It turns out to be convenient to partition our `half-strip' into four parts and to calculate the spectrum in
each part using different methods. Values close to $0$ are analyzed using the Fredholm properties of
$\mathcal{K}_{h;\overline{u}_h;c_h}$ exploiting many of the results from \S \ref{waveproperties}; values
with a large real part are considered using standard norm estimates, but values with a large imaginary
part are treated using a Fourier transform. The final set to consider is a compact set that is independent
of $h$ and bounded away from the origin. This allows us to apply a modified version of the procedure
described above that exploits the absence of spectrum in this region for the FitzHugh-Nagumo PDE.\par

Let us emphasize that our arguments here for bounded values of the spectral parameter $\lambda$
strongly use the fact that the PDE pulse is spectrally stable.
The main complication to establish the latter fact
is the presence of a secondary eigenvalue that is $O(\rho)$-close to the origin.
Intuitively, this eigenvalue arises as a consequence of the interaction between
the front and back solution to the Nagumo equation
that are both part of the singular pulse that arises in the $\rho \downarrow 0$ limit.
In the PDE case, Jones \cite{JONES1984} and Yanagida \cite{YAN1985}
essentially used shooting arguments to construct and analyze an Evans function $\mathcal{E}(\lambda)$
that vanishes precisely at eigenvalues. In particular,
they computed the sign of $\mathcal{E}'(0)$ and used counting arguments to show that the secondary eigenvalue
discussed above lies to the left of the origin. Currently, a program is underway to build a general
framework in this spirit based on the Maslov index
\cite{beck2018instability,chen2014stability,howard2016maslov}, which also works in
multi-dimensional spatial settings. In \cite{cornwell2017opening,cornwell2017existence}
this framework was applied to an equal-diffusion version
of the FitzHugh-Nagumo PDE. \par

An alternative approach involving Lin's method and exponential dichotomies was pioneered in
\cite{KRUSANSZM1997}. Based upon these ideas stability results have been obtained
for the LDE (\ref{ditishetalgemenealgemeneprobleem}) \cite{HJHSTBFHN}
and the PDE (\ref{firstfirstequation}) \cite{CARTER2016} in the non-hyperbolic regime
$r_0 \sim 0$. The first major advantage of this approach is that explicit bifurcation equations
can be formulated that allow asymptotic expansions to be developed for the location of the
interaction eigenvalue discussed above. The second major advantage is that
it allows us to avoid the use of the Evans function,
which cannot easily be defined in discrete settings
because MFDEs are ill-posed as initial value problems \cite{RUSHA1989}.
We believe that a direct approach along these lines should also be possible
for the infinite range system (\ref{infiniterangeversion}) as soon
as exponential dichotomies are available in this setting.

\subsection*{Nonlinear stability.}
The final step in our program is to leverage the spectral stability results to obtain a nonlinear stability result.
To do so, we follow \cite{HJHSTBFHN} and derive a formula that links the pointwise Green's function of our general problem
(\ref{infiniterangeversion}) to resolvents of the operator $\mathcal{K}_{h;\overline{u}_h;c_h}$ in \S \ref{sectiongreen}.
Since we have already analyzed the latter operator in detail, we readily obtain a spectral decomposition of this Green's
function into an explicit neutral part and a residual that decays exponentially in time and space. Therefore, we obtain detailed
estimates on the decay rate of the Green's function for the general problem. These Green's functions allow
us to establish   
the nonlinear stability of the family of
travelling pulse solutions $\overline{U}_h$. To be more precise, for each initial condition close to $\overline{U}_h(0)$,
we show that the solution with that initial condition converges at an exponential rate to the solution
$\overline{U}_h(\cdot+\tilde{\theta})$ for a small (and unique) phase shift $\tilde{\theta}$.\par


We emphasize that by now there are several techniques available to obtain nonlinear stability results
in the relatively simple spectral setting encountered in this paper. If a comparison principle is available,
which is not the case for the FitzHugh-Nagumo system,
one can follow the classic approach developed by Fife and McLeod \cite{Fife1977} to show that
travelling waves have a large basin of attraction. Indeed, one can construct explicit sub- and super-solutions
that trade-off additive perturbations at $t = 0$ to phase-shifts at $t = \infty$.
In fact, one can actually use this type of argument to establish the existence of travelling waves
by letting an appropriate initial condition evolve and tracking its asymptotic behaviour
\cite{CHEN1997, HJHNEGDIF}. For systems that can be written as gradient flows,
which is also not the case here, the existence and stability of travelling waves
can be obtained by using an elegant variational technique that was developed
by Gallay and Risler \cite{Gallay2007variational}.\par

In the spatially continuous setting, it is possible to freeze a travelling wave by passing
to a co-moving frame. In our setting, one can achieve this by simply
adding a convective term $-c_0 \partial_x (u,w)$ to the right hand side of
(\ref{firstfirstequation}). The main advantage is that
one can immediately use the semigroup $\mathrm{exp}[t \mathcal{L}_0]$
to describe the evolution of the linearised system in this co-moving frame,
which is temporally autonomous. Here $\mathcal{L}_0$ is the standard
linear operator associated to the linearisation of (\ref{firstfirstequation})
around $(\overline{u}_0, \overline{w}_0)$; see (\ref{eersteL0}).
For each $\vartheta \in \mathbb{R}$ one can subsequently construct the stable manifold
of $\big(\overline{u}_0( \cdot + \vartheta), \overline{w}_0(\cdot + \vartheta) \big)$
by applying a fixed point argument to Duhamel's formula. Upon varying $\vartheta$,
these stable manifolds span a tubular neighbourhood
of the family $(\overline{u}_0, \overline{w}_0)( \cdot + \mathbb{R})$.
This readily leads to the desired stability result; see e.g. \cite[{\S}4]{Kapitula2013spectral}.
We remark here that these stable manifolds are all related to each other via spatial shifts.\par

In the spatially discrete setting the wave can no longer be frozen. In particular,
the linearisation of (\ref{ditishetalgemenealgemeneprobleem}) around the pulse (\ref{eq:int:trvPulseAnsatz})
leads to an equation that is temporally shift-periodic.
In \cite{VL9} the authors attack this problem head-on by developing a shift-periodic version
of Floquet theory that leads to a nonlinear stability result in $\ell^\infty$.
However, they delicately exploit the geometric structure of $\ell^\infty$
and it is not clear how more degenerate spectral pictures can be fitted into the framework.
These issues are explained in detail in \cite[{\S}2]{HJHSTBFHN}.\par

In \cite{BGV2003} the authors found a way to express  the Green's function
of the temporally shift-periodic linear discrete equation in terms of resolvents
of the linear operator $\mathcal{L}_h$ associated to the pulse (\ref{eq:int:trvPulseAnsatz}).
Based on this procedure, it is possible to follow the spirit of the powerful pointwise Green's function
techniques pioneered by Zumbrun and Howard \cite{ZUMHOW1998}. Indeed, in \cite{HJHNLS}
a stability result is obtained in the setting of  discrete conservation laws, where one
encounters curves of essential spectrum that touch the imaginary axis. Using exponential dichotomies
in a setting with extended state-spaces $L^2([-h,h];\mathbb{R}^2) \times \mathbb{R}^2$,
pointwise $\lambda$-meromorphic expansions were obtained for the operators  $[\mathcal{L}_h - \lambda]^{-1}$.
This allowed the techniques from \cite{Beck2010nonlinear}
to be transferred from the continuous to the discrete setting.
A slightly more streamlined approach was developed in \cite{HJHSTBFHN}, which does not
need the extended state-space and avoids the use of a variation-of-constants formula.
However, exponential dichotomies are still used at certain key points.\par

In our paper we follow the spirit of the latter approach
and extend it to the present setting with infinite-range interactions.
In particular, we show how the use of exponential dichotomies
can be eliminated all together, which is a delicate task.
In addition, we need to be very careful in many computations
since integrals and sums over shifts as in (\ref{secondodequation})
can no longer be freely exchanged.
We emphasize here that our techniques do not depend on the specific LDE
that we are analyzing. All that is required is the spectral setting described above
and the fact that the shifts appearing in the problem are all rationally related.\par

Let us mention that it is also possible to bypass the construction of the stable manifolds
altogether and employ a direct phase-tracking approach along the lines of \cite{Zumbrun2009}.
In particular, one can couple the system with an extra equation for the phase.
To close the system, one chooses this extra equation in such a way
that the resulting nonlinear terms never encounter the non-decaying part of the relevant semigroup.
Such an approach
has been used in the current spectral setting to show that travelling waves
remain stable under the influence of a small stochastic noise term \cite{Hamster2017stability}.\par

\paragraph{Acknowledgements.}
Both authors acknowledge support from the Netherlands Organization for Scientific Research (NWO) (grant 639.032.612).

\section{Main results}\label{mainresults}

We consider the following system of equations
\begin{equation}\label{ditishetalgemeneprobleem}\begin{array}{lcl}\dot{u}_j&=&\frac{1}{h^2}\sum\limits_{k>0} \alpha_k[u_{j+k}+u_{j-k}-2u_j]+g(u_j)-w_j\\[0.2cm]
\dot{w}_j&=&\rho[u_j-\gamma w_j],\end{array}\end{equation}
which we refer to as the (spatially) discrete FitzHugh-Nagumo equation with infinite-range interactions. Often, for example in \cite{HJHFZHNGM,HJHSTBFHN}, it is assumed that only finitely many of these coefficients $\alpha_k$ are non-zero. However, we will impose the following much weaker conditions here.

\begin{assumption}{\text{H}\alpha 1}{}\label{$aannames$} The coefficients $\{\alpha_k\}_{k\in\Z_{>0}}$  satisfy the bound
\begin{equation}\begin{array}{lcl}\sum\limits_{k>0}|\alpha_k|k^2&<&\infty,\end{array}\end{equation}
as well as the identity
\begin{equation}\begin{array}{lcllcl}\sum\limits_{k>0}\alpha_k k^2&=&1.\end{array}\end{equation}
Finally, the inequality
\begin{equation}\label{conditionAz}\begin{array}{lcl}A(z):=\sum\limits_{k>0}\alpha_k\Big(1-\cos(kz)\Big)&>& 0\end{array}\end{equation}
holds for all $z\in(0,2\pi)$.\end{assumption}

We note that (\ref{conditionAz}) is automatically satisfied if $\alpha_1>0$ and $\alpha_k\geq 0$ for all $k\in\Z_{>1}$. The conditions in (\asref{aannames}{\text{H}\alpha}) ensure that for $\phi\in L^\infty(\R)$ with $\phi''\in L^2(\R)$, we have the limit
\begin{equation}\begin{array}{lcl}\lim\limits_{h\downarrow 0}\enskip\nrm{\frac{1}{h^2}\sum\limits_{k>0}\alpha_k\Big[\phi(\cdot+h k)+\phi(\cdot-hk)-2\phi(\cdot)\Big]-\phi''}_{L^2}=0,\end{array}\end{equation}
see Lemma \ref{eigenschappenDeltah}. In particular, we can see (\ref{ditishetalgemeneprobleem}) as the discretisation of the FitzHugh-Nagumo PDE (\ref{firstfirstequation}) on a grid with distance $h$. Additional remarks concerning the assumption (\asref{aannames}{\text{H}\alpha1}) can be found in \cite[\S 1]{BatesInfRange}.\par

Throughout this paper, we impose the following standard assumptions on the remaining parameters in (\ref{ditishetalgemeneprobleem}). The last condition on $\gamma$ in (\asref{aannamesconstanten}{\text{H}}) ensures that the origin is the only $j$-independent equilibrium of (\ref{ditishetalgemeneprobleem}).
\begin{assumption}{\text{HS}}{}\label{$aannamesconstanten$} The nonlinearity $g$ is given by $g(u)=u(1-u)(u-r_0)$, where $0<r_0<1$. In addition, we have $0<\rho<1$ and $0<\gamma<4(1-r_0)^{-2}$. \end{assumption}

Without explicitly mentioning it, we will allow all constants in this work to depend on $r_0,\rho$ and $\gamma$. Dependence on $h$ will always be mentioned explicitly. We will mainly work on the Sobolev spaces
\begin{equation}\begin{array}{lcl}H^1(\R)&=&\{f:\R\rightarrow\R|f,f'\in L^2(\R)\},\\[0.2cm]
H^2(\R)&=&\{f:\R\rightarrow\R|f,f',f''\in L^2(\R)\},\\[0.2cm]
\end{array}\end{equation}
with their standard norms
\begin{equation}\begin{array}{lcl}\nrm{f}_{H^1(\R)}&=&\left(\nrm{f}_{L^2(\R)}^2+\nrm{f'}_{L^2(\R)}^2\right)^\frac{1}{2},\\
\nrm{f}_{H^2(\R)}&=&\left(\nrm{f}_{L^2(\R)}^2+\nrm{f'}_{L^2(\R)}^2+\nrm{f''}_{L^2(\R)}^2\right)^\frac{1}{2}.
\end{array}\end{equation}\par

Our goal is to construct pulse solutions of (\ref{ditishetalgemeneprobleem}) as small perturbations to the travelling pulse solutions of the FitzHugh-Nagumo PDE. These latter pulses satisfy the system
\begin{equation}\label{firstequation}\begin{array}{lcl}c_0\overline{u}_0'&=&\overline{u}_0''+g(\overline{u}_0)-\overline{w}_0\\[0.2cm]
c_0\overline{w}_0'&=&\rho(\overline{u}_0-\gamma \overline{w}_0)\end{array}\end{equation}
with the boundary conditions
\begin{equation}\label{firstequationsboundary}\begin{array}{lcl}\lim\limits_{|\xi|\rightarrow\infty}(\overline{u}_0,\overline{w}_0)(\xi)&=&(0,0).\end{array}\end{equation}
If $(\overline{u}_0,\overline{w}_0)$ is a solution of (\ref{firstequation}) with wavespeed $c_0$, then the linearisation of (\ref{firstequation}) around this solution is characterized by the operator $\mathcal{L}_0:H^2(\R)\times H^1(\R)\rightarrow L^2(\R)\times L^2(\R)$, that acts as
\begin{equation}\label{eersteL0}\begin{array}{lcl}\mathcal{L}_0\left(\begin{array}{l}v\\ w\end{array}\right)&=&\left(\begin{array}{ll}c_0\frac{d}{d\xi}-\frac{d^2}{dx^2}-g_u(\overline{u}_0) & 1\\ -\rho & c_0\frac{d}{d\xi}+\gamma\rho\end{array}\right)\left(\begin{array}{l}v\\ w\end{array}\right).\end{array}\end{equation}
The existence of such pulse solutions for the case when $\rho$ is close to $0$ is established in \cite[\S 5.3]{JONES1995}. Here, we do not require $\rho>0$ to be small, but we simply impose the following condition.
\begin{assumption}{\text{HP}1}{}\label{$aannamespuls$} There exists a solution $(\overline{u}_0,\overline{w}_0)$ of (\ref{firstequation}) that satisfies the conditions (\ref{firstequationsboundary}) and has wavespeed $c_0\neq 0$. Furthermore, the operator $\mathcal{L}_0$ is Fredholm with index zero and it has a simple eigenvalue in zero.\end{assumption}

Recall that an eigenvalue $\lambda$ of a Fredholm operator $L$ is said to be \textit{simple} if the kernel of $L-\lambda$ is spanned by one vector $v$ and the equation $(L-\lambda)w=v$ does not have a solution $w$. Note that if $L$ has a formal adjoint $L^*$, this is equivalent to the condition that $\ip{v,w}\neq 0$ for all non-trivial $w\in \ker(L^*-\overline{\lambda})$.\par

We note that the conditions on $\mathcal{L}_0$ formulated in (\asref{aannamespuls}{\text{H}}) were established in \cite{JONES1984} for small $\rho>0$. In addition, these conditions imply that $\overline{u}_0'$ and $\overline{w}_0'$ decay exponentially.\par

In order to find travelling pulse solutions of (\ref{ditishetalgemeneprobleem}), we substitute the Ansatz
\begin{equation}\begin{array}{lcl}(u,w)_j(t)&=&(\overline{u}_h,\overline{w}_h)(hj+c_ht),\end{array}\end{equation}
into (\ref{ditishetalgemeneprobleem}) to obtain the system
\begin{equation}\label{ditishetprobleem}\begin{array}{lcl}c_h\overline{u}_h'(\xi)&=&\frac{1}{h^2}\sum\limits_{k>0}\alpha_k\Big[\overline{u}_h(\xi +h k)+\overline{u}_h(\xi -hk)-2\overline{u}_h(\xi)\Big]+g\big(\overline{u}_h(\xi)\big)-\overline{w}_h(\xi)\\[0.2cm]
c_h\overline{w}_h'(\xi)&=&\rho[\overline{u}_h(\xi)-\gamma \overline{w}_h(\xi)],\end{array}\end{equation}
in which $\xi=hj+c_ht$. The boundary conditions are given by
\begin{equation}\label{initialconditions}\begin{array}{lcl}\lim\limits_{|\xi|\rightarrow\infty}(\overline{u}_h,\overline{w}_h)(\xi)&=&(0,0).\end{array}\end{equation}
The existence of such solutions is established in our first main theorem.

\begin{theorem}[{see \S \ref{singularoperator}}]\label{Theorem1equivalent} Assume that (\asref{aannamespuls}{\text{H}}), (\asref{aannamesconstanten}{\text{H}}) and (\asref{aannames}{\text{H}\alpha1}) are satisfied. There exists a positive constant $h_*$ such that for all $h\in(0,h_*)$, the problem (\ref{ditishetprobleem}) with boundary conditions (\ref{initialconditions}) admits at least one solution $(c_h,\overline{u}_h,\overline{w}_h)$, which is locally unique in $\R\times H^1(\R)\times H^1(\R)$ up to translation and which has the property that
\begin{equation}\begin{array}{lcl}\lim\limits_{h\downarrow 0}\enskip(c_h-c_0,\overline{u}_h-\overline{u}_0,\overline{w}_h-\overline{w}_0)&=&(0,0,0)\qquad \text{in}\enskip\R\times H^1(\R)\times H^1(\R).\end{array}\end{equation}
\end{theorem}

Note that this result is very similar to \cite[Corollary 2.1]{Faye2015}. However, Faye and Scheel impose an extra assumption, similar to (\asref{extraaannames}{\text{H}\alpha}) below, which we do not need in our proof. This is a direct consequence of the strength of the method from \cite{BatesInfRange} that we described in \S \ref{introduction}.\par

Building on the existence of the travelling pulse solution, the natural next step is to show that our new-found pulse is asymptotically stable. However, we now do need to impose an extra condition on the coefficients $\{\alpha_k\}_{k>0}$.
\begin{assumption}{\text{H}\alpha 2}{}\label{$extraaannames$} The coefficients $\{\alpha_k\}_{k>0}$ satisfy the bound
\begin{equation}\begin{array}{lcl}\sum\limits_{k>0}|\alpha_k|e^{k\nu}&<&\infty\end{array}\end{equation}
for some $\nu>0$. \end{assumption}

Note that the prototype equation (\ref{infiniterangeversion}) indeed satisfies both assumptions (\asref{aannames}{\text{H}\alpha1}) and (\asref{extraaannames}{\text{H}\alpha}). An example of a system which satisfies (\asref{aannames}{\text{H}\alpha1}), but not (\asref{extraaannames}{\text{H}\alpha}) is given by
\begin{equation}\begin{array}{lcl}\label{infiniterangeversion2}\dot{u}_j&=&\frac{\kappa}{h^2}\sum\limits_{k>0}\frac{1}{k^4}[u_{j+k}+u_{j-k}-2u_j]+g(u_j)-w_j\\[0.2cm]
\dot{w}_j&=&\rho[u_j-\gamma w_j],\end{array}\end{equation}
in which $\kappa=\frac{6}{\pi^2}$ is the normalisation constant.\par

Moreover, we need to impose an extra condition on the operator $\mathcal{L}_0$ given by (\ref{eersteL0}). This spectral stability condition is established in \cite[Theorem 2]{Evans1972} together with \cite[Theorem 3.1]{YAN1985} for the case where $\rho$ is close to $0$.
\begin{assumption}{\text{HP}2}{}\label{$extraaannamespuls$} There exists a constant $\lambda_*>0$ such that for each $\lambda\in\C$ with $\Re \, \lambda\geq -\lambda_*$ and $\lambda\neq 0$, the operator
\begin{equation}\begin{array}{lcl}\mathcal{L}_0+\lambda:H^2(\R)\times H^1(\R)\rightarrow L^2(\R)\times L^2(\R)\end{array}\end{equation}
is invertible.\end{assumption}

To determine if the pulse solution described in Theorem \ref{Theorem1equivalent} is nonlinearly stable, we must first linearise (\ref{ditishetprobleem}) around this pulse and determine the spectral stability. The linearised operator now takes the form
\begin{equation}\label{definitieLh}\begin{array}{lcl}L_h\left(\begin{array}{l}v\\ w\end{array}\right)&=&\left(\begin{array}{ll}c_h\frac{d}{d\xi}-\Delta_h -g_u(\overline{u}_h)& 1\\ -\rho & c_h\frac{d}{d\xi}+\gamma\rho\end{array}\right)\left(\begin{array}{l}v\\ w\end{array}\right).\end{array}\end{equation}
Here the operator $\Delta_h$ is given by
\begin{equation}\begin{array}{lcl}\Delta_h\phi(\xi)&=&\frac{1}{h^2}\sum\limits_{k>0}\alpha_k\Big[\phi(\xi+h k)+\phi(\xi-hk)-2\phi(\xi)\Big].\end{array}\end{equation}
As usual, we define the \textit{spectrum}, $\sigma(L)$, of a bounded linear operator $L:H^1(\R)\times H^1(\R)\rightarrow L^2(\R)\times L^2(\R)$, as
\begin{equation}\begin{array}{lcl}\sigma(L)&=&\{\lambda\in\C : L-\lambda\text{ is not invertible}\}.\end{array}\end{equation}
Our second main theorem describes the spectrum of this operator $L_h$, or rather of $-L_h$, in a suitable half-plane.
\begin{theorem}[{see \S \ref{spectralanalysis}}]\label{totalespectrum}  Assume that (\asref{aannamespuls}{\text{H}\alpha 1}),(\asref{extraaannamespuls}{\text{H}}), (\asref{aannamesconstanten}{\text{H}}), (\asref{aannames}{\text{H}\alpha1}) and (\asref{extraaannames}{\text{H}\alpha}) are satisfied. There exist constants $\lambda_3>0$ and $h_{**}>0$ such that for all $h\in(0,h_{**})$, the spectrum of the operator $-L_h$ in the half-plane $\{z\in\C:\Re \,\text{ } z\geq-\lambda_3\}$ consists precisely of the points $k2\pi ic_h\frac{1}{h}$ for $k\in\Z$, which are all simple eigenvalues of $L_h$.\end{theorem}

We emphasize that $\lambda_3$ does not depend on $h$. The translational invariance of (\ref{ditishetprobleem}) guarantees that $\lambda=0$ is an eigenvalue of $-L_h$. In Lemma \ref{lemmaperiodiek} we show that the spectrum of the operator $L_h$ is periodic with period $2\pi ic_h\frac{1}{h}$, which means that the eigenvalues $k2\pi ic_h\frac{1}{h}$ for $k\in\Z$ all have the same properties as the zero eigenvalue.\par

Our final result concerns the nonlinear stability of our pulse solution, which we represent with the shorthand
\begin{equation}\label{definitieUh}\begin{array}{lcl}\Big[\overline{U}_h(t)\Big]_j&=&(\overline{u}_h,\overline{w}_h)(hj+c_ht).\end{array}\end{equation}
The perturbations are measured in the spaces $\ell^p$, which are defined by
\begin{equation}\begin{array}{lcl}\ell^p&=&\{V\in(\R^2)^\Z:\nrm{V}_{\ell^p}:=\Big[\sum\limits_{j\in\Z}|V_j|^p\Big]^{\frac{1}{p}}<\infty\}\end{array}\end{equation}
for $1\leq p<\infty$ and
\begin{equation}\begin{array}{lcl}\ell^\infty&=&\{V\in(\R^2)^\Z:||V||_{\ell^\infty}:=\sup\limits_{j\in\Z}|V_j|<\infty\}.\end{array}\end{equation}
\begin{theorem}[{see \S \ref{sectiongreen}}]\label{nonlinearstability}  Assume that (\asref{aannamespuls}{\text{H}}),(\asref{extraaannamespuls}{\text{H}}), (\asref{aannamesconstanten}{\text{H}}), (\asref{aannames}{\text{H}\alpha1}) and (\asref{extraaannames}{\text{H}\alpha}) are satisfied. Fix $0<h\leq h_{**}$ and $1\leq p\leq \infty$. Then there exist constants $\delta>0$, $C>0$ and $\beta>0$, which may depend on $h$ but not on $p$, such that for all initial conditions $U^0\in\ell^p$ with $\nrm{U^0-\overline{U}_h(0)}_{\ell^p}<\delta$, there exists an asymptotic phase shift $\tilde{\theta}\in\R$ such that the solution $U=(u,w)$ of (\ref{ditishetalgemeneprobleem}) with $U(0)=U^0$ satisfies the bound
\begin{equation}\begin{array}{lcl}\nrm{U(t)-\overline{U}_h(t +\tilde{\theta})}_{\ell^p}&\leq& Ce^{-\beta t}\nrm{U^0-\overline{U}_h(0)}_{\ell^p}\end{array}\end{equation}
for all $t>0$.\end{theorem}


\section{The singular perturbation}\label{singularoperator}

The main difficulty in analysing the travelling wave MFDE (\ref{ditishetprobleem}) is that it is a singular perturbation of the ODE (\ref{firstequation}). Indeed, the second derivative in (\ref{firstequation}) is replaced by the linear operator $\Delta_h:H^1(\R)\rightarrow L^2(\R)$ that acts as
\begin{equation}\begin{array}{lcl}\Delta_h\phi(\xi)&=&\frac{1}{h^2}\sum\limits_{k>0}\alpha_k\Big(\phi(\xi+h k)+\phi(\xi-hk)-2\phi(\xi)\Big).\end{array}\end{equation}
We will see in Lemma \ref{eigenschappenDeltah} that for all $\phi\in L^\infty(\R)$ with $\phi''\in L^2(\R)$, we have that $\lim\limits_{h\downarrow 0}\enskip\nrm{\Delta_h\phi-\phi''}_{L^2}=0$. Hence, the bounded operator $\Delta_h$ converges pointwise on a dense subset of $H^1(\R)$ to an unbounded operator on that same dense subset. In particular, the norm of the operator $\Delta_h$ grows to infinity as $h\downarrow 0$.\par

Since there are no second derivatives involved in (\ref{ditishetprobleem}), we have to view it as an equation posed on the space $H^1(\R)\times H^1(\R)$, while the ODE (\ref{firstequation}) is posed on the space $H^2(\R)\times H^1(\R)$. From now on we write
\begin{equation}\begin{array}{lcl}\textbf{H}^1&:=& H^1(\R)\times H^1(\R),\\[0.2cm]
\textbf{L}^2&:=&L^2(\R)\times L^2(\R).\end{array}\end{equation}

The main results in this section will be used in several different settings. In order to accommodate this, we introduce the following conditions.

\begin{assumption}{\text{hFam}}{}\label{$familyassumption$} For each $h>0$ there is a pair $(\tilde{u}_h,\tilde{w}_h)\in \textbf{H}^1$ and a constant $\tilde{c}_h$ such that $(\tilde{u}_h,\tilde{w}_h)-(\overline{u}_0,\overline{w}_0)\rightarrow 0$ in $\textbf{H}^1$ and $\tilde{c}_h\rightarrow c_0$ as $h\downarrow 0$.\end{assumption}

In the proof of Theorem \ref{Theorem1equivalent} we choose $(\tilde{u}_h,\tilde{w}_h)$ and $\tilde{c}_0$ to be $(\overline{u}_0,\overline{w}_0)$ and $c_0$ for all values of $h$. However, in \S \ref{waveproperties} we let $(\tilde{u}_h,\tilde{w}_h)$ be the travelling pulse $(\overline{u}_h,\overline{w}_h)$ from Theorem \ref{Theorem1equivalent} and we let $\tilde{c}_h$ be its wave speed $c_h$. \par

If (\asref{familyassumption}{\text{h}}) is satisfied, then for $\delta>0$ and $h>0$ we define the operators
\begin{equation}\label{overlinemathcallplus}\begin{array}{lcl}\overline{\mathcal{L}}^+_{h,\delta}&=&\left(\begin{array}{ll}\tilde{c}_h\frac{d}{d\xi}-\Delta_h-g_u(\tilde{u}_h)+\delta & 1\\ -\rho & \tilde{c}_h\frac{d}{d\xi}+\gamma\rho+\delta\end{array}\right)\end{array}\end{equation}
and
\begin{equation}\label{overlinemathcallmin}\begin{array}{lcl}\overline{\mathcal{L}}^-_{h,\delta}&=&\left(\begin{array}{ll}-\tilde{c}_h\frac{d}{d\xi}-\Delta_h-g_u(\tilde{u}_h)+\delta & -\rho\\ 1 & -\tilde{c}_h\frac{d}{d\xi}+\gamma\rho+\delta\end{array}\right).\end{array}\end{equation}
These operators are bounded linear functions from $\textbf{H}^1$ to $\textbf{L}^2$. We see that $\overline{\mathcal{L}}_{h,\delta}^-$ is the adjoint operator of $\overline{\mathcal{L}}_{h,\delta}^+$, in the sense that
\begin{equation}\begin{array}{lcl}\ip{(\phi,\psi),\overline{\mathcal{L}}_{h,\delta}^+(\theta,\chi)}&=&\ip{\overline{\mathcal{L}}_{h,\delta}^-(\phi,\psi),(\theta,\chi)}\end{array}\end{equation}
holds for all  $(\phi,\psi),(\theta,\chi)\in \textbf{L}^2$. Here we have introduced the notation
\begin{equation}\begin{array}{lcl}\ip{(\phi,\psi),(\theta,\chi)}&=&\ip{\phi,\theta}+\ip{\psi,\chi}\\[0.2cm]
&=&\int\limits_{-\infty}^\infty \Big(\phi(x)\theta(x)+\psi(x)\chi(x)\Big)dx\end{array}\end{equation}
for $(\phi,\psi),(\theta,\chi)\in \textbf{L}^2$.\par

Since at some point we want to consider complex-valued functions, we also work in the spaces $H^2_\C(\R)$, $H^1_\C(\R)$ and $L^2_\C(\R)$, which are given by
\begin{equation}
\begin{array}{lcl}
H^2_\C(\R)&=&\{f+gi|f,g\in H^2(\R)\},\\[0.2cm]
H^1_\C(\R)&=&\{f+gi|f,g\in H^1(\R)\},\\[0.2cm]
L^2_\C(\R)&=&\{f+gi|f,g\in L^2(\R)\}.
\end{array}
\end{equation}
These spaces are equipped with the inner product
\begin{equation}\begin{array}{lcl}\ip{\phi,\psi}&=&\int \Big(f_1(x)+ig_1(x)\Big)\Big(f_2(x)-ig_2(x)\Big)dx\end{array}\end{equation}
for $\phi=f_1+ig_1,\psi=f_2+ig_2$. As before, we write
\begin{equation}\begin{array}{lcl}\textbf{H}^1_\C&=&H^1_\C(\R)\times H^1_\C(\R)\\[0.2cm]
\textbf{L}^2_\C&=&L^2_\C(\R)\times L^2_\C(\R).\end{array}\end{equation}
Each operator $L$ from $\textbf{H}^1$ to $\textbf{L}^2$ can be extended to an operator from $\textbf{H}_\C^1$ to $\textbf{L}_\C^2$ by writing
\begin{equation}\begin{array}{lcl}L(f+ig)=Lf+iLg.\end{array}\end{equation}
It is well-known that this complexification preserves adjoints, invertibility, inverses, injectivity, surjectivity and boundedness, see for example \cite{SAB2007}. If $\lambda\in\C$ then the operators $\overline{\mathcal{L}}^\pm_{h,\lambda}$ are defined analogously to their real counterparts, but now we view them as operators from $H_\C^1(\R)\times H_\C^1(\R)$ to $L^2_\C(\R)\times L^2_\C(\R)$. Whenever it is clear that we are working in the complex setting we drop the subscript $\C$ from the spaces $\textbf{H}^1_\C$ and $\textbf{L}^2_\C$ and simply write $\textbf{H}^1$ and $\textbf{L}^2$. \par

We also introduce the operators $\mathcal{L}_0^\pm:H^2(\R)\times H^1(\R)\rightarrow L^2(\R)\times L^2(\R)$, that act as
\begin{equation}\begin{array}{lcl}\mathcal{L}_0^+&=&\left(\begin{array}{ll}c_0\frac{d}{d\xi}-\frac{d^2}{dx^2}-g_u(\overline{u}_0) & 1\\ -\rho & c_0\frac{d}{d\xi}+\gamma\rho\end{array}\right)\end{array}\end{equation}
and
\begin{equation}\begin{array}{lcl}\mathcal{L}_{0}^-&=&\left(\begin{array}{ll}-c_0\frac{d}{d\xi}-\frac{d^2}{dx^2}-g_u(\overline{u}_0) & -\rho\\ 1 & -c_0\frac{d}{d\xi}+\gamma\rho\end{array}\right).\end{array}\end{equation}
These operators can be viewed as the formal $h\downarrow 0$ limits of the operators $\overline{\mathcal{L}}_{h,0}^\pm$.
Upon introducing the notation
\begin{equation}\begin{array}{lcl}(\phi_0^+,\psi_0^+)&=&\frac{(\overline{u}_0',\overline{w}_0')}{\nrm{(\overline{u}_0',\overline{w}_0')}_{\textbf{L}^2}},\end{array}\end{equation}
we see that $\mathcal{L}_0^+(\phi_0^+,\psi_0^+)=0$ by differentiating (\ref{firstequation}).\par

To set the stage, we summarize several basic properties of $\mathcal{L}_0^\pm$. The last property references
a spectral set $M$, on which we impose the following
condition.
\begin{assumption}{\text{hM}}{}\label{$Massumption$} The set $M\subset \C$ is compact with $0\notin M$. In addition, recalling the constant $\lambda_*$ appearing in (\asref{extraaannamespuls}{\text{H}}), we have $\Re \, z\geq
-\lambda_*$ for all $z\in M$.\end{assumption}
In \S\ref{spectralanalysis} the set $M$ will be fixed
as the final region of our spectral analysis, which we
will refer to as $R_4$.
The proof of this result follows the standard procedure described in \cite[Lemma 5]{BatesInfRange} and, as such, will be omitted.

\begin{lemma}\label{eigenschappenL0} Assume that (\asref{aannamespuls}{\text{H}}), (\asref{aannamesconstanten}{\text{H}}) and (\asref{aannames}{\text{H}\alpha1}) are satisfied.
Then the following results hold.
\begin{enumerate}
\item We have that $(\phi_0^+,\psi_0^+)\in H^2(\R)\times H^1(\R)$ and $\ker(\mathcal{L}_0^+)=\span\{(\phi_0^+,\psi_0^+)\}$.
\item There exist $(\phi_0^-,\psi_0^-)\in H^2(\R)\times H^1(\R)$ with $\nrm{(\phi_0^-,\psi_0^-)}_{\textbf{L}^2}=1$, with
$\ip{(\overline{u}_0',\overline{w}_0'),(\phi_0^-,\psi_0^-)}>0$ and $\ker(\mathcal{L}_0^-)=\span\{(\phi_0^-,\psi_0^-)\}$.
\item For every $(\theta,\chi)\in \textbf{L}^2$ the problem $\mathcal{L}_0^{\pm}(\phi,\psi)=(\theta,\chi)$ with $(\phi,\psi)\in H^2(\R)\times H^1(\R)$ and $\ip{(\phi,\psi),(\phi_0^{\pm},\psi_0^{\pm})}=0$ has a unique solution $(\phi,\psi)$ if and only if $\ip{(\theta,\chi),(\phi_0^{\mp},\psi_0^{\mp})}=0$.
\item There exists a positive constant $C_1$ such that
\begin{equation}\begin{array}{lcl}\nrm{(\phi,\psi)}_{H^2(\R)\times H^1(\R)}&\leq& C_1\nrm{\mathcal{L}_0^{\pm}(\phi,\psi)}_{\textbf{L}^2}\end{array}\end{equation}
for all $(\phi,\psi)\in H^2(\R)\times H^1(\R)$ with $\ip{(\phi,\psi), (\phi_0^{\pm},\psi_0^{\pm})}=0$.
\item There exists a positive constant $C_2$ and a small constant $\delta_0>0$ such that for all $0<\delta<\delta_0$ we have
\begin{equation}\begin{array}{lcl}\nrm{(\mathcal{L}_0^\pm+\delta)^{-1}(\theta,\chi)}_{H^2(\R)\times H^1(\R)}\leq C_2\Big[\nrm{(\theta,\chi)}_{L^2(\R)\times L^2(\R)}+\frac{1}{\delta}|\ip{(\theta,\chi),(\phi_0^{\mp},\psi_0^{\mp})}|\Big]\end{array}\end{equation}
for all $(\theta,\chi)\in L^2(\R)\times L^2(\R)$.
\item If (\asref{extraaannamespuls}{\text{H}}) is also satisfied, then for each $M\subset \C$ that satisfies (\asref{Massumption}{\text{h}}), there exists a constant $C_3>0$ such that the uniform bound
\begin{equation}\begin{array}{lcl}\nrm{(\mathcal{L}_0^\pm+\lambda)^{-1}(\theta,\chi)}_{H^2_\C(\R)\times H^1_\C(\R)}\leq C_3\nrm{(\theta,\chi)}_{L^2_\C(\R)\times L^2_\C(\R)}\end{array}\end{equation}
holds for all $(\theta,\chi)\in L^2_\C(\R)\times L^2_\C(\R)$ and all $\lambda\in M$.
\end{enumerate}
\end{lemma}

The main goal of this section is to prove the following two propositions, which transfer parts (5) and (6) of Lemma \ref{eigenschappenL0} to the discrete setting.
\begin{proposition}\label{equivalenttheorem4} Assume that (\asref{familyassumption}{h}), (\asref{aannamespuls}{\text{H}}), (\asref{aannamesconstanten}{\text{H}}) and (\asref{aannames}{\text{H}\alpha1}) are satisfied. There exists a positive constant $C_0'$ and a positive function $h_0'(\cdot):\R^+\rightarrow \R^+$, depending only on the choice of $(\tilde{u}_h,\tilde{w}_h)$ and $\tilde{c_h}$, such that for every $0<\delta<\delta_0$ and every $h\in(0,h_0'(\delta))$, the operators $\overline{\mathcal{L}}_{h,\delta}^{\pm}$ are homeomorphisms from $\textbf{H}^1$ to $\textbf{L}^2$ that satisfy the bounds
\begin{equation}\begin{array}{lcl}\nrm{(\overline{\mathcal{L}}_{h,\delta}^{\pm})^{-1}(\theta,\chi)}_{\textbf{H}^1}&\leq& C_0'\Big[\nrm{(\theta,\chi)}_{\textbf{L}^2}+\frac{1}{\delta}|\ip{(\theta,\chi),(\phi_0^{\mp},\psi_0^{\mp})}|\Big]\end{array}\end{equation}
for all $(\theta,\chi)\in \textbf{L}^2$.
\end{proposition}

\begin{proposition}\label{equivalenttheorem4version2} Assume that (\asref{familyassumption}{h}), (\asref{aannamespuls}{\text{H}}),(\asref{extraaannamespuls}{\text{H}}), (\asref{aannamesconstanten}{\text{H}}) and (\asref{aannames}{\text{H}\alpha1}) are satisfied. Let $M\subset \C$ satisfy (\asref{Massumption}{\text{h}}). Then there exists a constant $h_M>0$, depending only on $M$ and the choice of $(\tilde{u}_h,\tilde{w}_h)$ and $\tilde{c_h}$, such that for all $0<h\leq h_M$ and all $\lambda\in M$ the operator $\overline{\mathcal{L}}_{h,\lambda}^{\pm}$ is a homeomorphism from $\textbf{H}^1$ to $\textbf{L}^2$.
\end{proposition}

\subsection{Strategy}
Our techniques here are inspired strongly by the approach developed in  \cite[\S 2-4]{BatesInfRange}. Indeed, Proposition \ref{equivalenttheorem4} and
Proposition \ref{lemma6equivalent} are the equivalents of \cite[Theorem 4]{BatesInfRange} and \cite[Lemma 6]{BatesInfRange} respectively.
The difference between our results and those in \cite{BatesInfRange} is that Bates, Chen and Chmaj study the discrete Nagumo equation, which can be seen as
the one-dimensional fast component of the FitzHugh-Nagumo equation by setting $\rho=0$ in (\ref{ditishetalgemeneprobleem}). In addition, the results
in \cite{BatesInfRange} are restricted to $\lambda\in\R$, while we allow $\lambda\in\C$ in Proposition \ref{equivalenttheorem4version2}.
These differences play a crucial role in the proof of Lemma \ref{lemma6bewijs3} below.
\par

Recall the constant $\delta_0>0$ appearing in Lemma \ref{eigenschappenL0}. For $0<\delta<\delta_0$ and $h>0$ we define the quantities
\begin{equation}\label{eqoverlinelambda}\begin{array}{lcl}\overline{\Lambda}^{\pm}(h,\delta)&=&\inf\limits_{\nrm{(\phi,\psi)}_{\textbf{H}^1}=1}\left[\nrm{\overline{\mathcal{L}}^{\pm}_{h,\delta}(\phi,\psi)}_{\textbf{L}^2}+\frac{1}{\delta}\left|\ip{\overline{\mathcal{L}}_{h,\delta}^{\pm}(\phi,\psi),(\phi_0^{\mp},\psi_0^{\mp})}\right|\right],\end{array}\end{equation}
together with
\begin{equation}\begin{array}{lcl}\label{eqoverlinelambda2}\overline{\Lambda}^{\pm}(\delta)&=&\liminf\limits_{h\downarrow 0}\overline{\Lambda}^{\pm}(h,\delta).\end{array}\end{equation}
Similarly for $M\subset \C$ that satisfies (\asref{Massumption}{\text{h}}) and $h>0$ we define
\begin{equation}\begin{array}{lcl}\overline{\Lambda}^{\pm}(h,M)&=&\inf\limits_{\nrm{(\phi,\psi)}_{\textbf{H}^1}=1,\ \lambda\in M}\left[\nrm{\overline{\mathcal{L}}^{\pm}_{h,\lambda}(\phi,\psi)}_{\textbf{L}^2}\right],\end{array}\end{equation}
together with
\begin{equation}\label{eqoverlinelambdam}\begin{array}{lcl}\overline{\Lambda}^{\pm}(M)&=&\liminf\limits_{h\downarrow 0}\overline{\Lambda}^{\pm}(h,M).\end{array}\end{equation}
The key ingredients that we need to establish Propositions \ref{equivalenttheorem4} and \ref{equivalenttheorem4version2} are lower bounds on the quantities $\overline{\Lambda}^\pm(\delta)$ and $\overline{\Lambda}^\pm(M)$. These are provided in the result below, which we consider to be the technical heart of this section.
\begin{proposition}\label{lemma6equivalent} Assume that (\asref{familyassumption}{h}), (\asref{aannamespuls}{\text{H}}), (\asref{aannamesconstanten}{\text{H}}) and (\asref{aannames}{\text{H}\alpha1}) are satisfied. There exists a positive constant $C_0'$, depending only on our choice of $(\tilde{u}_h,\tilde{w}_h)$ and $\tilde{c_h}$, such that $\overline{\Lambda}^{\pm}(\delta)>\frac{2}{C_0'}$ for all $0<\delta<\delta_0$. Similarly if $M\subset \C$ satisfies (\asref{Massumption}{\text{h}}), then there exists a positive constant $C_M'$, depending only on $M$ and our choice of $(\tilde{u}_h,\tilde{w}_h)$ and $\tilde{c_h}$, such that $\overline{\Lambda}^{\pm}(M)>\frac{2}{C_M'}$.
\end{proposition}
\noindent\textit{Proof of Proposition \ref{equivalenttheorem4}.}
Let $\delta>0$ be fixed. Since $\overline{\Lambda}^{\pm}(\delta)\geq \frac{2}{C_0'}$, the definition (\ref{eqoverlinelambda2}) implies that there exists $h_0'(\delta)$ such that $\overline{\Lambda}(h,\delta)\geq\frac{1}{C_0'}$ for all $h\in(0,h_0'(\delta)]$. Now pick $h\in(0,h_0'(\delta)]$.\par
First of all, $\overline{\mathcal{L}}_{h,\delta}^{\pm}$ is a bounded operator from $\textbf{H}^1$ to $\textbf{L}^2$. Since $\overline{\Lambda}^{\pm}(h,\delta)$ is strictly positive, this implies that $\overline{\mathcal{L}}_{h,\delta}^{\pm}$ is a homeomorphism from $\textbf{H}^1$ to its image $\overline{\mathcal{L}}_{h,\delta}^{\pm}(\textbf{H}^1)$. Furthermore, the norm of the inverse $(\overline{\mathcal{L}}_{h,\delta}^{\pm})^{-1}$ from $\overline{\mathcal{L}}_{h,\delta}^{\pm}(\textbf{H}^1)\subset \textbf{L}^2$ is bounded by $\frac{1}{\overline{\Lambda}^{\pm}(h,\delta)}\leq C_0'$. Since $\overline{\mathcal{L}}_{h,\delta}^\pm$ is bounded, it follows that $\overline{\mathcal{L}}_{h,\delta}^{\pm}(\textbf{H}^1)$ is closed in $\textbf{L}^2$.\par

For the remainder of this proof, we only consider the operators $\overline{\mathcal{L}}_{h,\delta}^+$, noting that their counterparts $\overline{\mathcal{L}}_{h,\delta}^-$ can be treated in an identical fashion.\par

Seeking a contradiction, let us assume that $\overline{\mathcal{L}}_{h,\delta}^{+}(\textbf{H}^1)\neq \textbf{L}^2$, which implies that there exists a non-zero $(\theta,\chi)\in \textbf{L}^2$ orthogonal to $\overline{\mathcal{L}}_{h,\delta}^{+}(\textbf{H}^1)$. For any $\phi\in C_c^\infty(\R)$, we hence obtain
\begin{equation}\begin{array}{lcl}0&=&\ip{\overline{\mathcal{L}}_{h,\delta}^+(\phi,0),(\theta,\chi)}\\[0.2cm]
&=&\ip{\tilde{c}_h\phi'-\Delta_h\phi-g_u(\tilde{u}_h)\phi+\delta\phi,\theta}+\ip{-\rho\phi,\chi}\\[0.2cm]
&=&\tilde{c}_h\ip{\phi',\theta}+\ip{\phi,-\Delta_h\theta-g_u(\tilde{u}_h)\theta+\delta\theta-\rho\chi}.\end{array}\end{equation}
By definition this implies that $\theta$ has a weak derivative and that $\tilde{c}_h\theta'=-\Delta_h \theta-g_u(\tilde{u}_h)\theta+\delta\theta-\rho\chi\in L^2(\R)$. In particular, we see that $\theta\in H^1(\R)$.\par

For any $\psi\in C_c^\infty(\R)$ a similar computation yields
\begin{equation}\begin{array}{lcl}0&=&\ip{\overline{\mathcal{L}}_{h,\delta}^+(0,\psi),(\theta,\chi)}\\[0.2cm]
&=&\ip{\psi,\theta}+\ip{\tilde{c_h}\psi'+(\gamma\rho+\delta)\psi,\chi}\\[0.2cm]
&=&\tilde{c}_h\ip{\psi',\chi}+\ip{\psi,\theta+(\gamma\rho+\delta)\chi}.\end{array}\end{equation}
Again, this means that $\chi$ has a weak derivative and in fact $\tilde{c}_h\chi'=\theta+(\gamma\rho+\delta)\chi$. In particular, it follows that $\chi\in H^1(\R)$.\par

We, therefore, conclude that
\begin{equation}\begin{array}{lcl}0&=&\ip{\overline{\mathcal{L}}_{h,\delta}^{+}(\phi,\psi),(\theta,\chi)}\\[0.2cm]
&=&\ip{(\phi,\psi),(\overline{\mathcal{L}}_{h,\delta}^{-}(\theta,\chi)}\end{array}\end{equation}
holds for all $(\phi,\psi)\in \textbf{H}^1$. Since $\textbf{H}^1$ is dense in $\textbf{L}^2$ this implies that $\overline{\mathcal{L}}_{h,\delta}^{-}(\theta,\chi)=0$. Since we already know that $\overline{\mathcal{L}}_{h,\delta}^{-}$ is injective, this means that $(\theta,\chi)=0$, which gives a contradiction. Hence, we must have $\overline{\mathcal{L}}_{h,\delta}^{+}(\textbf{H}^1)=\textbf{L}^2$, as desired.
\qed

\vspace*{4pt}\noindent\textit{Proof of Proposition \ref{equivalenttheorem4version2}.} The result follows in the same
fashion as outlined in the proof of Proposition \ref{equivalenttheorem4} above.\qed


\vspace*{4pt}\noindent\textit{Proof of Theorem \ref{Theorem1equivalent}.}
Recall the notation $(c_0, \overline{u}_0,\overline{w}_0)$ for the pulse solution of the PDE
(\ref{firstfirstequation}).
The desired solutions to (\ref{ditishetprobleem}) can be constructed by writing
\begin{equation}\label{equivalentprobleem}\begin{array}{lcl}(\overline{u}_h,\overline{w}_h)&=&(\overline{u}_0,\overline{w}_0)+(\phi_h,\psi_h)\end{array}\end{equation}
and setting up a fixed-point argument to find the small perturbations $(\phi_h,\psi_h)$
along with the wavespeed $c_h$. This can be done by following the standard procedure
described in \cite[{\S}4]{BatesInfRange}.
The relevant linear operator can be found by
applying Proposition \ref{equivalenttheorem4} to the constant family
$(\tilde{u}_h,\tilde{w}_h)=(\overline{u}_0,\overline{w}_0)$ and $\tilde{c}_h=c_0$. \qed

\subsection{Preliminaries}
\label{sec:sp:conv:prlm}

Our goal here is to establish some basic facts concerning the operator $\Delta_h$. In particular, we extend the real-valued results from \cite{BatesInfRange} to complex-valued functions. We emphasize that the inequalities in Lemma \ref{eigenschappencomplexeDeltah} in general do not hold for the imaginary parts of these inner products.
\begin{lemma}\label{eigenschappenDeltah} (see \cite[Lemma 3]{BatesInfRange}) Assume that (\asref{aannames}{\text{H}\alpha1}) is satisfied. The following three properties hold.
\begin{enumerate}
\item For all $\phi\in L^\infty(\R)$ with $\phi''\in L^2(\R)$ we have $\lim\limits_{h\downarrow 0}\enskip\nrm{\Delta_h\phi-\phi''}_{L^2}=0$.
\item For all $\phi\in H^1(\R)$ and $h>0$ we have $\ip{\Delta_h\phi,\phi'}=0$.
\item For all $\phi,\psi\in L^2(\R)$ and $h>0$ we have $\ip{\Delta_h\phi,\psi}=\ip{\phi,\Delta_h\psi}$ and $\ip{\Delta_h\phi,\phi}\leq 0$.
\end{enumerate}\end{lemma}
\begin{lemma}\label{eigenschappencomplexeDeltah} Assume that (\asref{aannames}{\text{H}\alpha1}) is satisfied and pick $f\in H^1_\C(\R)$. Then the following properties hold.
\begin{enumerate}
\item For all $h>0$ we have $\Re \,\ip{-\Delta_hf,f}\geq 0$.
\item For all $h>0$ we have $\Re \,\ip{\Delta_h f,f'}=0$.
\item We have $\Re \,\ip{f,f'}=0$.
\item For all $\lambda\in\C$ we have $\Re \,\ip{\lambda f,f'}=2\ (\Im \,\lambda)\ip{\Re \, f,\Im \, f'}$.
\end{enumerate}
\end{lemma}
\vspace*{4pt}\noindent\textit{Proof.} Write $f=\phi+i\psi$ with $\phi,\psi\in H^1(\R)$. Lemma \ref{eigenschappenDeltah} implies that
\begin{equation}\begin{array}{lcl} \Re \,\ip{-\Delta_hf,f}&=&\Re \,\int \Big(-\Delta_h\phi -i\Delta_h \psi\Big)(x)\Big(\phi-i\psi\Big)(x)dx\\[0.2cm]
&=&\int (-\Delta_h\phi)(x)\phi(x)+(-\Delta_h\psi)(x)\psi(x)dx\\[0.2cm]
&=&\ip{-\Delta_h\phi,\phi}+\ip{-\Delta_h\psi,\psi}\\[0.2cm]
&\geq &0.\end{array}\end{equation}
Similarly we have
\begin{equation}\begin{array}{lcl} \Re \,\ip{\Delta_h f,f'}&=&\ip{-\Delta_h\phi,\phi'}+\ip{-\Delta_h\psi,\psi'}\\[0.2cm]
&= &0.\end{array}\end{equation}
For $\lambda\in\C$ we may compute
\begin{equation}\begin{array}{lcl} \Re \,\ip{\lambda f,f'}&=&\Re \,\int \Big(\lambda \phi(x)+\lambda i\psi(x)\Big)\Big(\phi'(x)-i\psi'(x)\Big)dx\\[0.2cm]
&=&( \Re \,\lambda)\ip{\phi,\phi'}+(\Im \,\lambda)\ip{\phi,\psi'}-(\Im \,\lambda) \ip{\psi,\phi'}+(\Re \,\lambda)\ip{\psi,\psi'}\\[0.2cm]
&=&0+2\ (\Im \,\lambda)\ip{\phi,\psi'}+0\\[0.2cm]
&=&2\ (\Im \,\lambda)\ip{\phi,\psi'}.\end{array}\end{equation}
Taking $\lambda=1$ gives the third property.\qed

\subsection{Proof of Proposition \ref{lemma6equivalent}}
We now set out to prove Proposition \ref{lemma6equivalent}. In Lemmas \ref{lemma6bewijs1} and \ref{lemma6bewijs1.5} we construct weakly converging sequences that realize the infima in (\ref{eqoverlinelambda})-(\ref{eqoverlinelambdam}). In Lemmas \ref{lemma6bewijs2}-\ref{lemma6bewijs4} we exploit the structure of our operators (\ref{overlinemathcallplus}) and (\ref{overlinemathcallmin}) to recover bounds on the derivatives of these sequences that are typically lost when taking weak limits. Recall the constant $C_2>0$ defined in Lemma \ref{eigenschappenL0}, which does not depend on $\delta>0$.
\begin{lemma}\label{lemma6bewijs1} Assume that (\asref{familyassumption}{h}), (\asref{aannamespuls}{\text{H}}), (\asref{aannamesconstanten}{\text{H}}) and (\asref{aannames}{\text{H}\alpha1}) are satisfied. Consider the setting of Proposition \ref{lemma6equivalent} and fix $0<\delta<\delta_0$. Then there exists a sequence $\{(h_j,\phi_j,\psi_j)\}_{j\geq 0}$ in $(0,1)\times \textbf{H}^1$ with the following properties.
\begin{enumerate}
\item We have $\lim_{j\rightarrow\infty}h_j=0$ and $ \nrm{(\phi_j,\psi_j)}_{\textbf{H}^1}=1$ for all $j\geq 0$.
\item The sequence $(\theta_j,\chi_j)=\overline{\mathcal{L}}_{h_j,\delta}^+(\phi_j,\psi_j)$ satisfies
\begin{equation}\begin{array}{lcl}\lim_{j\rightarrow\infty} \Big[\nrm{(\theta_j,\chi_j)}_{\textbf{L}^2}+\frac{1}{\delta}|\ip{(\theta_j,\chi_j),(\phi_0^{-},\psi_0^{-})}|\Big]&=&\overline{\Lambda}^+(\delta).\end{array}\end{equation}
\item There exist $(\phi,\psi)\in \textbf{H}^1$ and $(\theta,\chi)\in \textbf{L}^2$ such that $(\phi_j,\psi_j)\rightharpoonup (\phi,\psi)$ weakly in $\textbf{H}^1$ and such that $(\theta_j,\chi_j)\rightharpoonup (\theta,\chi)$ weakly in $\textbf{L}^2$ as $j\rightarrow\infty$.
\item We have $(\phi_j,\psi_j)\rightarrow (\phi,\psi)$ in $L^2_{\mathrm{loc}}(\R)\times L^2_{\mathrm{loc}}(\R)$ as $j\rightarrow \infty$.
\item The pair $(\phi,\psi)$ is a weak solution to $(\overline{\mathcal{L}}_0^{+}+\delta)(\phi,\psi)=(\theta,\chi)$.
\item We have the bound \begin{equation}\begin{array}{lcl}\nrm{(\phi,\psi)}_{H^2(\R)\times H^1(\R)}&\leq &C_2\overline{\Lambda}^{+}(\delta).\end{array}\end{equation}
\end{enumerate}
The same statements hold upon replacing $\overline{\mathcal{L}}_{h,\delta}^+$, $\overline{\Lambda}^+$ and $\overline{\mathcal{L}}_{0}^+$ by $\overline{\mathcal{L}}_{h,\delta}^-$, $\overline{\Lambda}^-$ and $\overline{\mathcal{L}}_{0}^-$.
\end{lemma}

\vspace*{4pt}\noindent\textit{Proof.} Let $0<\delta<\delta_0$ be fixed. By definition of $\overline{\Lambda}^{+}(\delta)$ there exists a sequence $\{(h_j,\phi_j,\psi_j)\}$ in $(0,1)\times \textbf{H}^1$ such that (1) and (2) hold. Taking a subsequence if necessary, we may assume that there exist $(\phi,\psi)\in \textbf{H}^1$ and $(\theta,\chi)\in \textbf{L}^2$ such that $(\phi_j,\psi_j)\rightarrow (\phi,\psi)$ in $L^2_{\mathrm{loc}}(\R)\times L^2_{\mathrm{loc}}(\R)$ and weakly in $\textbf{H}^1$ as $j\rightarrow\infty$ and such that $(\theta_j,\chi_j)\rightharpoonup (\theta,\chi)$ weakly in $\textbf{L}^2$. By the weak lower-semicontinuity of the $\textbf{L}^2$-norm we obtain
\begin{equation}\label{lowersemicontinuity}\begin{array}{lcl}\nrm{(\theta,\chi)}_{\textbf{L}^2}+\frac{1}{\delta}|\ip{(\theta,\chi),(\phi_0^{-},\psi_0^{-})}|&\leq&\overline{\Lambda}^{+}(\delta).\end{array}\end{equation}
For any pair of test functions $(\zeta_1,\zeta_2)\in C_c^\infty(\R)\times C_c^\infty(\R)$ we have
\begin{equation}\label{weaksolutions}\begin{array}{lcl}\ip{(\theta_j,\chi_j),(\zeta_1,\zeta_2)}&=&\ip{\overline{\mathcal{L}}_{h_j,\delta}^{+}(\phi_j,\psi_j),(\zeta_1,\zeta_2)}\\[0.2cm]
&=&\ip{(\phi_j,\psi_j),\overline{\mathcal{L}}_{h_j,\delta}^{-}(\zeta_1,\zeta_2)}.\end{array}\end{equation}
Since $\overline{u}_0$ is a bounded function, the limit $\tilde{u}_h-\overline{u}_0\rightarrow 0$ in $H^1$ implies that also $\tilde{u}_h\rightarrow \overline{u}_0$ in $L^\infty$. In particular, we can choose $h'>0$ and $N>0$ in such a way that $|\tilde{u}_h|\leq N$ and $|\overline{u}_0|\leq N$ for all $0<h\leq h'$. Since $g_u$ is Lipschitz continuous on $[-N,N]$, there is a constant $K>0$ such that $|g_u(x)-g_u(y)|\leq K|x-y|$ for all $x,y\in[-N,N]$. We obtain
\begin{equation}\begin{array}{lcl}\lim\limits_{h\downarrow 0}\  \nrm{g_u(\tilde{u}_h)-g_u(\overline{u}_0)}_{L^2}^2&=&\lim\limits_{h\downarrow 0}\  \int(g_u(\tilde{u}_h)-g_u(\overline{u}_0))^2
dx\\[0.2cm]
&\leq & \lim\limits_{h\downarrow 0}\ \int K^2(\tilde{u}_h-\overline{u}_0)^2dx\\[0.2cm]
&\leq& \lim\limits_{h\downarrow 0}\ K^2\nrm{\tilde{u}_h-\overline{u}_0}_{\textbf{L}^2}^2\\[0.2cm]
&=&0,\end{array}\end{equation}
together with
\begin{equation}\begin{array}{lcl}\lim\limits_{h\downarrow 0}\enskip\nrm{g_u(\tilde{u}_h)\zeta_1-g_u(\overline{u}_0)\zeta_1}_{L^2}&
\leq &\lim\limits_{h\downarrow 0}\enskip\nrm{\zeta_1}_\infty\nrm{g_u(\tilde{u}_h)-g_u(\overline{u}_0)}_{L^2}\\[0.2cm]
&=&0.\end{array}\end{equation}
Furthermore, we know that $\tilde{c}_h\rightarrow c_0$ as $h\downarrow 0$, which gives
\begin{equation}\begin{array}{lcl}\lim\limits_{h\downarrow 0}\enskip\nrm{\tilde{c}_h\zeta_1'-c_0\zeta_1'}_{L^2}&=&\lim\limits_{h\downarrow 0}\enskip\nrm{\tilde{c}_h\zeta_2'-c_0\zeta_2'}_{L^2}\\[0.2cm]
&=&0.\end{array}\end{equation}
Finally, Lemma \ref{eigenschappenDeltah} implies
\begin{equation}\begin{array}{lcl}\lim\limits_{h\downarrow 0}\enskip\nrm{\Delta_h\zeta_1-\zeta_1''}_{L^2}&=&0,\end{array}\end{equation}
which means that
\begin{equation}
\begin{array}{lcl}
\nrm{\overline{\mathcal{L}}_{h_j,\delta}^{-}(\zeta_1,\zeta_2)-(\overline{\mathcal{L}}_0^{-}+\delta)(\zeta_1,\zeta_2)}_{\textbf{L}^2}&\rightarrow & 0
\end{array}
\end{equation}
as $j\rightarrow\infty$. Sending $j\rightarrow\infty$ in (\ref{weaksolutions}), this yields
\begin{equation}\begin{array}{lcl}\ip{(\theta,\chi),(\zeta_1,\zeta_2)}&=&\ip{(\phi,\psi),(\overline{\mathcal{L}}_0^{-}+\delta)(\zeta_1,\zeta_2)}.\end{array}\end{equation}
In particular, we see that $(\phi,\psi)$ is a weak solution to $(\overline{\mathcal{L}}_0^{+}+\delta)(\phi,\psi)=(\theta,\chi)$. Since $\phi\in H^1$, $\psi\in L^2$, $\theta\in L^2$ and
\begin{equation}\begin{array}{lcl}\phi''&=&c_0\phi'-g_u(\overline{u}_0)\phi+\delta\phi+\psi-\theta,\end{array}\end{equation}
we get $\phi''\in L^2$ and hence $\phi\in H^2$. Since we already know that $\psi\in H^1$, we may apply Lemma \ref{eigenschappenL0} and (\ref{lowersemicontinuity}) to obtain
\begin{equation}\begin{array}{lcl}\nrm{(\phi,\psi)}_{H^2(\R)\times H^1(\R)}&\leq& C_2[\nrm{(\theta,\chi)}_{\textbf{L}^2}+\frac{1}{\delta}|\ip{(\theta,\chi),(\phi_0^{-},\psi_0^{-})}|]\\[0.2cm]
&\leq& C_2\overline{\Lambda}^{+}(\delta).\end{array}\end{equation}
\qed

The next result is the analogue
of Lemma \ref{lemma6bewijs1} for the setting
where we are considering a spectral
set $M \subset \mathbb{C}$ that satisfies
(\asref{Massumption}{\text{h}}).
The proof is omitted as it is almost identical to that of Lemma \ref{lemma6bewijs1}. We recall the constant
$C_3>0$ from Lemma \ref{eigenschappenL0}, which only depends on the choice of the set $M\subset\C$.

\begin{lemma}\label{lemma6bewijs1.5} Assume that (\asref{aannamespuls}{\text{H}}),(\asref{extraaannamespuls}{\text{H}}), (\asref{aannamesconstanten}{\text{H}}) and (\asref{aannames}{\text{H}\alpha1}) are satisfied. Let $M\subset \C$ satisfy (\asref{Massumption}{\text{h}}). There exists a sequence $\{(\lambda_j,h_j,\phi_j,\psi_j)\}$ in $M\times (0,1)\times \textbf{H}^1$ with the following properties.
\begin{enumerate}
\item We have $\lim\limits_{j\rightarrow\infty}h_j=0$, $\lim\limits_{j\rightarrow\infty}\lambda_j=\lambda$ for some $\lambda\in M$ and $ \nrm{(\phi_j,\psi_j)}_{\textbf{H}^1}=1$ for all $j$.
\item The pair $(\theta_j,\chi_j)=\overline{\mathcal{L}}_{h_j,\lambda_j}^{+}(\phi_j,\psi_j)$ satisfies
\begin{equation}\begin{array}{lcl}\lim\limits_{j\rightarrow\infty} \nrm{(\theta_j,\chi_j)}_{\textbf{L}^2}&=&\overline{\Lambda}^{+}(M).\end{array}\end{equation}
\item There exist $(\phi,\psi)\in \textbf{H}^1$ and $(\theta,\chi)\in \textbf{L}^2$ such that as $j\rightarrow\infty$ $(\phi_j,\psi_j)\rightarrow (\phi,\psi)$ in $L^2_{\mathrm{loc}}(\R)\times L^2_{\mathrm{loc}}(\R)$ and weakly in $\textbf{H}^1$ and such that $(\theta_j,\chi_j)\rightharpoonup (\theta,\chi)$ weakly in $\textbf{L}^2$.
\item The pair $(\phi,\psi)$ is a weak solution to $(\overline{\mathcal{L}}_0^{+}+\lambda)(\phi,\psi)=(\theta,\chi)$.
\item We have the bound \begin{equation}\begin{array}{lcl}\nrm{(\phi,\psi)}_{H^2(\R)\times H^1(\R)}&\leq &C_3\overline{\Lambda}^{+}(M).\end{array}\end{equation}
\end{enumerate}
The same statements hold upon replacing $\overline{\mathcal{L}}_{h,\lambda_j}^+$, $\overline{\Lambda}^+(M)$ and $\overline{\mathcal{L}}_{0}^+$ by $\overline{\mathcal{L}}_{h,\lambda_j}^-$, $\overline{\Lambda}^-$ and $\overline{\mathcal{L}}_{0}^-$.
\end{lemma}

 In our arguments below, we often consider the sequences $\{(h_j,\phi_j,\psi_j)\}$ and $\{(\lambda_j,h_j,\phi_j,\psi_j)\}$ defined in Lemmas \ref{lemma6bewijs1} and \ref{lemma6bewijs1.5} in a similar fashion. To streamline our notation, we simply write $\{(\lambda_j,h_j,\phi_j,\psi_j)\}$ for all these sequences, with the understanding that $\lambda_j=\delta$ when referring to Lemma \ref{lemma6bewijs1}.
As argued in the proof of Lemma \ref{lemma6bewijs1}, it is possible to choose $\overline{h}>0$ in such a way that
\begin{equation}\label{definitieoverlineh}\begin{array}{lclcl}c_*&:=&\inf_{0\leq h\leq \overline{h}}|\tilde{c}_h|&>& 0,\\[0.2cm]
g_*&:=&\sup_{0\leq h\leq \overline{h}}\nrm{g_u(\tilde{u}_h)}_\infty&<&\infty.\end{array}\end{equation}
By taking a subsequence if necessary, we assume from now on that $h_j<\overline{h}$ for all $j$.

It remains to find a positive lower bound for $\nrm{(\phi,\psi)}_{\textbf{L}^2}$. An essential step
to accomplish this is to keep the derivatives
$(\phi_j', \psi_j')$ under control. This can be achieved
by exploiting the results for $\Delta_h$ derived
in {\S}\ref{sec:sp:conv:prlm}.
\begin{lemma}\label{lemma6bewijs2} Assume that (\asref{familyassumption}{h}), (\asref{aannamespuls}{\text{H}}), (\asref{aannamesconstanten}{\text{H}}) and (\asref{aannames}{\text{H}\alpha1}) are satisfied. Consider the setting of Proposition \ref{lemma6equivalent} and Lemma \ref{lemma6bewijs1} or Lemma \ref{lemma6bewijs1.5}. Then there exists a constant $B>0$, such that for all $j$ we have the bound
\begin{equation}\label{equation9equivalent'''}\begin{array}{lcl}B\nrm{(\phi_j,\psi_j)}_{\textbf{L}^2}^2&\geq& c_*^2\nrm{(\phi_j',\psi_j')}_{\textbf{L}^2}^2-4\nrm{(\theta_j,\chi_j)}_{\textbf{L}^2}^2.\end{array}\end{equation}
\end{lemma}
\vspace*{4pt}\noindent\textit{Proof.} We first consider the sequence for $\overline{\Lambda}^+$. Using $\overline{\mathcal{L}}_{h_j,\lambda_j}^{+}(\phi_j,\psi_j)=(\theta_j,\chi_j)$ and $\Re \,\ip{\Delta_{h_j}\phi_j,\phi_j'}=0=\Re \,\ip{\phi_j,\phi_j'}=\Re \,\ip{\psi_j,\psi_j'}$, which follow from Lemma \ref{eigenschappencomplexeDeltah}, we obtain
\begin{equation}\begin{array}{lcl}\Re \,\ip{(\theta_j,\chi_j),(\phi_j',\psi_j')}&=&\Re \,\ip{\overline{\mathcal{L}}_{h_j,\lambda_j}^{+}(\phi_j,\psi_j),(\phi_j',\psi_j')}\\[0.2cm]
&=&\Re \,\ip{ \tilde{c}_{h_j}\phi_j'-\Delta_{h_j}\phi_j-g_u(\tilde{u}_{h_j})\phi_j+\lambda_j\phi_j+\psi_j,\phi_j'}\\[0.2cm]
&&\qquad +\Re \,\ip{-\rho \phi_j+ \tilde{c}_{h_j}\psi_j'+\gamma\rho\psi_j+\lambda_j\psi_j,\psi_j'}\\[0.2cm]
&=& \tilde{c}_{h_j}\nrm{\phi_j'}_{L^2}^2-\Re \,\ip{g_u(\tilde{u}_{h_j})\phi_j,\phi_j'}+\Re \,\ip{\psi_j,\phi_j'}\\[0.2cm]
&&\qquad +\Re \,\ip{\lambda_j\phi_j,\phi_j'}-\rho\Re \,\ip{\phi_j,\psi_j'}\\[0.2cm]
&&\qquad + \tilde{c}_{h_j}\nrm{\psi_j'}^2_{L^2}+\Re \,\ip{\lambda_j\psi_j,\psi_j'}\\[0.2cm]
&=& \tilde{c}_{h_j}\nrm{(\phi_j',\psi_j')}_{\textbf{L}^2}^2-\Re \,\ip{g_u(\tilde{u}_{h_j})\phi_j,\phi_j'}+(1+\rho)\ip{\psi_j,\phi_j'}\\[0.2cm]
&&\qquad +\Re \,\ip{\lambda_j(\phi_j,\psi_j),(\phi_j',\psi_j')}.\end{array}\end{equation}

We write $\lambda_{\mathrm{max}}=\delta_0$ in the setting of Lemma \ref{lemma6bewijs1} or $\lambda_{\mathrm{max}}=\max\{|z|:z\in M\}$ in the setting of Lemma \ref{lemma6bewijs1.5}. We write
\begin{equation}\begin{array}{lcl}G&=&\lambda_{\mathrm{max}}\nrm{(\phi_j,\psi_j)}_{\textbf{L}^2}\nrm{(\phi_j',\psi_j')}_{\textbf{L}^2}
+g_*\nrm{\phi_j}_{L^2}
\nrm{(\phi_j',\psi_j')}_{\textbf{L}^2}.\end{array}\end{equation}
Using the Cauchy-Schwarz inequality we now obtain
\begin{equation}\begin{array}{lcl}G&\geq&\lambda_{\mathrm{max}}\nrm{(\phi_j,\psi_j)}_{\textbf{L}^2}\nrm{(\phi_j',\psi_j')}_{\textbf{L}^2} +\nrm{g_u(\tilde{u}_{h_j})}_{L^\infty}\nrm{\phi_j}_{L^2}\nrm{\phi_j'}_{L^2}\\[0.2cm]
&\geq & \mathrm{sign}(\tilde{c}_{h_j})\Big(-\Re \,\ip{\lambda_j(\phi_j,\psi_j),(\phi_j',\psi_j')}+\Re \,\ip{g_u(\tilde{u}_{h_j})\phi_j,\phi_j'}\Big) \\[0.2cm]
&=& \mathrm{sign}(\tilde{c}_{h_j})\Big(\tilde{c}_{h_j}\nrm{(\phi_j',\psi_j')}_{\textbf{L}^2}^2+(1+\rho)\Re \,\ip{\psi_j,\phi_j'} -\Re \,\ip{(\theta_j,\chi_j),
(\phi_j',\psi_j')}\Big)\\[0.2cm]
&\geq &|\tilde{c}_{h_j}|\nrm{(\phi_j',\psi_j')}^2_{\textbf{L}^2}-(1+\rho)\nrm{\psi_j}_{L^2}\nrm{\phi_j'}_{L^2}-\nrm{(\theta_j,\chi_j)}_{\textbf{L}^2}\nrm{(\phi_j',\psi_j')}_{\textbf{L}^2}\\[0.2cm]
&\geq &c_*\nrm{(\phi_j',\psi_j')}_{\textbf{L}^2}^2-(1+\rho)\nrm{\psi_j}_{L^2}\nrm{(\phi_j',\psi_j')}_{\textbf{L}^2}-\nrm{(\theta_j,\chi_j)}_{\textbf{L}^2}\nrm{(\phi_j',\psi_j')}_{\textbf{L}^2}.\end{array}\end{equation}
This implies
\begin{equation}\begin{array}{lcl}c_*\nrm{(\phi_j',\psi_j')}_{\textbf{L}^2}&\leq &g_*\nrm{\phi_j}_{L^2}+(1+\rho)\nrm{\psi_j}_{L^2}+\nrm{(\theta_j,\chi_j)}_{\textbf{L}^2}+\lambda_{\mathrm{max}}\nrm{(\phi_j,\psi_j)}_{\textbf{L}^2}.\end{array}\end{equation}
Squaring this equation and
using the standard inequality $2\mu\omega\leq \mu^2+\omega^2$, this implies that
\begin{equation}\begin{array}{lcl}c_*^2\nrm{(\phi_j',\psi_j')}_{\textbf{L}^2}^2&\leq &4g_*^2\nrm{\phi_j}_{L^2}^2
+4(1+\rho)^2\nrm{\psi_j}_{L^2}^2\\[0.2cm]
&&\qquad +4\nrm{(\theta_j,\chi_j)}_{\textbf{L}^2}^2+4\lambda_{\mathrm{max}}^2\nrm{(\phi_j,\psi_j)}_{\textbf{L}^2}^2.\end{array}\end{equation}
In particular, we see
\begin{equation}\begin{array}{lcl}4\Big(\max\{
g_*^2,(1+\rho)^2\}
+\lambda_{\mathrm{max}}^2\Big)\nrm{(\phi_j,\psi_j)}_{\textbf{L}^2}^2&\geq& c_*^2\nrm{(\phi_j',\psi_j')}_{\textbf{L}^2}^2-4\nrm{(\theta_j,\chi_j)}_{\textbf{L}^2}^2.\end{array}\end{equation}

We now look at the sequence for $\overline{\Lambda}^-$. Using $\overline{\mathcal{L}}_{h_j,\lambda_j}^{-}(\phi_j,\psi_j)=(\theta_j,\chi_j)$ and $\Re \,\ip{\Delta_{h_j}\phi_j,\phi_j'}=0=\Re \,\ip{\phi_j,\phi_j'}=\Re \,\ip{\psi_j,\psi_j'}$, which follow from Lemma \ref{eigenschappencomplexeDeltah}, we obtain
\begin{equation}\begin{array}{lcl}\Re \,\ip{(\theta_j,\chi_j),(\phi_j',\psi_j')}&=&\Re \,\ip{\overline{\mathcal{L}}_{h_j,\lambda_j}^{-}(\phi_j,\psi_j),(\phi_j',\psi_j')}\\[0.2cm]
&=&\Re \,\ip{ -\tilde{c}_{h_j}\phi_j'-\Delta_{h_j}\phi_j-g_u(\tilde{u}_h)\phi_j+\lambda_j\phi_j-\rho\psi_j,\phi_j'}\\[0.2cm]
&&\qquad+\Re \,\ip{ \phi_j- \tilde{c}_h\psi_j'+\gamma\rho\psi_j+\lambda_j\psi_j,\psi_j'}\\[0.2cm]
&=& -\tilde{c}_{h_j}\nrm{\phi_j'}_{L^2}^2-\Re \,\ip{g_u(\tilde{u}_h)\phi_j,\phi_j'}-\rho\Re \,\ip{\psi_j,\phi_j'}\\[0.2cm]
&&\qquad +\Re \,\ip{\lambda_j\phi_j,\phi_j'}+\Re \,\ip{\phi_j,\psi_j'}\\[0.2cm]
&&\qquad - \tilde{c}_{h_j}\nrm{\psi_j'}^2_{L^2}+\Re \,\ip{\lambda_j\psi_j,\psi_j'}\\[0.2cm]
&=& -\tilde{c}_{h_j}\nrm{(\phi_j',\psi_j')}_{\textbf{L}^2}^2-\Re \,\ip{g_u(\tilde{u}_h)\phi_j,\phi_j'}+(1+\rho)\ip{\psi_j,\phi_j'}\\[0.2cm]
&&\qquad +\Re \,\ip{\lambda_j(\phi_j,\psi_j),(\phi_j',\psi_j')}.\end{array}\end{equation}
We write
\begin{equation}\begin{array}{lcl}G&=&\lambda_{\mathrm{max}}\nrm{(\phi_j,\psi_j)}_{\textbf{L}^2}\nrm{(\phi_j',\psi_j')}_{\textbf{L}^2}
+g_*\nrm{\phi_j}_{L^2}
\nrm{(\phi_j',\psi_j')}_{\textbf{L}^2}.\end{array}\end{equation}
Using the Cauchy-Schwarz inequality we now obtain
\begin{equation}\begin{array}{lcl}G&\geq &\lambda_{\mathrm{max}}\nrm{(\phi_j,\psi_j)}_{\textbf{L}^2}\nrm{(\phi_j',\psi_j')}_{\textbf{L}^2} +\nrm{g_u(\tilde{u}_{h_j})}_{L^\infty}\nrm{\phi_j}_{L^2}\nrm{\phi_j'}_{L^2}\\[0.2cm]
&\geq & -\mathrm{sign}(\tilde{c}_{h_j})\Big(-\Re \,\ip{\lambda_j(\phi_j,\psi_j),(\phi_j',\psi_j')}+\Re \,\ip{g_u(\tilde{u}_{h_j})\phi_j,\phi_j'}\Big) \\[0.2cm]
&=&  -\mathrm{sign}(\tilde{c}_{h_j})\Big(-\tilde{c}_{h_j}\nrm{(\phi_j',\psi_j')}_{\textbf{L}^2}^2-(1+\rho)\Re \,\ip{\psi_j,\phi_j'}-\Re \,\ip{(\theta_j,\chi_j),
(\phi_j',\psi_j')}\Big)\\[0.2cm]
&\geq &|\tilde{c}_{h_j}|\nrm{(\phi_j',\psi_j')}^2_{\textbf{L}^2}-(1+\rho)\nrm{\psi_j}_{L^2}\nrm{\phi_j'}_{L^2}-\nrm{(\theta_j,\chi_j)}_{\textbf{L}^2}\nrm{(\phi_j',\psi_j')}_{\textbf{L}^2}\\[0.2cm]
&\geq &c_*\nrm{(\phi_j',\psi_j')}_{\textbf{L}^2}^2-(1+\rho)\nrm{\psi_j}_{L^2}\nrm{(\phi_j',\psi_j')}_{\textbf{L}^2} -\nrm{(\theta_j,\chi_j)}_{\textbf{L}^2}\nrm{(\phi_j',\psi_j')}_{\textbf{L}^2}.\end{array}\end{equation}This is the same equation that we derived for $\overline{\Lambda}^+$. Hence, we again obtain
\begin{equation}\begin{array}{lcl}B\nrm{(\phi_j,\psi_j)}_{\textbf{L}^2}^2&\geq& c_*^2\nrm{(\phi_j',\psi_j')}_{\textbf{L}^2}^2-4\nrm{(\theta_j,\chi_j)}_{\textbf{L}^2}^2,\end{array}\end{equation}
where
\begin{equation}\begin{array}{lcl}B&=&4\Big(\max\{
g_*^2,(1+\rho)^2\}
+\lambda_{\mathrm{max}}^2\Big).\end{array}\end{equation}\qed

The next step is to show that the $L^2$-mass of $\phi_j$
can be concentrated in a compact interval.
We heavily exploit the bistable structure of the non-linearity $g$ to accomplish this. Moreover, we are aided by the fact that the off-diagonal elements are constant, which allows us to keep the cross-terms under control. In fact, one might be tempted to think that it is sufficient to note that the eigenvalues of the matrix $\left(\begin{array}{ll}
-g_u(0)& 1\\ -\rho &\gamma\rho
\end{array}\right)$ all have positive real part, as then one would be able to find a basis in which this matrix is positive definite. However, passing over to another basis destroys the structure of the diffusion
terms and, therefore, does not give any insight.
\begin{lemma}\label{lemma6bewijs3} Assume that (\asref{familyassumption}{h}), (\asref{aannamespuls}{\text{H}}), (\asref{aannamesconstanten}{\text{H}}) and (\asref{aannames}{\text{H}\alpha1}) are satisfied. Consider the setting of Proposition \ref{lemma6equivalent} and Lemma \ref{lemma6bewijs1} or Lemma \ref{lemma6bewijs1.5}. There exist positive constants $a$ and $m$, depending only on our choice of $(\tilde{u}_h,\tilde{w}_h)$, such that we have the following inequality for all $j$
\begin{equation}\label{equation10equivalent}\begin{array}{lcl}\frac{1}{\rho}(a+g_*)\int\limits_{|x|\leq m}|\phi_j(x)|^2dx&\geq &\Big(\frac{1}{2}\min\{a,\frac{1}{2}\rho\gamma\}+\lambda_{\mathrm{min}}\Big)\nrm{(\phi_j,\psi_j)}_{\textbf{L}^2}^2\\[0.2cm]
&&\qquad -\frac{1}{2\min\{a,\frac{1}{2}\rho\gamma\}}\nrm{(\theta_j,\chi_j)}_{\textbf{L}^2}^2 -\beta
\nrm{(\theta_j,\chi_j)}_{\textbf{L}^2}^2.\end{array}\end{equation}
Here we write $\lambda_{\mathrm{min}}=0$ in the setting of Lemma \ref{lemma6bewijs1} or $\lambda_{\mathrm{min}}=\min\{\Re \,\lambda:\lambda\in M\}$ in the setting of Lemma \ref{lemma6bewijs1.5}, together with
\begin{equation}
\begin{array}{lcl}
\beta&=&\frac{1-\rho}{\rho}\frac{1}{4(\frac{\rho}{1-\rho}\frac{1}{2}\gamma\rho+\gamma\rho+\lambda_{\mathrm{min}})}.
\end{array}
\end{equation}
\end{lemma}
\vspace*{4pt}\noindent\textit{Proof.} Again we first look at the sequence for $\overline{\Lambda}^+$. We know that $\tilde{u}_h-\overline{u}_0\rightarrow 0$ in $H^1$ as $h\downarrow 0$. Hence, it follows that $\tilde{u}_h-\overline{u}_0\rightarrow 0$ in $L^\infty$ and, therefore, also $g_u(\tilde{u}_h)-g_u(\overline{u}_0)\rightarrow 0$ in $L^\infty$ as $h\downarrow 0$.  By the bistable nature of our non-linearity $g$, we can choose $m$ to be a positive constant such that for all $h\in[0,\overline{h}]$ (by making $\overline{h}$ smaller if necessary)
\begin{equation}\min_{|x|\geq m}[-g_u(\tilde{u}_h(x))]\geq a:=\frac{1}{2}r_0>0.\end{equation}
Here $r_0$ is the constant appearing in the choice of our function $g$ in (\asref{aannamesconstanten}{\text{H}}). Then we obtain, using $\Re \,\ip{\phi_j',\phi_j}=\Re \,\ip{\psi_j',\psi_j}=0$ and $\Re \,\ip{-\Delta_{h_j}\phi_j,\phi_j}\geq 0$, which we know from Lemma \ref{eigenschappencomplexeDeltah}, that
\begin{equation}\label{afschbates1}\begin{array}{lcl}\Re \,\ip{(\theta_j,\chi_j),(\phi_j,\psi_j)}&=&\Re \,\ip{\overline{\mathcal{L}}_{h_j,\lambda_j}^{+}(\phi_j,\psi_j),(\phi_j,\psi_j)}\\[0.2cm]
&\geq&\Re \,\ip{-g_u(\tilde{u}_{h_j})\phi_j,\phi_j}+\Re \,\ip{\psi_j,\phi_j}\\[0.2cm]
&&\qquad -\rho\Re \,\ip{\psi_j,\phi_j}+\gamma\rho\nrm{\psi_j}_{L^2}^2+\lambda_{\mathrm{min}}\nrm{(\phi_j,\psi_j)}^2\\[0.2cm]
&\geq &\min_{|x|\geq m}\{-g_u(\tilde{u}_{h_j}(x))\}\int_{|x|\geq m}|\phi_j(x)|^2dx\\[0.2cm]
&&\qquad -\nrm{g_u(\tilde{u}_{h_j})}_{L^\infty}\int\limits_{|x|\leq m}|\phi_j(x)|^2dx+(1-\rho)\Re \,\ip{\psi_j,\phi_j}\\[0.2cm]
&&\qquad +\gamma\rho\nrm{\psi_j}_{L^2}^2+\lambda_{\mathrm{min}}\nrm{(\phi_j,\psi_j)}^2\\[0.2cm]
&\geq &a\nrm{\phi_j}_{L^2}^2-(a+g_*)\int\limits_{|x|\leq m}|\phi_j(x)|^2dx+(1-\rho)\Re \,\ip{\psi_j,\phi_j}\\[0.2cm]
&&\qquad +\gamma\rho\nrm{\psi_j}_{L^2}^2+\lambda_{\mathrm{min}}\nrm{(\phi_j,\psi_j)}^2.\end{array}\end{equation}
We assumed that $0<\rho<1$ so we see that $\frac{1-\rho}{-\rho}<0$. We set
\begin{equation}
\begin{array}{lcl}
\beta_j^+&=&\frac{1}{4(\frac{\rho}{1-\rho}\frac{1}{2}\gamma\rho+\gamma\rho+\Re \,\lambda_j)}.
\end{array}
\end{equation}
Now we obtain
\begin{equation}\label{crosstermafschatting}\begin{array}{lcl}\Re \,\ip{\chi_j,\psi_j}&\leq&\nrm{\chi_j}_{L^2}\nrm{\psi_j}_{L^2}\\[0.2cm]
&=&\frac{1}{\sqrt{2(\frac{\rho}{1-\rho}\frac{1}{2}\gamma\rho+\gamma\rho+\Re \,\lambda_j)}}\nrm{\chi_j}_{L^2}
\sqrt{2(\frac{\rho}{1-\rho}\frac{1}{2}\gamma\rho+\gamma\rho+\Re \,\lambda_j)}\nrm{\psi_j}_{L^2}\\[0.2cm]
&\leq&\frac{1}{4(\frac{\rho}{1-\rho}\frac{1}{2}\gamma\rho+\gamma\rho+\Re \,\lambda_j)}
\nrm{\chi_j}_{L^2}^2+(\frac{\rho}{1-\rho}\frac{1}{2}\gamma\rho+\gamma\rho+\Re \,\lambda_j)\nrm{\psi_j}_{L^2}^2\\[0.2cm]
&=&\beta_j^+\nrm{\chi_j}_{L^2}^2+(\frac{\rho}{1-\rho}\frac{1}{2}\gamma\rho+\gamma\rho+\Re \,\lambda_j)\nrm{\psi_j}_{L^2}^2.
\end{array}\end{equation}
Note that the denominator $4(\frac{\rho}{1-\rho}\frac{1}{2}\gamma\rho+\gamma\rho+\Re \,\lambda_j)$ is never zero since we can assume that $\lambda_*$ is small enough to have $\Re \,\lambda_j\geq -\lambda_*>-\gamma\rho$. Using the identity
\begin{equation}\begin{array}{lcl}\chi_j&=&-\rho\phi_j+\tilde{c}_{h_j}\psi_j'+\gamma\rho\psi_j+\lambda_j\psi_j\end{array}\end{equation} and the fact that $\Re \,\ip{\psi_j',\psi_j}=0$, we also have
\begin{equation}\begin{array}{lcl}\Re \,\ip{\chi_j,\psi_j}&=
&-\rho\Re \,\ip{\phi_j,\psi_j}
+(\gamma\rho+\Re \,\lambda_j)\nrm{\psi_j}_{L^2}^2.\end{array}\end{equation}
Hence, we must have that
\begin{equation}\begin{array}{lcl}(1-\rho)\Re \,\ip{\phi_j,\psi_j}&=&\frac{1-\rho}{\rho}\Big(-\Re \,\ip{\chi_j,\psi_j}+(\gamma\rho+\Re \,\lambda_j)\nrm{\psi_j}_{L^2}^2\Big)\\[0.2cm]
&\geq &\frac{1-\rho}{\rho}\Big(-\beta_j^+\nrm{\chi_j}_{L^2}^2-(\frac{\rho}{1-\rho}\frac{1}{2}\gamma\rho+\gamma\rho+\Re \,\lambda_j)\nrm{\psi_j}_{L^2}^2\\[0.2cm]
&&\qquad +(\gamma\rho+\Re \,\lambda_j)\nrm{\psi_j}_{L^2}^2\Big) \\[0.2cm]
&=&-\frac{1-\rho}{\rho}\beta_j^+
\nrm{\chi_j}_{L^2}^2-\frac{1}{2}\gamma\rho\nrm{\psi_j}_{L^2}^2.\end{array}\end{equation}
Combining this bound with (\ref{afschbates1}) yields the estimate
\begin{equation}\label{afschbates2}\begin{array}{lcl}&&\Re \,\ip{(\theta_j,\chi_j),(\phi_j,\psi_j)}\\[0.2cm]
&\geq &a\nrm{\phi_j}_{L^2}^2-(a+g_*)\int\limits_{|x|\leq m}|\phi_j(x)|^2dx+(1-\rho)\Re \,\ip{\psi_j,\phi_j}\\[0.2cm]
&&\qquad +\gamma\rho\nrm{\psi_j}_{L^2}^2+\lambda_{\mathrm{min}}\nrm{(\phi_j,\psi_j)}^2\\[0.2cm]
&\geq &a\nrm{\phi_j}_{L^2}^2-(a+g_*)\int\limits_{|x|\leq m}|\phi_j(x)|^2dx\\[0.2cm]
&&\qquad +\frac{1}{2}\gamma\rho\nrm{\psi_j}_{L^2}^2+\lambda_{\mathrm{min}}\nrm{(\phi_j,\psi_j)}^2-\frac{1-\rho}{\rho}\beta_j^+
\nrm{\chi_j}_{L^2}^2.\end{array}\end{equation}

We now look at the sequence for $\overline{\Lambda}^-$. Let $m$ and $a$ be as before. Then we obtain, using $\overline{\mathcal{L}}_{h_j,\lambda_j}^-(\phi_j,\psi_j)=(\theta_j,\chi_j)$, $\Re \,\ip{\phi_j',\phi_j}=\Re \,\ip{\psi_j',\psi_j}=0$ and $\Re \,\ip{-\Delta_{h_j}\phi_j,\phi_j}\geq 0$ that
\begin{equation}\label{beginafschbatesadj}\begin{array}{lcl}\Re \, \ip{(\theta_j,\chi_j),(\phi_j,\psi_j)}&=&\Re \,\ip{\overline{\mathcal{L}}_{h_j,\delta}^{-}(\phi_j,\psi_j),(\phi_j,\psi_j)}\\[0.2cm]
&\geq&\Re \,\ip{-g_u(\tilde{u}_h)\phi_j,\phi_j}+(1-\rho)\Re \,\ip{\psi_j,\phi_j}\\[0.2cm]
&&\qquad+\gamma\rho\nrm{\psi_j}_{L^2}^2+\lambda_{\mathrm{min}}\nrm{(\phi_j,\psi_j)}_{\textbf{L}^2}^2.\end{array}\end{equation}
We set
\begin{equation}
\begin{array}{lcl}
\beta_j^-&=&\frac{1}{4(\frac{1}{1-\rho}\frac{1}{2}\gamma\rho+\gamma\rho+\Re \,\lambda_j)}.
\end{array}
\end{equation}
Arguing as in (\ref{crosstermafschatting}) with different constants, we obtain
\begin{equation}\begin{array}{lcl}
\Re \,\ip{\theta_j,\phi_j}&\geq&-\nrm{\theta_j}_{L^2}\nrm{\phi_j}_{L^2}\\[0.2cm]

&\geq&-\frac{1}{4(a+\Re \,\lambda_j)}
\nrm{\theta_j}_{L^2}^2-(a+\Re\lambda_j)\nrm{\phi_j}_{L^2}^2\\[0.2cm]
&=&-\beta_j^-\nrm{\theta_j}_{L^2}^2-(a+\Re \,\lambda_j)\nrm{\phi_j}_{L^2}^2. \end{array}\end{equation}
Note that the denominator $4(a+\Re \,\lambda_j)$ is never zero since we can assume that $\lambda_*$ is small enough to have $\Re \,\lambda_j\geq -\lambda_*>-a$. Using the identity \begin{equation}\begin{array}{lcl}\theta_j&=&-\tilde{c}_{h_j}\phi_j'-\Delta_h\phi_j-g_u(\tilde{u}_h)\phi_j+\lambda_j\psi_j-\rho\phi_j\end{array}\end{equation} and the fact that $\Re \,\ip{\phi_j',\phi_j}=0$, we also have
\begin{equation}\begin{array}{lcl}\Re \,\ip{\theta_j,\phi_j}&=
&\Re \,\ip{-\Delta_h\phi_j,\phi_j}+\Re \,\ip{-g_u(\tilde{u}_h)\phi_j,\phi_j}\\[0.2cm]&&
+\Re \,\lambda_j\nrm{\psi_j}_{L^2}^2-\rho\Re \,\ip{\phi_j,\psi_j}.\end{array}\end{equation}
Hence, we must have that
\begin{equation}\begin{array}{lcl}(1-\rho)\Re \,\ip{\phi_j,\psi_j}&=&\frac{1-\rho}{\rho}\Big(-\Re \,\ip{\theta_j,\phi_j}+\Re \,\ip{-\Delta_h\phi_j,\phi_j}\\[0.2cm]
&&\qquad+\Re \,\ip{-g_u(\tilde{u}_h)\phi_j,\phi_j}
+\Re \,\lambda_j\nrm{\psi_j}_{L^2}^2\Big)\\[0.2cm]
&\geq &\frac{1-\rho}{\rho}\Big(-\beta_j^-\nrm{\theta_j}_{L^2}^2-\big(a+\Re\lambda_j\big)\nrm{\phi_j}_{L^2}^2\\[0.2cm]
&&\qquad +\Re \,\ip{-g_u(\tilde{u}_h)\phi_j,\phi_j}+\Re \,\lambda_j\nrm{\psi_j}_{L^2}^2\Big)\\[0.2cm]
&=&\frac{1-\rho}{\rho}\Big(-\beta_j^-\nrm{\theta_j}_{L^2}^2-a\nrm{\phi_j}_{L^2}^2+\Re \,\ip{-g_u(\tilde{u}_h)\phi_j,\phi_j}\Big).\end{array}\end{equation}
Combining this with the estimate (\ref{beginafschbatesadj}) and noting that $\frac{1-\rho}{\rho}+1=\frac{1}{\rho}$ yields
\begin{equation}\label{afschbates2adj}
    \begin{array}{lcl}\Re \,\ip{(\theta_j,\chi_j),(\phi_j,\psi_j)}
  &\geq& \frac{1}{\rho}\Re \,\ip{-g_u(\tilde{u}_h)\phi_j,\phi_j}+\lambda_{\mathrm{min}}\nrm{(\phi_j,\psi_j)}_{\textbf{L}^2}^2\\[0.2cm]
  &&\qquad +\gamma\rho\nrm{\psi_j}_{L^2}^2-a\frac{1-\rho}{\rho}\nrm{\phi_j}_{L^2}^2-\frac{1-\rho}{\rho}\beta_j^-\nrm{\theta_j}_{L^2}^2\\[0.2cm]
&\geq &\frac{1}{\rho}\Big(\min_{|x|\geq m}\{-g_u(\tilde{u}_h(x))\}\int_{|x|\geq m}|\phi_j|^2dx\\[0.2cm]
&&\qquad -\nrm{g_u(\tilde{u}_h)}_{L^\infty}\int\limits_{|x|\leq m}|\phi_j|^2dx\Big)+\lambda_{\mathrm{min}}\nrm{(\phi_j,\psi_j)}_{\textbf{L}^2}^2\\[0.2cm]
&&\qquad +\gamma\rho\nrm{\psi_j}_{L^2}^2-a\frac{1-\rho}{\rho}\nrm{\phi_j}_{L^2}^2-\frac{1-\rho}{\rho}\beta_j^-\nrm{\theta_j}_{L^2}^2\\[0.2cm]
&\geq& a\nrm{\phi_j}_{L^2}^2-\frac{1}{\rho}(a+g_*)\int\limits_{|x|\leq m}|\phi_j|^2dx+\gamma\rho\nrm{\psi_j}_{L^2}^2\\[0.2cm]
&&\qquad +\lambda_{\mathrm{min}}\nrm{(\phi_j,\psi_j)}_{\textbf{L}^2}^2  -\frac{1-\rho}{\rho}\beta_j^-\nrm{\theta_j}_{L^2}^2
    \end{array}
\end{equation}

Upon setting
\begin{equation}
\begin{array}{lcl}
\beta&=&\frac{1-\rho}{\rho}\min\big\{\frac{1}{4(\frac{\rho}{1-\rho}\frac{1}{2}\gamma\rho+\gamma\rho+\lambda_{\mathrm{min}})},\frac{1}{4(a+\lambda_{\mathrm{min}})}\big\},
\end{array}
\end{equation} we note that $\frac{1-\rho}{\rho}\beta_j^+\leq \beta$ and $\frac{1-\rho}{\rho}\beta_j^-\leq \beta$ for all $j$ since $\rho<1$ and since $\beta_j^+$ and $\beta_j^-$ are maximal for $\Re \,\lambda=\lambda_{\mathrm{min}}$. Therefore, in both cases, we obtain
\begin{equation}\begin{array}{lcl}\frac{1}{\rho}(a+g_*)\int\limits_{|x|\leq m}|\phi_j(x)|^2dx&\geq &a\nrm{\phi_j}_{L^2}^2+\frac{1}{2}\rho\gamma\nrm{\psi_j}_{L^2}^2
-\Re \,\ip{(\theta_j,\chi_j),(\phi_j,\psi_j)}\\[0.2cm]
&&\qquad -\beta
\nrm{(\theta_j,\chi_j)}_{\textbf{L}^2}^2+\lambda_{\mathrm{min}}\nrm{(\phi_j,\psi_j)}_{\textbf{L}^2}^2\\[0.2cm]
&\geq &\Big(\min\{a,\frac{1}{2}\rho\gamma\}+\lambda_{\mathrm{min}}\Big)\nrm{(\phi_j,\psi_j)}_{\textbf{L}^2}^2\\[0.2cm]
&&\qquad -\frac{\nrm{(\theta_j,\chi_j)}_{\textbf{L}^2}}{\sqrt{\min\{a,\frac{1}{2}\rho\gamma\}}}\sqrt{\min\{a,\frac{1}{2}\rho\gamma\}}\nrm{(\phi_j,\psi_j)}_{\textbf{L}^2}\\[0.2cm]
&&\qquad  -\beta
\nrm{(\theta_j,\chi_j)}_{\textbf{L}^2}\end{array}\end{equation}
and thus, again using the inequality $2\mu\omega\leq \mu^2+\omega^2$ for $\mu,\omega\in\R$,
\begin{equation}\begin{array}{lcl}\frac{1}{\rho}(a+g_*)\int\limits_{|x|\leq m}|\phi_j(x)|^2dx&\geq &\Big(\min\{a,\frac{1}{2}\rho\gamma\}+\lambda_{\mathrm{min}}\Big)\nrm{(\phi_j,\psi_j)}_{\textbf{L}^2}^2\\[0.2cm]
&&\qquad -\frac{1}{2}\frac{1}{\min\{a,\frac{1}{2}\rho\gamma\}}\nrm{(\theta_j,\chi_j)}_{\textbf{L}^2}^2\\[0.2cm] &&\qquad -\frac{1}{2}\min\{a,\frac{1}{2}\rho\gamma\}\nrm{(\phi_j,\psi_j)}_{\textbf{L}^2}^2\\[0.2cm]
&&\qquad -\beta
\nrm{(\theta_j,\chi_j)}_{\textbf{L}^2}^2\\[0.2cm]
&=&\Big(\frac{1}{2}\min\{a,\frac{1}{2}\rho\gamma\}+\lambda_{\mathrm{min}}\Big)\nrm{(\phi_j,\psi_j)}_{\textbf{L}^2}^2 \\[0.2cm]
&&\qquad -\frac{1}{2\min\{a,\frac{1}{2}\rho\gamma\}}\nrm{(\theta_j,\chi_j)}_{\textbf{L}^2}^2 -\beta
\nrm{(\theta_j,\chi_j)}_{\textbf{L}^2}^2,\end{array}\end{equation}
as desired.
\qed

\begin{lemma}\label{lemma6bewijs4} Assume that (\asref{familyassumption}{h}), (\asref{aannamespuls}{\text{H}}), (\asref{aannamesconstanten}{\text{H}}) and (\asref{aannames}{\text{H}\alpha1}) are satisfied. Consider the setting of Proposition \ref{lemma6equivalent} and Lemma \ref{lemma6bewijs1} or Lemma \ref{lemma6bewijs1.5}. There exist positive constants $C_4$ and $C_5$, depending only on our choice of $(\tilde{u}_h,\tilde{w}_h)$, such that for all $j$ we have
\begin{equation}\label{vergelijking}\begin{array}{lcl}\frac{1}{\rho}(a+g_*)\int\limits_{|x|\leq m}|\phi_j^2(x)|dx &\geq &C_4-C_5\nrm{(\theta_j,\chi_j)}_{\textbf{L}^2}^2.\end{array}\end{equation}
\end{lemma}
\vspace*{4pt}\noindent\textit{Proof.}
Without loss of generality we assume that $\frac{1}{2}\min\{a,\frac{1}{2}\rho\gamma\}+\lambda_{\mathrm{min}}>0$. Write
\begin{equation}\begin{array}{lcl}\mu&=&\frac{\frac{1}{2}\min\{a,\frac{1}{2}\rho\gamma\}+\lambda_{\mathrm{min}}}{c_*^2+B}.\end{array}\end{equation}
Adding $\mu$ times equation (\ref{equation9equivalent'''}) to equation (\ref{equation10equivalent}) gives
\begin{equation}\begin{array}{lcl}\frac{1}{\rho}(a+g_*)\int\limits_{|x|\leq m}|\phi_j(x)|^2 dx
&\geq & \mu c_*^2\nrm{(\phi_j',\psi_j')}_{\textbf{L}^2}^2 -4\mu \nrm{(\theta_j,\chi_j)}_{\textbf{L}^2}^2 \\[0.2cm]
&&\qquad +\frac{1}{2}(\min\{a,\frac{1}{2}\rho\gamma\}+\lambda_{\mathrm{min}})\nrm{(\phi_j,\psi_j)}_{\textbf{L}^2}^2 \\[0.2cm]
&&\qquad -\frac{1}{2(\min\{a,\frac{1}{2}\rho\gamma\}+\lambda_{\mathrm{min}})}\nrm{(\theta_j,\chi_j)}_{\textbf{L}^2}^2 \\[0.2cm]
&&\qquad -\beta
\nrm{(\theta_j,\chi_j)}_{\textbf{L}^2}^2 -B\mu\nrm{(\phi_j,\psi_j)}_{\textbf{L}^2}^2.\end{array}\end{equation}
We hence obtain
\begin{equation}\begin{array}{lcl}\frac{1}{\rho}(a+g_*)\int\limits_{|x|\leq m}|\phi_j(x)|^2dx&\geq &-C_5\nrm{(\theta_j,\chi_j)}_{\textbf{L}^2}^2+\mu c_*^2\nrm{(\phi_j',\psi_j')}_{\textbf{L}^2}^2 \\[0.2cm]
&&\qquad +\frac{1}{2}(\min\{a,\frac{1}{2}\rho\gamma\}+\lambda_{\mathrm{min}})\nrm{(\phi_j,\psi_j)}_{\textbf{L}^2}^2 \\[0.2cm]
&&\qquad  -B\mu\nrm{(\phi_j,\psi_j)}_{\textbf{L}^2}^2,\end{array}\end{equation}
where
\begin{equation}\begin{array}{lcl}
C_5& =& 4\mu+\frac{1}{2(\min\{a,\frac{1}{2}\rho\gamma\}+\lambda_{\mathrm{min}})}+\beta
\\[0.2cm]
&>&0.\end{array}\end{equation}
This allows us to compute
\begin{equation}\begin{array}{lcl}\frac{1}{\rho}(a+g_*)\int\limits_{|x|\leq m}|\phi_j(x)|^2dx&\geq &-C_5\nrm{(\theta_j,\chi_j)}_{\textbf{L}^2}^2+\mu c_*^2\nrm{(\phi_j',\psi_j')}_{\textbf{L}^2}^2 \\[0.2cm]
&&\qquad +\frac{1}{2}(\min\{a,\frac{1}{2}\rho\gamma\}+\lambda_{\mathrm{min}})\nrm{(\phi_j,\psi_j)}_{\textbf{L}^2}^2 \\[0.2cm]
&&\qquad  -B\mu\nrm{(\phi_j,\psi_j)}_{\textbf{L}^2}^2 \\[0.2cm]
&=& -C_5\nrm{(\theta_j,\chi_j)}_{\textbf{L}^2}^2+\mu c_*^2\nrm{(\phi_j',\psi_j')}_{\textbf{L}^2} \\[0.2cm]
&&\qquad +(\mu(c_*^2+B) -B\mu)\nrm{(\phi_j,\psi_j)}_{\textbf{L}^2}^2 \\[0.2cm]
&=&\mu c_*^2\nrm{(\phi_j,\psi_j)}_{\textbf{H}^1}^2-C_5\nrm{(\theta_j,\chi_j)}_{\textbf{L}^2}^2 \\[0.2cm]
&=&C_4-C_5\nrm{(\theta_j,\chi_j)}_{\textbf{L}^2}^2,\end{array}\end{equation}
where $C_4=\mu c_*^2>0$.\qed

\vspace*{4pt}\noindent\textit{Proof of Proposition \ref{lemma6equivalent}.} We first choose $0<\delta<\delta_0$ and consider the setting of Lemma \ref{lemma6bewijs1}. Sending $j\rightarrow\infty$ in (\ref{vergelijking}), Lemma \ref{lemma6bewijs1} implies
\begin{equation}\begin{array}{lcl} C_4-C_5\overline{\Lambda}^{\pm}(\delta) & \leq &C_4-C_5\lim\limits_{j\rightarrow\infty}\nrm{(\theta_j,\chi_j)}_{\textbf{L}^2}^2\\[0.2cm]
&\leq &\frac{1}{\rho}(a+g_*)\int\limits_{|x|\leq m}|\phi|^2dx\\[0.2cm]
&\leq &\frac{1}{\rho}(a+g_*)\nrm{(\phi,\psi)}_{H^2(\R)\times H^1(\R)}^2\\[0.2cm]
&\leq &\frac{1}{\rho}(a+g_*)C_2^2\overline{\Lambda}^{+}(\delta)^2.\end{array}\end{equation}
Solving this quadratic inequality, we obtain
\begin{equation}\begin{array}{lcl}\overline{\Lambda}^{\pm}(\delta)&\geq &\frac{-C_5+\sqrt{C_5^2+\frac{4}{\rho}(a+g_*)C_2^2C_4}}{\frac{2}{\rho}(a+g_*)C_2^2}\\[0.2cm]
&:=&\frac{2}{C_0'}.\end{array}\end{equation}

\noindent The analogous computation in the setting of Lemma \ref{lemma6bewijs1.5}  yields
\begin{equation}\begin{array}{lcl}\overline{\Lambda}^{+}(M)&\geq &\frac{-C_5+\sqrt{C_5^2+\frac{4}{\rho}(a+g_*)C_3^2C_4}}{\frac{2}{\rho}(a+g_*)C_3^2}\\[0.2cm]
& :=&\frac{2}{C_M'}.\end{array}\end{equation}
\qed

\vspace*{-10pt}

\section{The point and essential spectrum}\label{waveproperties}
In this section we discuss several properties of the operator that arises after linearising the travelling pulse MFDE (\ref{ditishetprobleem}) around our wave solution $(\overline{u}_h,\overline{w}_h)$. The main goals are to determine the Fredholm properties of this operator. In particular, we show that both the linearised operator and its adjoint have Fredholm index $0$ and that they both have a one-dimensional kernel. Moreover, we construct a family of kernel elements of the adjoint operator that converges to $(\phi_0^-,\psi_0^-)$, the kernel element of the operator $\mathcal{L}_0^-$.\par

Pick $0<h<\min\{h_*,\overline{h}\}$, where $h_*$ is given in Theorem \ref{Theorem1equivalent} and $\overline{h}$ is characterized by (\ref{definitieoverlineh}). We recall the operator $L_h:\textbf{H}^1\rightarrow \textbf{L}^2$, introduced in \S \ref{mainresults}, which acts as
\begin{equation}\begin{array}{lcl}L_h&=&\left(\begin{array}{ll}c_h\frac{d}{d\xi}-\Delta_h -g_u(\overline{u}_h)& 1\\ -\rho & c_h\frac{d}{d\xi}+\gamma\rho\end{array}\right).\end{array}\end{equation}
In addition, we write $L_h^*:\textbf{H}^1\rightarrow \textbf{L}^2$ for the formal adjoint of $L_h$, which is given by
\begin{equation}\begin{array}{lcl}L_h^*&=&\left(\begin{array}{ll}-c_h\frac{d}{d\xi}-\Delta_h -g_u(\overline{u}_h)& -\rho\\ 1 & -c_h\frac{d}{d\xi}+\gamma\rho\end{array}\right).\end{array}\end{equation}
We emphasize that $L_h$ and $L_h^*$ correspond to the operators $\overline{\mathcal{L}}_{h,0}^+$ and $\overline{\mathcal{L}}_{h,0}^-$ defined in \S \ref{singularoperator} respectively upon writing
\begin{equation}\begin{array}{lcl}(\tilde{u}_h,\tilde{w}_h)&=&(\overline{u}_h,\overline{w}_h),\\[0.2cm]
\tilde{c}_h&=&c_h\end{array}\end{equation}
for the family featuring in (\asref{familyassumption}{\text{h}}). Finally, we introduce the notation
\begin{equation}
\begin{array}{lclcl}
\Phi_h^+&=&(\phi_h^+,\psi_h^+)
&=&\frac{1}{\nrm{(\overline{u}_h',\overline{w}_h')}_{\textbf{L}^2}}(\overline{u}_h',\overline{w}_h'),\\[0.2cm]
\Phi_0^+&=&(\phi_0^+,\psi_0^+),\\[0.2cm]
\Phi_0^-&=&(\phi_0^-,\psi_0^-).
\end{array}
\end{equation}

The results of this section should be seen as a bridge between the singular perturbation theory developed in \S \ref{singularoperator} and the spectral analysis preformed in \S \ref{spectralanalysis}. Indeed, one might be tempted to think that most of the work required for the spectral analysis of the operator $L_h$ can already be found in Proposition \ref{equivalenttheorem4} and Proposition \ref{equivalenttheorem4version2}. However, the problem is that we have no control over the $\delta$-dependence of the interval $(0,h_0'(\delta))$, which contains all values of $h$ for which $L_h+\delta=\overline{\mathcal{L}}_{h,\delta}^+$ is invertible. In particular, for fixed $h>0$ we cannot directly conclude that $\overline{\mathcal{L}}_{h,\delta}^+$ is invertible for all $\delta$ in a subset of the positive real axis.\par

Our main task in this section is, therefore, to remove the $\delta$-dependence and study $L_h$ and $L_h^*$ directly. The main conclusions are summarized in the results below.
\begin{proposition}\label{LhlambdaFredholm} Assume that (\asref{aannamespuls}{\text{H}}), (\asref{aannamesconstanten}{\text{H}}), (\asref{aannames}{\text{H}\alpha1}) and (\asref{extraaannames}{\text{H}\alpha}) are satisfied. Then there exists a constant $\tilde{\lambda}>0$ such that for all $\lambda\in\C$ with $\Re \,\lambda>-\tilde{\lambda}$ and all $0<h<\min\{h_*,\overline{h}\}$ the operator $L_h+\lambda$ is Fredholm with index $0$.\end{proposition}
\begin{proposition}\label{samenvattingh5} Assume that (\asref{aannamespuls}{\text{H}}), (\asref{aannamesconstanten}{\text{H}}), (\asref{aannames}{\text{H}\alpha1}) and (\asref{extraaannames}{\text{H}\alpha}) are satisfied. Then there exists a constant $h_{**}>0$, together with a family $\Phi_h^-=(\phi_h^-,\psi^-)\in \textbf{H}^1$, defined for $0<h<h_{**}$, such that the following properties hold.
\begin{enumerate}
\item For each $0<h<h_{**}$ we have the identities
\begin{equation}\label{kernvanlh}\begin{array}{lcl}\ker (L_h)&=&\span\{\Phi_h^+\}\\[0.2cm]
&=&\{\Psi\in \textbf{L}^2:\ip{\Psi,\Theta}=0\text{ for all }\Theta\in \Range (L_h^*)\}\end{array}\end{equation}
and
\begin{equation}\label{kernvanlhster}\begin{array}{lcl}\ker(L_h^*)&=&\span\{\Phi_h^-\}\\[0.2cm]
&=&\{\Psi\in \textbf{L}^2:\ip{\Psi,\Theta}=0\text{ for all }\Theta\in \Range (L_h)\}.\end{array}\end{equation}
\item The family $\Phi_h^-$ converges to $\Phi_0^-$ in $\textbf{H}^1$ as $h\downarrow 0$.

\item Upon introducing the spaces
\begin{equation}\begin{array}{lcl}X_h&=&\{\Theta\in \textbf{H}^1: \ip{\Phi_h^-,\Theta}=0\}\end{array}\end{equation}
and
\begin{equation}\begin{array}{lcl}Y_h&=&\{\Theta\in \textbf{L}^2: \ip{\Phi_h^-,\Theta}=0\},\end{array}\end{equation}
the operator $L_h:X_h\rightarrow Y_h$ is invertible and there exists a constant $C_{\mathrm{unif}}>0$ such that for each $0<h< h_{**}$ we have the uniform bound
\begin{equation}
\begin{array}{lcl}
\nrm{L_h^{-1}}_{\mathcal{B}(Y_h,X_h)}\leq C_{\mathrm{unif}}.
\end{array}
\end{equation}
\end{enumerate}
\end{proposition}

A direct consequence of these results is that the zero eigenvalue of $L_h$ is simple. In addition, these results allow us to construct a quasi-inverse for $L_h$ that we use heavily in \S \ref{spectralanalysis} and \S \ref{sectiongreen}.
\begin{corollary}\label{simpleeigenvalue} Assume that (\asref{aannamespuls}{\text{H}}), (\asref{aannamesconstanten}{\text{H}}), (\asref{aannames}{\text{H}\alpha1}) and (\asref{extraaannames}{\text{H}\alpha}) are satisfied. Then for any $0<h< h_{**}$ the zero eigenvalue of $L_h$ is simple.\end{corollary}
\vspace*{4pt}\noindent\textit{Proof.} We can assume that $\ip{\Phi_h^-,\Phi_h^+}\neq 0$ for all $0<h< h_{**}$, since by Proposition \ref{samenvattingh5} $\ip{\Phi_h^-,\Phi_h^+}\rightarrow\ip{\Phi_0^-,\Phi_0^+}\neq 0$. Equation (\ref{kernvanlhster}) now implies that $\Phi_h^+\notin\Range(L_h)$, which together with (\ref{kernvanlh}) completes the proof. \qed

\begin{corollary}\label{vhwh} Assume that (\asref{aannamespuls}{\text{H}}), (\asref{aannamesconstanten}{\text{H}}), (\asref{aannames}{\text{H}\alpha1}) and (\asref{extraaannames}{\text{H}\alpha}) are satisfied. There exist linear maps
\begin{equation}\begin{array}{lcl}\gamma_h:\textbf{L}^2&\rightarrow&\R\\[0.2cm] L_h^{\mathrm{qinv}}:\textbf{L}^2&\rightarrow &\textbf{H}^1,\end{array}\end{equation}
such that for all $\Theta\in \textbf{L}^2$ and each $0<h<h_{**}$ the pair
\begin{equation}\begin{array}{lcl}(\gamma,\Psi)&=&(\gamma_h\Theta,L_h^{\mathrm{qinv}}\Theta)\end{array}\end{equation}
is the unique solution to the problem
\begin{equation}\begin{array}{lcl}L_h\Psi&=&\Theta+\gamma\Phi_h^+\end{array}\end{equation}
that satisfies the normalisation condition
\begin{equation}\begin{array}{lcl}\ip{\Phi_h^-,\Psi}&=&0.\end{array}\end{equation}\end{corollary}
\vspace*{4pt}\noindent\textit{Proof.} Fix $0<h< h_{**}$ and $\Theta\in \textbf{L}^2$. Upon defining
\begin{equation}\begin{array}{lcl}\gamma_h[\Theta]&=&-\frac{\ip{\Phi_h^-,\Theta}}{\ip{\Phi_h^-,\Phi_h^+}},\end{array}\end{equation}
we see that $\Theta+\gamma_h[\Theta]\Phi_h^+\in Y_h$. In particular, Proposition \ref{samenvattingh5} implies that the problem
\begin{equation}\begin{array}{lcl}L_h\Psi&=&\Theta+\gamma_h[\Theta]\Phi_h^+\end{array}\end{equation}
has a unique solution $\Psi\in X_h$, which we refer to as $L_h^{\mathrm{qinv}}\Theta$.\qed

The results in \cite{Faye2014,MPA} allow us to read off the Fredholm properties of $L_h$ from the behaviour of this operator in the limits $\xi\rightarrow\pm\infty$. In particular, we let $L_{h,\infty}$ be the operator defined by
\begin{equation}\begin{array}{lcl}L_{h,\infty} &=&\left(\begin{array}{ll}c_h\frac{d}{d\xi}-\Delta_h -\lim\limits_{\xi\rightarrow\infty}g_u(\overline{u}_h(\xi))& 1\\ -\rho & c_h\frac{d}{d\xi}+\gamma\rho\end{array}\right)\\[0.4cm]
&=&\left(\begin{array}{ll}c_h\frac{d}{d\xi}-\Delta_h -g_u(0)& 1\\ -\rho & c_h\frac{d}{d\xi}+\gamma\rho\end{array}\right).\end{array}\end{equation}
This system has constant coefficients. For $\lambda\in\C$ we introduce the notation
\begin{equation}
\begin{array}{lcl}
L_{h,\infty;\lambda}&=&L_{h,\infty}+\lambda.
\end{array}
\end{equation} We show that for $\lambda$ in a suitable right half-plane the operator $L_{h,\infty;\lambda}$ is hyperbolic in the sense of \cite{Faye2014,MPA}, i.e. we write
\begin{equation}\label{definitiedeltalh}\begin{array}{lcl}&&\Delta_{L_{h,\infty;\lambda}}(z)\\&=&\Big[L_{h,\infty;\lambda}e^{z\xi}\Big](0)\\[0.4cm]
&=&\left(\begin{array}{ll}c_hz-\frac{1}{h^2}\Big[\sum\limits_{k>0}\alpha_k\Big(e^{khz}+e^{-khz}-2\Big)\Big] -g_u(0)+\lambda& 1\\ -\rho & c_hz+\gamma\rho+\lambda\end{array}\right)\end{array}\end{equation}
and show that $\det(\Delta_{L_{h,\infty;\lambda}}(iy))\neq 0$ for all $y\in\R$. In the terminology of \cite{Faye2014,MPA}, this means that $L_h+\lambda$ is asymptotically hyperbolic. This allows us to compute the Fredholm index of $L_h+\lambda$.

\begin{remark} From this section onward we assume that (\asref{extraaannames}{\text{H}\alpha}) is satisfied. This is done for technical reasons, allowing us to apply the results from \cite{Faye2014}. In particular, this condition implies that the function $\Delta_{L_{h,\infty;\lambda}}(z)$ defined in (\ref{definitiedeltalh}) is well-defined in a vertical strip $|\mathrm{Re}\ (z)|< \nu$ .\end{remark}

\begin{lemma}\label{hyperbolic} Assume that (\asref{aannamespuls}{\text{H}}), (\asref{aannamesconstanten}{\text{H}}), (\asref{aannames}{\text{H}\alpha1}) and (\asref{extraaannames}{\text{H}\alpha}) are satisfied. There exists a constant $\tilde{\lambda}>0$ such that for all $0<h<\min\{h_*,\overline{h}\}$ and all $\lambda\in\C$ with $\Re \,\lambda>-\tilde{\lambda}$ the operator $L_{h,\infty;\lambda}$ is hyperbolic and
thus the operator $L_h+\lambda$ is asymptotically hyperbolic.\end{lemma}
\vspace*{4pt}\noindent\textit{Proof.} Remembering that $-g_u(0)=r_0>0$ and picking $y\in\R$, we compute
\begin{equation}\begin{array}{lcl}\Delta_{L_{h,\infty;\lambda}}(iy)&=&\left(\begin{array}{ll}c_hiy+\frac{1}{h^2}\Big[\sum\limits_{k>0}\alpha_k\Big(2-2\cos(khy)\Big)\Big]+ r_0+\lambda& 1\\ -\rho & c_hiy+\gamma\rho\end{array}\right)\\[0.4cm]
&=&\left(\begin{array}{ll}c_hiy+\frac{1}{h^2}A(hy)+ r_0+\lambda & 1\\ -\rho & c_hiy+\gamma\rho+\lambda\end{array}\right),\end{array}\end{equation}
where $A(hy)\geq 0$ is defined in (\asref{aannames}{\text{H}\alpha1}). We hence see
\begin{equation}\begin{array}{lcl}\det(\Delta_{L_{h,\infty;\lambda}}(iy))&=& \Big(c_hiy+\frac{1}{h^2}A(hy)+ r_0+\lambda\Big)\Big(c_hiy+\gamma\rho+\lambda\Big)+\rho.\end{array}\end{equation}
Let $\tilde{\lambda}=\frac{1}{4}\min\{\gamma\rho,r_0\}$ and assume that $\Re \,\lambda>-\tilde{\lambda}$. If $y\neq -\frac{\Im \,\lambda}{c_h}$ then we obtain
\begin{equation}\begin{array}{lcl}\Im \,\Big(\det(\Delta_{L_{h,\infty;\lambda}}(iy))\Big)&=&(c_hy+\Im \,\lambda)(\gamma\rho+\Re \,\lambda)\\[0.2cm]
&&\qquad +(\frac{1}{h^2}A(hy)+r_0+\Re \,\lambda)(c_hy+\Im \,\lambda) \\[0.2cm]
&=& (c_hy+\Im \,\lambda)(\gamma\rho+\frac{1}{h^2}A(hy)+r_0+2\Re \,\lambda)\\[0.2cm]
&\neq&0,\end{array}\end{equation}
since $\gamma\rho+\frac{1}{h^2}A(hy)+r_0+\Re \,\lambda>0$. For $y= -\frac{\Im \,\lambda}{c_h}$ we obtain
\begin{equation}\begin{array}{lcl}\Re \,\Big(\det(\Delta_{L_{h,\infty;\lambda}}( y))\Big)&=& \Big(\frac{1}{h^2}A(hy)+ r_0+\Re \,\lambda\Big)\Big(\gamma\rho+\Re \,\lambda\Big)+ \rho\\[0.2cm]
&>&\rho\\[0.2cm]
&> &0.\end{array}\end{equation}
In particular, we see that $\det(\Delta_{L_{h,\infty;\lambda}}(iy))\neq 0$ for all $y\in\R$, as desired.\qed

Before we consider the Fredholm properties of $L_h+\lambda$, we establish a technical estimate for the function $\Delta_{L_{h,\infty;\lambda}}$, which we need in \S \ref{sectiongreen} later on.
\begin{lemma}\label{hyperbolic2} Assume that (\asref{aannamespuls}{\text{H}}), (\asref{aannamesconstanten}{\text{H}}), (\asref{aannames}{\text{H}\alpha1}) and (\asref{extraaannames}{\text{H}\alpha}) are satisfied. Fix $0<h<\min\{h_*,\overline{h}\}$ and $S\subset \C$ compact in such a way that $\Re \,\lambda>-\tilde{\lambda}$ for all $\lambda\in S$. Then there exist constants $\kappa>0$ and $\Gamma>0$, possibly depending on $h$ and $S$, such that for all $z=x+iy\in\C$ with $|x|\leq \kappa$ and all $\lambda\in S$ we have the bound
\begin{equation}
\begin{array}{lcl}
|\det(\Delta_{L_{h,\infty;\lambda}}(z))|&\geq & \frac{1}{\Gamma}.
\end{array}
\end{equation} \end{lemma}
\vspace*{4pt}\noindent\textit{Proof.} Using assumption (\asref{extraaannames}{\text{H}\alpha}) we can pick $\kappa_1>0$ and $\Gamma_1>0$ in such a way that the bound
\begin{equation}
\begin{array}{lcl}
|\frac{1}{h^2}A(hz)|&:=&
\Big|\frac{1}{h^2}\Big[\sum\limits_{k>0}\alpha_k\Big(2-e^{khz}-e^{-khz}\Big)\Big]\Big|\\[0.2cm]
&\leq & \frac{1}{h^2}\sum\limits_{k>0}|\alpha_k|\Big(e^{hk|x|}+3\Big)\\[0.2cm]
&\leq &\Gamma_1
\end{array}
\end{equation}
holds for all $z=x+iy\in\C$ with $|x|\leq \kappa_1$.\par

Observe that for $z=x+iy\in\C$ and $\lambda\in S$ we have
\begin{equation}\begin{array}{lcl}\Re \,\Big(\det(\Delta_{L_{h,\infty;\lambda}}( z))\Big)&=& \Big(c_hx+\frac{1}{h^2}\Re \, A(hz)+ r_0+\Re \,\lambda\Big)\Big(c_hx+\gamma\rho+\Re \,\lambda\Big)\\[0.2cm]
&&\qquad -(c_hy+\Im \,\lambda)^2
-(c_hy+\Im \,\lambda)\frac{1}{h^2}(\Im \, A(y)) + \rho
.\end{array}\end{equation}
Since $S$ is compact we can find $Y>0$ such that for all $z=x+iy\in\C$ with $|y|\geq Y$ and $|x|\leq k_1$ and all $\lambda \in S$ we have
\begin{equation}\begin{array}{lcl}\Big|\Re \,\Big(\det(\Delta_{L_{h,\infty;\lambda}}( z))\Big)\Big|&\geq & \frac{1}{2}c_h^2 y^2\\[0.2cm]
&\geq & \frac{1}{2}c_h^2 Y^2
.\end{array}\end{equation}

Seeking a contradiction, let us assume that for each $0<\kappa\leq \kappa_1$ and each $\Gamma>0$ there exist $\lambda\in S$ and $z=x+iy\in\C$ with $|x|\leq\kappa$ and $|y|\leq Y$ for which
\begin{equation}\begin{array}{lcl}|\det(\Delta_{L_{h,\infty;\lambda}}( z))|&<&\frac{1}{\Gamma}.\end{array}\end{equation}
Then we can construct a sequence $\{\kappa_n,z_n,\lambda_n\}$ with $0<\kappa_n\leq \kappa_1$ for each $n$, $\kappa_n\rightarrow 0$, $\lambda_n\in S$ for each $n$ and $z_n=x_n+iy_n\in\C$ with $|x_n|\leq \kappa_n$ and $|y_n|\leq Y$ in such a way that $|\det(\Delta_{L_{h,\infty;\lambda_n}}( z_n))|<\frac{1}{n}$ for each $n$. By taking a subsequence if necessary we see that $\lambda_n\rightarrow \lambda$ for some $\lambda\in S$ and $z_n\rightarrow iy$ for some $y\in\R$ with $|y|\leq Y$. Since $\det(\Delta_{L_{h,\infty;\lambda}}( z))$ is continuous as a function of $\lambda$ and $z$, it follows that
\begin{equation}
\begin{array}{lcl}
\det(\Delta_{L_{h,\infty;\lambda}}( iy))&=&\lim\limits_{n\rightarrow\infty}\det(\Delta_{L_{h,\infty;\lambda_n}}( z_n))\\[0.2cm]
&=&0,
\end{array}
\end{equation}
which contradicts Lemma \ref{hyperbolic}. Hence, we can find $\kappa>0$ and $\Gamma>0$ as desired.\qed

\vspace*{4pt}\noindent\textit{Proof of Proposition \ref{LhlambdaFredholm}.} We have already seen in Lemma \ref{hyperbolic} that $L_h+\lambda$ is asymptotically hyperbolic in the sense of \cite{Faye2014,MPA}. Now according to \cite[Theorem 1.6]{Faye2014}, we obtain that $L_h+\lambda$ is a Fredholm operator and that the following identities hold
\begin{equation}\label{dimensiesfredholm}\begin{array}{lcl}\dim\Big(\ker (L_h+\lambda)\Big)&=&\codim\Big(\Range (L_h^*+\overline{\lambda})\Big),\\[0.2cm]
\codim\Big(\Range (L_h+\lambda)\Big)&=&\dim\Big(\ker (L_h^*+\overline{\lambda})\Big),\\[0.2cm]
\text{ind}(L_h+\lambda)&=&-\text{ind}(L_h^*+\overline{\lambda}),\end{array}\end{equation}
where
\begin{equation}\label{fredholmindex}\begin{array}{lcl}\text{ind}(L_h+\lambda)&=&\dim\Big(\ker(L_h+\lambda)\Big)
-\codim\Big(\Range(L_h+\lambda)\Big)\end{array}\end{equation}
is the Fredholm index of $L_h+\lambda$.\par

We follow the proof of \cite[Theorem B]{MPA}. For $0\leq\vartheta\leq 1$, we let the operator $L^\vartheta(h)$ be defined by
\begin{equation}\begin{array}{lcl}L^\vartheta(h)&=&(1-\vartheta)(L_h+\lambda)+\vartheta (L_{h,\infty}+\lambda).\end{array}\end{equation}
Note that the operator $L^\vartheta(h)$ is asymptotically hyperbolic for each $\vartheta$ and thus \cite[Theorem 1.6]{Faye2014} implies that these operators $L^\vartheta(h)$ are Fredholm. Moreover, the family $L^\vartheta(h)$ varies continuously with $\vartheta$ in $\mathcal{B}(\textbf{H}^1, \textbf{L}^2)$, which means the Fredholm index is constant. In particular, we see that
\begin{equation}\begin{array}{lcl}\text{ind}(L_h+\lambda)&=&\text{ind}(L_{h,\infty}+\lambda)\\[0.2cm]
&=&0,\end{array}\end{equation}
where the last equality follows from \cite[Theorem 1.7]{Faye2014}.\qed

We can now concentrate on the kernel of $L_h$. The goal is to exclude kernel elements other than $\Phi_h^+$. In order to accomplish this, we construct a quasi-inverse for $L_h$ by mimicking the approach of \cite[Proposition 3.2]{HJHBDF}. As a preparation, we obtain the following technical result.
\begin{lemma}\label{nietnul} Assume that (\asref{aannamespuls}{\text{H}}), (\asref{aannamesconstanten}{\text{H}}) and (\asref{aannames}{\text{H}\alpha1}) are satisfied. Recall the constant $\delta_0$ from Lemma \ref{eigenschappenL0}. Let $0<\lambda<\min\{\frac{1}{2},\delta_0\}$ be given. Then there exist constants $0<h_1^*\leq \min\{h_*,\overline{h}\}$ and $\kappa>0$ such that for all $0<h\leq h_1^*$ we have
\begin{equation}\begin{array}{lcl}\ip{\Phi_0^-,(L_h+\lambda)^{-1}\Phi_0^+}&>&\frac{1}{2}\lambda^{-1}\ip{\Phi_0^-,\Phi_h^+}\\[0.2cm]
&>&\frac{1}{2}\lambda^{-1}\kappa\\[0.2cm]
&>&0.\end{array}\end{equation}\end{lemma}
\vspace*{4pt}\noindent\textit{Proof.} We know from Lemma \ref{eigenschappenL0} that $\ip{\Phi_0^-,\Phi_0^+}>0$. Since $\Phi_h^+$ converges to $\Phi_0^+$ in $\textbf{L}^2$, it follows that $\ip{\Phi_0^-,\Phi_h^+}$ converges to $\ip{\Phi_0^-,\Phi_0^+}>0$. Fix $h_1^*\leq \min\{h_*,\overline{h},h_0'(\lambda)\}$ in such a way that
\begin{equation}\begin{array}{lcl}\nrm{\Phi_0^+-\Phi_h^+}_{\textbf{L}^2}&<&\frac{1}{2}\frac{\ip{\Phi_0^-,\Phi_h^+}}{2C_{\mathrm{unif}}}\end{array}\end{equation}
holds for all $0\leq h\leq h_1^*$, where
\begin{equation}\label{definitiecunif}\begin{array}{lcl}C_{\mathrm{unif}}&=&4C_0'\end{array}\end{equation}
and $C_0'$ is defined in Proposition \ref{equivalenttheorem4}. The factor $4$ in the definition is for technical reasons in a later proof. We assume from now on that $0<h\leq h_1^*$. Using $L_h\Phi_h^+=0$ we readily see
\begin{equation}\begin{array}{lcl}(L_h+\lambda)^{-1}\Phi_h^+&=&\lambda^{-1}\Phi_h^+.\end{array}\end{equation}
Recall that $\nrm{\Phi_0^-}_{\textbf{L}^2}=1$. Since $1<\lambda^{-1}$, we may use Proposition \ref{equivalenttheorem4} to obtain
\begin{equation}\begin{array}{lcl}\nrm{(L_h+\lambda)^{-1}\Phi_0^+-\lambda^{-1}\Phi_h^+}_{\textbf{L}^2}&=&
\nrm{(L_h+\lambda)^{-1}[\Phi_0^+-\Phi_h^+]}_{\textbf{L}^2}\\[0.2cm]
&\leq&C_{\mathrm{unif}}\Big[\nrm{\Phi_h^+-\Phi_0^+}_{\textbf{L}^2}+\lambda^{-1}|\ip{\Phi_h^+-\Phi_0^+,\Phi_0^-}|\Big]\\[0.2cm]
&< & C_{\mathrm{unif}}\lambda^{-1}\nrm{\Phi_0^+-\Phi_h^+}_{\textbf{L}^2}\Big(1+\nrm{\Phi_0^-}_{\textbf{L}^2}\Big)\\[0.2cm]
&=&2C_{\mathrm{unif}}\lambda^{-1}\nrm{\Phi_0^+-\Phi_h^+}_{\textbf{L}^2}.\end{array}\end{equation}
Remembering $\ip{\Phi_0^-,\Phi_h^+}>0$ and using Cauchy-Schwarz, we see that
\begin{equation}\label{tegenspraak}\begin{array}{lcl}|\ip{\frac{\Phi_0^-}{\ip{\Phi_0^-,\Phi_h^+}},(L_h+\lambda)^{-1}\Phi_0^+}-\lambda^{-1}|&=&|\ip{\frac{\Phi_0^-}{\ip{\Phi_0^-,\Phi_h^+}},(L_h+\lambda)^{-1}\Phi_0^+-\lambda^{-1}\Phi_h^+}|\\[0.2cm]
&<& \frac{\nrm{\Phi_0^-}_{\textbf{L}^2}}{\ip{\Phi_0^-,\Phi_h^+}}2C_{\mathrm{unif}}\lambda^{-1}\nrm{\Phi_0^+-\Phi_h^+}_{\textbf{L}^2}\\[0.2cm]
&\leq &\frac{1}{\ip{\Phi_0^-,\Phi_h^+}}2C_{\mathrm{unif}}\lambda^{-1}\frac{1}{2}\frac{\ip{\Phi_0^-,\Phi_h^+}}{2C_{\mathrm{unif}}}\\[0.2cm]
&=&\frac{1}{2}\lambda^{-1}.\end{array}\end{equation}
Hence, we must have
\begin{equation}\begin{array}{lcl}\ip{\Phi_0^-,(L_h+\lambda)^{-1}\Phi_0^+}&>&\frac{1}{2}\lambda^{-1}\ip{\Phi_0^-,\Phi_h^+}>0.\end{array}\end{equation} \qed

\begin{lemma}\label{prop3.2HJHBDF} Assume that (\asref{aannamespuls}{\text{H}}), (\asref{aannamesconstanten}{\text{H}}) and (\asref{aannames}{\text{H}\alpha1}) are satisfied. There exists $0<h_{**}\leq \min\{h_*,\overline{h}\}$ together with linear maps
\begin{equation}\begin{array}{lcl}\tilde{\gamma}_h^+:\textbf{L}^2&\rightarrow&\R\\[0.2cm] \tilde{L}_h^{\mathrm{qinv}}:\textbf{L}^2&\rightarrow &\textbf{H}^1,\end{array}\end{equation}
defined for all $0<h<h_{**}$, such that for all $\Theta\in \textbf{L}^2$ the pair
\begin{equation}\begin{array}{lcl}(\gamma,\Psi)&=&(\tilde{\gamma}_h^+\Theta,\tilde{L}_h^{\mathrm{qinv}}\Theta)\end{array}\end{equation}
is the unique solution to the problem
\begin{equation}\label{HJHBDF1}\begin{array}{lcl}L_h\Psi&=&\Theta+\gamma\Phi_0^+\end{array}\end{equation}
that satisfies the normalisation condition
\begin{equation}\label{HJHBDF2}\begin{array}{lcl}\ip{\Phi_0^-,\Psi}&=&0.\end{array}\end{equation}
In addition, there exists $C>0$ such that for all $0<h<h_{**}$ and all $\Theta\in \textbf{L}^2$ we have the bound
\begin{equation}\begin{array}{lcl}|\tilde{\gamma}_h^+\Theta|+\nrm{\tilde{L}_h^{\mathrm{qinv}}\Theta}_{\textbf{H}^1}&\leq &C\nrm{\Theta}_{\textbf{L}^2}.\end{array}\end{equation}
\end{lemma}
\vspace*{4pt}\noindent\textit{Proof.} Fix $0<\lambda<\min\{\frac{1}{2},\delta_0\}$ and let $0<h\leq \min\{h_*,\overline{h},h_0'(\lambda)\}$ be given, where $h_0'(\lambda)$ is defined in Proposition \ref{equivalenttheorem4}. For now, all constants will not depend on our choice of $\lambda$. We define the set
\begin{equation}\begin{array}{lcl}Z^1&=&\{\Psi\in \textbf{H}^1:\ip{\Phi_0^-,\Psi}=0\}.\end{array}\end{equation}
Pick $\Theta\in \textbf{L}^2$. We look for a solution $(\gamma,\Psi)\in\R\times Z^1$ of the problem
\begin{equation}\label{3.149equivalent}\begin{array}{lcl}\Psi&=&(L_h +\lambda)^{-1}[\Theta+\gamma \Phi_0^++\lambda \Psi].\end{array}\end{equation}
By Lemma \ref{nietnul} we have $\ip{\Phi_0^-,(L_h+\lambda)^{-1}\Phi_0^+}\neq 0$. Hence, for given $\Theta\in \textbf{L}^2,\Psi\in Z^1,h,\lambda$, we may write
\begin{equation}\begin{array}{lcl}\gamma(\Psi,\Theta,h,\lambda)&=&-\frac{\ip{\Phi_0^-,(L_h+\lambda)^{-1}(\Theta+\lambda \Psi)}}{\ip{\Phi_0^-,(L_h+\lambda)^{-1}\Phi_0^+}},\end{array}\end{equation}
which is the unique value for $\gamma$ for which
\begin{equation}\begin{array}{lcl}(L_h+\lambda)^{-1}[\Theta+\gamma \Phi_0^++\lambda \Psi]&\in &Z^1.\end{array}\end{equation}
Recall the constant $C_{\mathrm{unif}}$ from (\ref{definitiecunif}). With Proposition \ref{equivalenttheorem4} we obtain
\begin{equation}\begin{array}{lcl}|\ip{\Phi_0^-,(L_h+\lambda)^{-1}(\Theta+\lambda \Psi)}|&\leq &\nrm{\Phi_0^-}_{\textbf{L}^2} C_{\mathrm{unif}}\Big[\nrm{\Theta+\lambda \Psi}_{\textbf{L}^2}+\frac{1}{\lambda}|\ip{\Theta+\lambda \Psi,\Phi_0^-}|\Big]\\[0.2cm]
&\leq &\nrm{\Phi_0^-}_{\textbf{L}^2} C_{\mathrm{unif}}\Big[(1+\frac{1}{\lambda})\nrm{\Theta}_{\textbf{L}^2}+\lambda\nrm{\Psi}_{\textbf{L}^2}\Big]\\[0.2cm]
&\leq &C_1\Big[\lambda^{-1}\nrm{\Theta}_{\textbf{L}^2}
+\lambda\nrm{\Psi}_{\textbf{L}^2}\Big]\end{array}\end{equation}
for some $C_1$ that is independent of $h,\lambda$. Here we used that $\lambda<1$ and thus $1+\frac{1}{\lambda}<\frac{2}{\lambda}$. Exploiting $\lambda<\frac{1}{2}$ and applying Lemma \ref{nietnul}, we see that
\begin{equation}\begin{array}{lcl}|\gamma(\Psi,\Theta,h,\lambda)|&=&|\ip{\Phi_0^-,(L_h+\lambda)^{-1}(\Theta+\lambda \Psi)}|\frac{1}{|\ip{\Phi_0^-,(L_h+\lambda)^{-1}\Phi_0^+}|}\\[0.2cm]
&\leq &C_1\Big[\lambda^{-1}\nrm{\Theta}_{\textbf{L}^2}+\lambda\nrm{\Psi}_{\textbf{L}^2}\Big]\frac{1}{\frac{1}{2}\lambda^{-1}\ip{\Phi_0^-,\Phi_h^+}}\\[0.2cm]
&\leq &C_1\Big[\kappa\nrm{\Theta}_{\textbf{L}^2}+\kappa\lambda^2\nrm{\Psi}_{\textbf{L}^2}\Big]\\[0.2cm]
&\leq  &C_2\Big[\nrm{\Theta}_{\textbf{L}^2}+\lambda^2\nrm{\Psi}_{\textbf{L}^2}\Big].\end{array}\end{equation}
Here we used that $\ip{\Phi_0^-,\Phi_h^+}$ converges to $\ip{\Phi_0^-,\Phi_0^+}>0$, which means that $\ip{\Phi_0^-,\Phi_h^+}$ can be bounded away from zero. For $\Psi\in Z^1$ we write
\begin{equation}\begin{array}{lcl}t(\Psi)&=&\Theta+\gamma(\Psi,\Theta,h,\lambda) \Phi_0^++\lambda \Psi\end{array}\end{equation}
and
\begin{equation}\begin{array}{lcl}T(\Psi)&=&(L_h +\lambda)^{-1}t(\Psi).\end{array}\end{equation}
For $\Psi\in Z^1$ Proposition \ref{equivalenttheorem4} implies
\begin{equation}\begin{array}{lcl}\nrm{T(\Psi)}_{\textbf{H}^1}&\leq &C_{\mathrm{unif}}\Big[\nrm{\Theta+\gamma(\Psi,\Theta,h,\lambda)\Phi_0^++\lambda \Psi}_{\textbf{L}^2}\\[0.2cm]
&&\qquad +\frac{1}{\lambda}|\ip{\Theta+\gamma(\Psi,\Theta,h,\lambda)\Phi_0^++\lambda \Psi,\Phi_0^-}|\Big]\\[0.2cm]
&\leq &C_3\Big[\frac{1}{\lambda}\nrm{\Theta}_{\textbf{L}^2}+\lambda\nrm{\Psi}_{\textbf{L}^2}\Big]\\[0.2cm]
&\leq &C_3\Big[\frac{1}{\lambda}\nrm{\Theta}_{\textbf{L}^2}+\lambda\nrm{\Psi}_{\textbf{H}^1}\Big].\end{array}\end{equation}
For $\Psi_1,\Psi_2\in Z^1$, a second application of Proposition \ref{equivalenttheorem4} yields
\begin{equation}\begin{array}{lcl}|\gamma(\Psi_1,\Theta,h,\lambda)-\gamma(\Psi_2,\Theta,h,\lambda)|&=&\left|\frac{\ip{\Phi_0^-,(L_h+\lambda)^{-1}(\lambda \Psi_1-\lambda \Psi_2)}}{\ip{\Phi_0^-,(L_h+\lambda)^{-1}\Phi_0^+}}\right|\\[0.2cm]
&\leq & \frac{1}{\ip{\Phi_0^-,(L_h+\lambda)^{-1}\Phi_0^+}}C_{\mathrm{unif}}\Big[\lambda\nrm{\Psi_1-\Psi_2}_{\textbf{L}^2}\\[0.2cm]
&&\qquad +\frac{1}{\lambda}|\ip{\lambda \Psi_1-\lambda \Psi_2,\Phi_0^-}|\Big]\\[0.2cm]
&\leq& C_4\lambda\Big[\lambda\nrm{\Psi_1-\Psi_2}_{\textbf{L}^2}+0\Big]\\[0.2cm]
&\leq &C_4 \lambda^2 \nrm{\Psi_1-\Psi_2}_{\textbf{H}^1}.\end{array}\end{equation}
Applying Proposition \ref{equivalenttheorem4} for the final time, we see
\begin{equation}\begin{array}{lcl}\nrm{T(\Psi_1)-T(\Psi_2)}_{\textbf{H}^1}&\leq &C_{\mathrm{unif}}\Big[\nrm{t(\Psi_1)-t(\Psi_2)}_{\textbf{L}^2} +\frac{1}{\lambda}|\ip{t(\Psi_1)-t(\Psi_2),\Phi_0^-}|\Big]\\[0.2cm]
&\leq & C_{\mathrm{unif}}\Big[\nrm{t(\Psi_1)-t(\Psi_2)}_{\textbf{L}^2}\\[0.2cm]
&&\qquad+\frac{1}{\lambda}\ip{\big(\gamma(\Psi_1,\Theta,h,\lambda)-\gamma(\Psi_2,\Theta,h,\lambda)\big)\Phi_0^+,\Phi_0^-}\\[0.2cm]
&&\qquad +\frac{1}{\lambda}\ip{\lambda(\Psi_1-\Psi_2),\Phi_0^-}\Big]\\[0.2cm]
&\leq & C_{\mathrm{unif}}\Big[\nrm{t(\Psi_1)-t(\Psi_2)}_{\textbf{L}^2}+\frac{1}{\lambda}\Big(C_4\lambda^2\nrm{\Psi_1-\Psi_2}_{\textbf{H}^1}+0\Big)\Big]\\[0.2cm]
&\leq &C_{\mathrm{unif}}C_4\lambda^2 \nrm{\Psi_1-\Psi_2}_{\textbf{H}^1}+C_{\mathrm{unif}}\lambda\nrm{\Psi_1-\Psi_2}_{\textbf{H}^1}\\[0.2cm]
&&\qquad +C_4\lambda \nrm{\Psi_1-\Psi_2}_{\textbf{H}^1}\\[0.2cm]
&\leq &C_5\lambda \nrm{\Psi_1-\Psi_2}_{\textbf{H}^1}.
\end{array}\end{equation}
In view of these bounds, we pick $\lambda$ to be small enough to have $C_3\lambda<\frac{1}{2}$ and $C_5\lambda<\frac{1}{2}$. Since this $\lambda$ is now fixed, we can allow the constants in the final part of the proof to depend on $\lambda$. In addition, we write $h_{**}=\min\{h_1^*,h_0'(\lambda)\}$ and pick $0<h< h_{**}$. Then $T:Z^1\rightarrow Z^1$ is a contraction, so the fixed point theorem implies that there is a unique $\tilde{L}_h^{\mathrm{qinv}}(\Theta)\in Z^1$ for which
\begin{equation}\begin{array}{lcl}\tilde{L}_h^{\mathrm{qinv}}(\Theta)&=&(L_h+\lambda)^{-1}\Big[\Theta+\gamma(\tilde{L}_h^{\mathrm{qinv}}(\Theta),\Theta,h,\lambda)\Phi_0^++\lambda \tilde{L}_h^{\mathrm{qinv}}(\Theta)\Big].\end{array}\end{equation}
Furthermore, we have
\begin{equation}\begin{array}{lcl}\frac{1}{2}\nrm{\tilde{L}_h^{\mathrm{qinv}}(\Theta)}_{\textbf{H}^1}&\leq &(1-\lambda C_3)\nrm{\tilde{L}_h^{\mathrm{qinv}}(\Theta)}_{\textbf{H}^1}\\[0.2cm]
&\leq &C_3\lambda^{-1}\nrm{\Theta}_{\textbf{L}^2}\\[0.2cm]
&\leq &C_6\nrm{\Theta}_{\textbf{L}^2}.\end{array}\end{equation}
Writing $\tilde{\gamma}_h^+(\Theta)=\gamma(\tilde{L}_h^{\mathrm{qinv}}(\Theta),\Theta,h,\lambda)$, we compute
\begin{equation}\begin{array}{lcl}|\tilde{\gamma}_h^+(\Theta)|&\leq &C_2[\nrm{\Theta}_{\textbf{L}^2}+\lambda^2\nrm{\Theta}_{\textbf{L}^2}]\\[0.2cm]
&\leq &C_7\nrm{\Theta}_{\textbf{L}^2}.\end{array}\end{equation}
Finally we see that (\ref{3.149equivalent}) is in fact equivalent to (\ref{HJHBDF1})-(\ref{HJHBDF2}), so in fact $\tilde{L}_h^{\mathrm{qinv}}(\Theta)$ and $\tilde{\gamma}_h^+(\Theta)$ do not depend on $\lambda$. In addition, the constants $h_{**}$, $C_6$ and $C_7$ above only depend on the one fixed $\lambda$ and, as such, do not depend on $h$ or $\Theta$.\qed

\begin{lemma}\label{FredholmeigenschappenLhklein2} Assume that (\asref{aannamespuls}{\text{H}}), (\asref{aannamesconstanten}{\text{H}}), (\asref{aannames}{\text{H}\alpha1}) and (\asref{extraaannames}{\text{H}\alpha}) are satisfied. Let $0<h< h_{**}$ be given. Then we have the inclusion
\begin{equation}\begin{array}{lcl}\span\{\Phi_h^+\}&\subset&\ker (L_h)\\[0.2cm]
&=&\{\Psi\in \textbf{L}^2:\ip{\Psi,\Theta}=0\text{ for all }\Theta\in \Range (L_h^*)\},\end{array}\end{equation}
where $L_h^*$ is the formal adjoint of $L_h$.\end{lemma}
\vspace*{4pt}\noindent\textit{Proof.} By differentiating the differential equation (\ref{ditishetprobleem}) we see that $L_h\Phi_h^+=0$. We know that $(\overline{u}_h,\overline{w}_h)-(\overline{u}_0,\overline{w}_0)\rightarrow 0\in \textbf{H}^1$. Since $(\overline{u}_0',\overline{w}_0')$ decays exponentially, we get $(\overline{u}_0',\overline{w}_0')\in\textbf{L}^2$. Hence, we can assume that $h_{**}$ is small enough such that $\Phi_h^+\in\textbf{L}^2$ for all $0<h< h_{**}$. Since $L_h\Phi_h^+=0$ we obtain from the differential equation that also $(\Phi_h^+)'\in\textbf{L}^2$. In particular, we see that $\Phi_h^+\in\textbf{H}^1$ and hence $\Phi_h^+\in\ker(L_h)$. \qed

\begin{lemma}\label{FredholmeigenschappenLh} Assume that (\asref{aannamespuls}{\text{H}}), (\asref{aannamesconstanten}{\text{H}}), (\asref{aannames}{\text{H}\alpha1}) and (\asref{extraaannames}{\text{H}\alpha}) are satisfied. Let $0<h< h_{**}$ be given. Then we have
\begin{equation}\begin{array}{lcl}\ker (L_h)&=&\span\{\Phi_h^+\}\\[0.2cm]
&=&\{\Psi\in \textbf{L}^2:\ip{\Psi,\Theta}=0\text{ for all }\Theta\in \Range (L_h^*)\},\end{array}\end{equation}
where $L_h^*$ is the formal adjoint of $L_h$.\end{lemma}
\noindent\textit{Proof.} We show that $\dim(\ker(L_h))=1$. Since $\Phi_h^+\in\ker(L_h)$, we assume that there exists $\Psi\in\ker(L_h)$ in such a way that $\Psi$ is not a scalar multiple of $\Phi_h^+$.\par

Suppose first that $\ip{\Psi,\Phi_0^-}=0$. Then Lemma \ref{prop3.2HJHBDF} yieds by linearity of $\tilde{L}_h^{\text{qinv}}$ that
\begin{equation}
\begin{array}{lcl}
\Psi&=&\tilde{L}_h^{\text{qinv}}[0]\\[0.2cm]
&=&0,
\end{array}
\end{equation}
which gives a contradiction. Hence, we suppose that $\ip{\Psi,\Phi_0^-}\neq 0$. In the proof of Lemma \ref{nietnul} we saw that $\ip{\Phi_h^+,\Phi_0^-}\neq 0$. As such, we can pick $a,b\in\R\setminus\{0\}$ in such a way that
\begin{equation}
\begin{array}{lcl}
\ip{a\Phi_h^++b\Psi,\Phi_0^-}&=&0.
\end{array}
\end{equation}
Again it follows from Lemma \ref{prop3.2HJHBDF} that $a\Phi_h^++b\Psi=0$ which gives a contradiction. Therefore,
such a kernel element $\Psi$ does not exist. \par

Since we already know that $\Phi_h^+\in \ker(L_h)$, we must have $\dim\big(\ker(L_h)\big)=1$, which completes the proof.\qed

The remaining major goal of this section is to find a family of elements $\Phi_h^-\in\ker(L_h^*)$ which satisfies $\Phi_h^-\rightarrow\Phi_0^-$ as $h\downarrow 0$. To establish this, we repeat part of the process above for the adjoint operator $L_h^*$. The key difference is that we must construct the family $\Phi_h^-$ by hand. This requires a significant refinement of the argument used above to characterize $\ker(L_h^*)$.\par

First we need a technical result, similar to Lemma \ref{nietnul}.
\begin{lemma}\label{nietnul2} Assume that (\asref{aannamespuls}{\text{H}}), (\asref{aannamesconstanten}{\text{H}}) and (\asref{aannames}{\text{H}\alpha1}) are satisfied. Fix $0<\lambda<\frac{1}{2}$ and $0< h\leq \min\{h_{**},h_0'(\lambda)\}$, where $h_0'(\lambda)$ is defined in Proposition \ref{equivalenttheorem4}. Then we have
\begin{equation}\begin{array}{lcl}\ip{\Phi_0^+,(L_h^*+\lambda)^{-1}\Phi_0^-}&>&\frac{\ip{\Phi_0^+,\Phi_0^-}}{2}\lambda^{-1}.\end{array}\end{equation}
\end{lemma}
\vspace*{4pt}\noindent\textit{Proof.} Lemma \ref{eigenschappenL0} implies that $\ip{\Phi_0^+,\Phi_0^-}>0$. Remembering that
\begin{equation}\begin{array}{lcl}L_h^*-L_0^*&=&\left(\begin{array}{ll} (c_0-c_h)\frac{d}{d\xi}-(\Delta_h-\frac{d^2}{dx^2})+(g_u(\overline{u}_0)-g_u(\overline{u}_h))&0\\ 0 & (c_0-c_h)\frac{d}{d\xi}\end{array}\right)\end{array}\end{equation}
and that $L_0^*\Phi_0^-=0$, we obtain
\begin{equation}\begin{array}{lcl}&& (L_h^*+\lambda)\Big[(L_h^*+\lambda)^{-1}\Phi_0^--(L_0^*+\lambda)^{-1}\Phi_0^-\Big]\\[0.2cm]
&=&\Phi_0^--\Phi_0^-+(L_h^*-L_0^*)(L_0^*+\lambda)^{-1}\Phi_0^-\\[0.2cm]
&=&(L_h^*-L_0^*)\lambda^{-1}\Phi_0^-.\end{array}\end{equation}
Recall the constant $C_{\mathrm{unif}}$ from (\ref{definitiecunif}). Proposition \ref{equivalenttheorem4} yields
\begin{equation}\begin{array}{lcl}\nrm{(L_h^*+\lambda)^{-1}\Phi_0^--(L_0^*+\lambda)^{-1}\Phi_0^-}_{\textbf{L}^2}&\leq& C_{\mathrm{unif}}\Big[\nrm{(L_h^*-L_0^*)\lambda^{-1}\Phi_0^-}_{\textbf{L}^2}\\[0.2cm]
&&\qquad +|\ip{(L_h^*-L_0^*)\lambda^{-1}\Phi_0^-,\Phi_0^+}|\Big]\\[0.2cm]
&\leq &C_{\mathrm{unif}}(1+\lambda^{-1})\nrm{(L_h^*-L_0^*)\lambda^{-1}\Phi_0^-}_{\textbf{L}^2}.\end{array}\end{equation}
Using Lemma \ref{eigenschappenDeltah} and the fact that $c_h$ converges to $c_0$ and $g_u(\overline{u}_h)$ to $g_u(\overline{u}_0)$, it follows that
\begin{equation}\begin{array}{lcl}C_{\mathrm{unif}}(1+\lambda^{-1})\nrm{(L_h^*-L_0^*)\lambda^{-1}\Phi_0^-}_{\textbf{L}^2}&\rightarrow & 0\end{array}\end{equation}
as $h\downarrow 0$. Possibly after decreasing $h_0'(\lambda)>0$, we hence see that
\begin{equation}\begin{array}{lcl}\ip{\Phi_0^+,(L_h^*+\lambda)^{-1}\Phi_0^-}&=&\ip{\Phi_0^+,(L_0^*+\lambda)^{-1}\Phi_0^-}\\[0.2cm]
&&\qquad +\ip{\Phi_0^+,(L_h^*+\lambda)^{-1}\Phi_0^--(L_0^*+\lambda)^{-1}\Phi_0^-}\\[0.2cm]
&=&\lambda^{-1}\ip{\Phi_0^+,\Phi_0^-}+\ip{\Phi_0^+,(L_h^*+\lambda)^{-1}\Phi_0^--(L_0^*+\lambda)^{-1}\Phi_0^-}\\[0.2cm]
&>&\frac{\ip{\Phi_0^+,\Phi_0^-}}{2}\lambda^{-1}\end{array}\end{equation}
holds for all $0<h<\min\{h_{**},h_0'(\lambda)\}$.
\qed

\begin{lemma}\label{lemmaphih-bewijs2} Assume that (\asref{aannamespuls}{\text{H}}), (\asref{aannamesconstanten}{\text{H}}), (\asref{aannames}{\text{H}\alpha1}) and (\asref{extraaannames}{\text{H}\alpha}) are satisfied. Fix $0< h< h_{**}$. There exist linear maps
\begin{equation}\begin{array}{lcl}\tilde{\gamma}_h^-:\textbf{L}^2&\rightarrow&\R,\\[0.2cm]
\tilde{L}_h^{*,\mathrm{qinv}}:\textbf{L}^2&\rightarrow & \textbf{H}^1\end{array}\end{equation}
such that for all $\Theta\in \textbf{L}^2$ the pair
\begin{equation}\begin{array}{lcl}(\gamma,\Psi)&=&(\tilde{\gamma}_h^-\Theta,\tilde{L}_h^{*,\mathrm{qinv}}\Theta)\end{array}\end{equation}
is the unique solution to the problem
\begin{equation}\begin{array}{lcl}L_h^*\Psi&=&\Theta+\gamma\Phi_0^-\end{array}\end{equation}
that satisfies the normalisation condition
\begin{equation}\begin{array}{lcl}\ip{\Phi_0^+,\Psi}&=&0.\end{array}\end{equation}
Furthermore, there exists $C^*>0$, such that for all $0<h<h_{**}$ and all $\Theta\in \textbf{L}^2$ we have the bound
\begin{equation}\begin{array}{lcl}|\tilde{\gamma}_h^-\Theta|+\nrm{\tilde{L}_h^{*,\mathrm{qinv}}\Theta}_{\textbf{H}^1}&\leq &C^*\nrm{\Theta}_{\textbf{L}^2}.\end{array}\end{equation}\end{lemma}
\vspace*{4pt}\noindent\textit{Proof.} We define the set
\begin{equation}\begin{array}{lcl} Z^1&=&\{\Psi\in \textbf{H}^1:\ip{\Phi_0^+,\Psi}=0\}.\end{array}\end{equation}
Pick $\Theta\in \textbf{L}^2$. We look for a solution $(\gamma,\Psi)\in\R\times Z^1$ of the problem
\begin{equation}\begin{array}{lcl}\Psi&=&(L_h^* +\lambda)^{-1}[\Theta+\gamma \Phi_0^-+\lambda \Psi].\end{array}\end{equation}
Lemma \ref{nietnul2} implies that $\ip{\Phi_0^+,(L_h^*+\lambda)^{-1}\Phi_0^-}\neq 0$. Hence, for given $\Theta\in \textbf{L}^2,\Psi\in Z^1,h,\lambda$, we may write
\begin{equation}\begin{array}{lcl}\gamma(\Psi,\Theta,h,\lambda)&=&-\frac{\ip{\Phi_0^+,(L_h^*+\lambda)^{-1}(\Theta+\lambda \Psi)}}{\ip{\Phi_0^+,(L_h^*+\lambda)^{-1}\Phi_0^-}},\end{array}\end{equation}
which is the unique value for $\gamma$ for which
\begin{equation}\begin{array}{lcl}(L_h^*+\lambda)^{-1}[\Theta+\gamma \Phi_0^-+\lambda \Psi]&\in &Z^1.\end{array}\end{equation}
From now on the proof is identical to that of Lemma \ref{prop3.2HJHBDF}, so we omit it.\qed

\begin{lemma}\label{lemmaphih-} Assume that (\asref{aannamespuls}{\text{H}}), (\asref{aannamesconstanten}{\text{H}}), (\asref{aannames}{\text{H}\alpha1}) and (\asref{extraaannames}{\text{H}\alpha}) are satisfied. For each $0<h< h_{**}$ there exists an element $\Phi_h^-\in \ker(L_h^*)$ such that the family $\Phi_h^-$ converges to $\Phi_0^-$ in $\textbf{H}^1$ as $h\downarrow 0$.\end{lemma}
\noindent\textit{Proof.} We repeat some of the steps of the proof of Lemma \ref{FredholmeigenschappenLh}, but now for $L_h^*$.\par

Fix $0<h< h_{**}$. Since $\dim(\ker(L_h^*))=1$ by Proposition \ref{LhlambdaFredholm} and Lemma \ref{FredholmeigenschappenLh}, we can pick $\Phi\in \ker(L_h^*)$ with $\Phi\neq 0$. If we would have $\ip{\Phi,\Phi_0^+}=0$, then we would obtain
\begin{equation}
\begin{array}{lcl}
0&=&L_h^{*,\mathrm{qinv}}[0]\\[0.2cm]
&=&\Phi,
\end{array}
\end{equation}
which leads to a contradiction. Hence, we can define the kernel element $\Phi_h^-$ of $L_h^*$ as follows: $\Phi_h^-$ is the unique kernel element of $L_h^*$ with $\ip{\Phi_h^-,\Phi_0^+}=\ip{\Phi_0^-,\Phi_0^+}$. Since we see that
\begin{equation}
\begin{array}{lcl}
\ip{\Phi_0^--\Phi_h^-,\Phi_0^+}&=&0,
\end{array}
\end{equation}
we obtain, upon defining
\begin{equation}
\begin{array}{lcl}
\Theta_h&:=&L_h^*\Phi_0^-,
\end{array}
\end{equation}
that
\begin{equation}
\begin{array}{lcl}
\Phi_0^--\Phi_h^-&=&L_h^{*,\mathrm{qinv}}[\Theta_h].
\end{array}
\end{equation}
Using Lemma \ref{lemmaphih-bewijs2}, we can estimate
\begin{equation}
\begin{array}{lcl}
\nrm{\Phi_0^--\Phi_h^-}_{\mathbf{H}^1}&=&\nrm{L_h^{*,\mathrm{qinv}}[\Theta_h]}_{\mathbf{H}^1}\\[0.2cm]
&\leq & C_-\nrm{\Theta_h}_{\mathbf{L}^2}.
\end{array}
\end{equation}
From the proof of Lemma \ref{nietnul2} we know that $\Theta_h\rightarrow 0$ as $h\downarrow 0$ in $\mathbf{L}^2$. Therefore, we see that $\Phi_h^-\rightarrow \Phi_0^-$ as $h\downarrow 0$ in $\mathbf{H}^1$. \qed

In the final part of this section we establish item (3) of Proposition \ref{samenvattingh5}. To this end, we recall the spaces
\begin{equation}\begin{array}{lcl}X_h&=&\{\Theta\in \textbf{H}^1: \ip{\Phi_h^-,\Theta}=0\}\end{array}\end{equation}
and
\begin{equation}\begin{array}{lcl}Y_h&=&\{\Theta\in \textbf{L}^2: \ip{\Phi_h^-,\Theta}=0\},\end{array}\end{equation}
together with the constant $C_{\mathrm{unif}}$ from (\ref{definitiecunif}).
\begin{lemma}\label{lemmaXh} Assume that (\asref{aannamespuls}{\text{H}}), (\asref{aannamesconstanten}{\text{H}}), (\asref{aannames}{\text{H}\alpha1}) and (\asref{extraaannames}{\text{H}\alpha}) are satisfied.
For each $0<h< h_{**}$ we have that $L_h:X_h\rightarrow Y_h$ is invertible and we have the uniform bound
\begin{equation}
\begin{array}{lcl}
\nrm{L_h^{-1}}\leq C_{\mathrm{unif}}.
\end{array}
\end{equation} \end{lemma}

\vspace*{4pt}\noindent\textit{Proof.} Fix $0<h< h_{**}$. Clearly $L_h:X_h\rightarrow Y_h$ is a bounded bijective linear map, so the Banach isomorphism theorem implies that $L_h^{-1}:Y_h\rightarrow X_h$ is bounded. Now let $\delta>0$ be a small constant such that $\delta C_{\mathrm{unif}}<1$. Without loss of generality we assume that $0<h_{**}\leq h_0'(\delta)$ and that $\nrm{\Phi_h^--\Phi_0^-}_{\textbf{H}^1}\leq\delta$ for all $0<h< h_{**}$. This is possible by Lemma \ref{lemmaphih-}.\par

Pick any $\Psi\in X_h$. Remembering that $\ip{\Psi,\Phi_h^-}=0$ and $\ip{L_h\Psi,\Phi_h^-}=0$, we obtain the estimate
\begin{equation}
\begin{array}{lcl}
\frac{1}{\delta}|\ip{(L_h+\delta)\Psi,\Phi_0^-}|&=&\frac{1}{\delta}|\ip{(L_h+\delta)\Psi,\Phi_0^--\Phi_h^-}|\\[0.2cm]
&\leq & \frac{1}{\delta}\nrm{(L_h+\delta)\Psi}_{\textbf{L}^2}\delta\\[0.2cm]
&\leq & \nrm{L_h\Psi}_{\textbf{L}^2}+\delta\nrm{\Psi}_{\textbf{H}^1}.
\end{array}
\end{equation}
Applying Proposition \ref{equivalenttheorem4}, we hence see
\begin{equation}\begin{array}{lcl}\nrm{\Psi}_{\textbf{H}^1}&\leq &\frac{1}{4}C_{\mathrm{unif}}[\nrm{(L_h+\delta)\Psi}_{\textbf{L}^2}+\frac{1}{\delta}|\ip{(L_h+\delta)\Psi,\Phi_0^-}|]\\[0.2cm]
&\leq & \frac{1}{4}C_{\mathrm{unif}}[2\nrm{L_h\Psi}_{\textbf{L}^2}+2\delta\nrm{\Psi}_{\textbf{H}^1}]\\[0.2cm]
&\leq &\frac{1}{2}C_{\mathrm{unif}}\nrm{L_h\Psi}_{\textbf{L}^2}+\frac{1}{2}\nrm{\Psi}_{\textbf{H}^1}.\end{array}\end{equation}
We, therefore, get the bound
\begin{equation}\begin{array}{lcl}\nrm{\Psi}_{\textbf{H}^1}&\leq &C_{\mathrm{unif}}\nrm{L_h\Psi}_{\textbf{L}^2},\end{array}\end{equation}
which yields the desired estimate $\nrm{L_h^{-1}}\leq C_{\mathrm{unif}}$.\qed
%
%

\vspace*{4pt}\noindent\textit{Proof of Proposition \ref{samenvattingh5}.} This result follows directly from Lemmas \ref{FredholmeigenschappenLh}, \ref{lemmaphih-} and \ref{lemmaXh}.\qed.\par


\section{The resolvent set}\label{spectralanalysis}

In this section we prove Theorem \ref{totalespectrum} by explicitly determining the spectrum of the operator $-L_h$ defined in (\ref{definitieLh}). Our approach hinges on the periodicity of this spectrum, which we describe in our first result.
\begin{lemma}\label{lemmaperiodiek} Assume that (\asref{aannamespuls}{\text{H}}), (\asref{aannamesconstanten}{\text{H}}), (\asref{aannames}{\text{H}\alpha1}) and (\asref{extraaannames}{\text{H}\alpha}) are satisfied. Fix $0<h< h_{**}$. Then the spectrum of $L_h$ is invariant under the operation $\lambda\mapsto \lambda+2\pi i c_h\frac{1}{h}$.\end{lemma}
In particular, we can restrict our attention to values with imaginary part in between $-\frac{\pi c_h}{h}$ and $\frac{\pi c_h}{h}$. We divide our `half-strip' into four regions and in each region we calculate the spectrum. Values close to $0$ (region $R_1$) will be treated in Proposition \ref{lemmaspectrumklein}; values with a large real part (region $R_2$) in Proposition \ref{spectrumgroot} and values with a large imaginary part (region $R_3$) in Proposition \ref{spectrumhoog}. In Corollary \ref{spectrumintermediate} we discuss the remaining intermediate subset (region $R_4$), which is compact and independent of $h$. The regions are illustrated in Figure \ref{plaatjespectrum} below.\par

\begin{figure}[h]
\begin{center}
\includegraphics[scale=0.18]{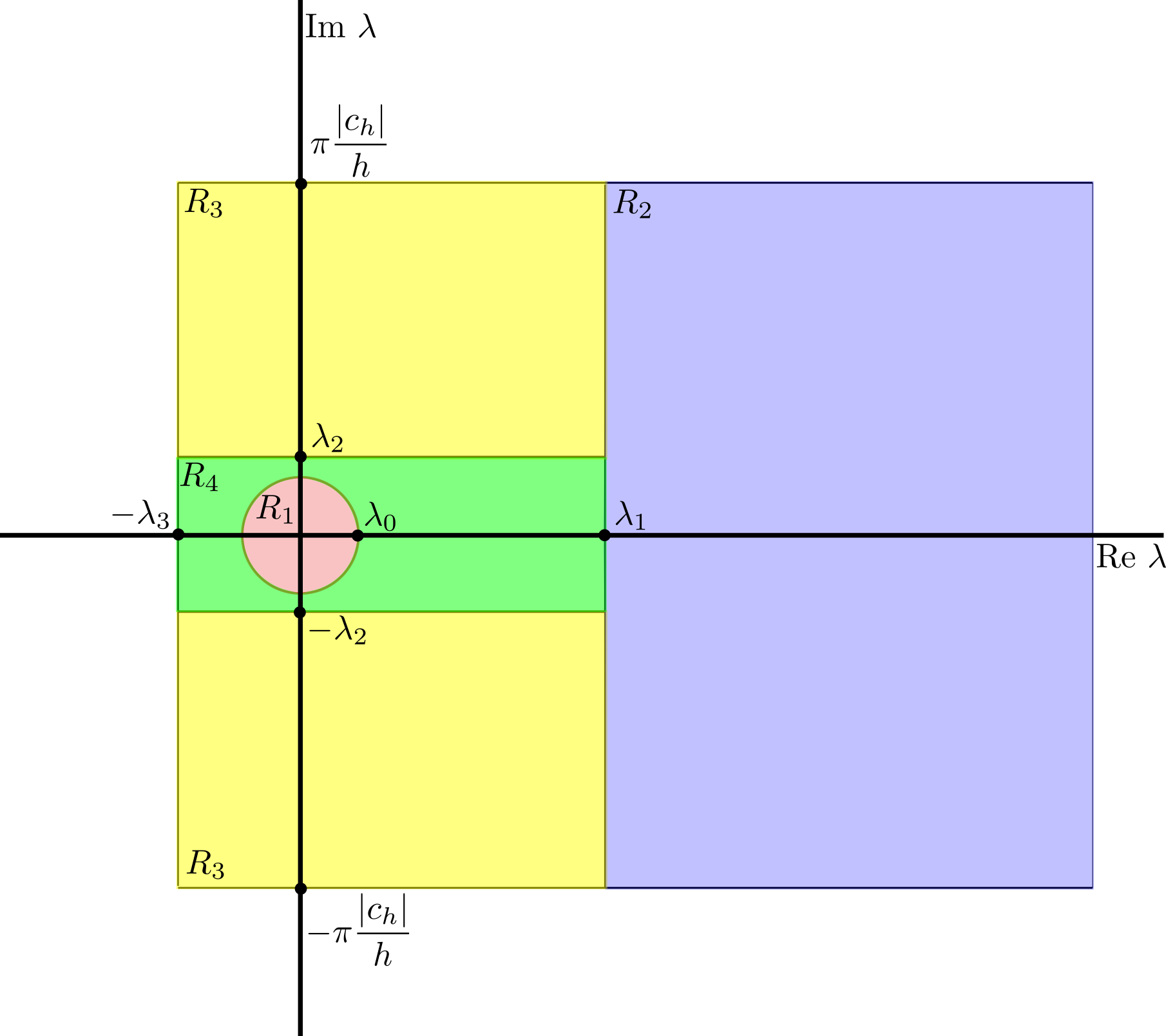}
\caption{Illustration of the regions $R_1,R_2,R_3$ and $R_4$. Note that the regions $R_2$ and $R_3$ grow when $h$ decreases, while the regions $R_1$ and $R_4$ are independent of $h$.}
\label{plaatjespectrum}
\end{center}
\end{figure}

From this section onward we need to assume that (\asref{extraaannamespuls}{\text{H}}) is satisfied. Indeed, this allows us to lift the invertibility of $L_0+\lambda$ to $L_h+\lambda$ simultaneously for all $\lambda$ in appropriate compact sets.\vspace*{4pt}

\noindent\textit{Proof of Lemma \ref{lemmaperiodiek}.} Fix $k\in\Z$ and write $p=2\pi i k\frac{1}{h}$. We define the exponential shift operator $e_\omega$ by
\begin{equation}\begin{array}{lcl}[e_\omega V](x)&=&e^{\omega x}V(x).\end{array}\end{equation}
For any $\lambda\in\C$, $\Psi=(\phi,\psi)\in \textbf{H}^1$ and $x\in\R$ we obtain
\begin{equation}\begin{array}{lcl}(e_{-p}\Delta_h e_{p})\phi(x)&=&e^{-px}\Delta_h(e_{p}\phi)(x)\\[0.2cm]
&=&\frac{1}{h^2}\sum\limits_{l>0}\alpha_l(e^{p lh}\phi(x+lh)+e^{-p lh}\phi(x-lh)-2\phi(x))\\[0.2cm]
&=&\frac{1}{h^2}\sum\limits_{l>0}\alpha_l(\phi(x+lh)+\phi(x-lh)-2\phi(x))\\[0.2cm]
&=&\Delta_h \phi(x),\end{array}\end{equation}
since $plh\in 2\pi i\Z$ for all $l>0$. In particular, we can compute
\begin{equation}\begin{array}{lcl}[e_{-p}(L_h-\lambda)e_{p}\Psi](x)&=&e^{-px}[(L_h-\lambda)e_{p}\Psi](x)\\[0.2cm]
&=&e^{-px}\Big(\begin{array}{l} c_h\frac{d}{d\xi}(e^{px}\phi(x))-\Delta_h(e_{p}\phi)(x)\\ -\rho e^{px}+c_h\frac{d}{d\xi}(e^{px}\psi(x))\end{array}\Big)\\[0.4cm]
&&\qquad +e^{-px}\Big(\begin{array}{l}-g_u(\overline{u}_h)e^{px}\phi(x)+e^{px}\psi(x)-\lambda e^{px}\phi(x)\\ +\gamma\rho e^{px}\psi(x)-\lambda e^{px}\psi(x)\end{array}\Big)\\[0.4cm]
&=&\left(\begin{array}{l}pc_h\phi(x)+c_h\phi'(x)-g_u(\overline{u}_h)\phi(x)+\psi(x)\\ -\rho \phi(x)+pc_h\psi(x)+c_h\psi'(x)+\gamma\rho \psi(x)-\lambda \psi(x)\end{array}\right)\\[0.4cm]
&&\qquad +\left(\begin{array}{l}-\Delta_h\phi(x)-\lambda \phi(x)\\ 0\end{array}\right)\\[0.4cm]
&=&(L_h-\lambda+pc_h)\Psi(x).\end{array}\end{equation}
Since $e_{p}$ and $e_{-p}$ are invertible operators on $\textbf{H}^1$ and $\textbf{L}^2$ respectively, we know that the spectrum of $L_h$ equals that of $e_{-p}L_he_{p}$ and thus that of $L_h+pc_h$.\qed

\vspace*{4pt}\Large{\textbf{Region $R_1$.}}\vspace*{4pt}
\normalsize

Since $L_h$ has a simple eigenvalue at zero, it is relatively straightforward to construct a small neighbourhood around the origin that contains no other part of the spectrum. Exploiting the results from \S \ref{waveproperties}, it is possible to control the size of this neighbourhood as $h\downarrow 0$.
\begin{proposition}\label{lemmaspectrumklein} Assume that (\asref{aannamespuls}{\text{H}}),(\asref{extraaannamespuls}{\text{H}}), (\asref{aannamesconstanten}{\text{H}}), (\asref{aannames}{\text{H}\alpha1}) and (\asref{extraaannames}{\text{H}\alpha}) are satisfied. There exists a constant $\lambda_0>0$ such that for all $0<h< h_{**}$ the operator $L_h+\lambda:\textbf{H}^1\rightarrow \textbf{L}^2$ is invertible for all $\lambda\in\C$ with $0<|\lambda|<\lambda_0$.\end{proposition}
\vspace*{4pt}\noindent\textit{Proof.} Fix $0<h< h_{**}$ and $\Theta\in\textbf{L}^2$. We recall the notation $(\gamma_h[\Theta],L_h^{\mathrm{qinv}}\Theta)$ from Corollary \ref{vhwh} for the unique solution $(\gamma,\Psi)$ of the equation
\begin{equation}\begin{array}{lcl}L_h\Psi&=&\Theta+\gamma\Phi_h^+\end{array}\end{equation}
in the space
\begin{equation}\begin{array}{lcl}X_h&=&\{\Theta\in \textbf{H}^1: \ip{\Phi_h^-,\Theta}=0\}.\end{array}\end{equation}
Also recall the space
\begin{equation}\begin{array}{lcl}Y_h&=&\{\Theta\in \textbf{L}^2: \ip{\Phi_h^-,\Theta}=0\}.\end{array}\end{equation}
Now for $\lambda\in\C$ with $|\lambda|$ small enough, but $\lambda\neq 0$, we want to solve the equation $L_h\Psi=\lambda\Psi+\Theta$. Upon writing
\begin{equation}\begin{array}{lcl}\Psi&=&L_h^{\mathrm{qinv}}\Theta+\lambda^{-1}\gamma_h[\Theta]\Phi_h^++\tilde{\Psi},\end{array}\end{equation}
with $\tilde{\Psi}\in X_h$, we see that
\begin{equation}\begin{array}{lcl}(L_h-\lambda)\Psi&=&(L_h-\lambda)L_h^{\mathrm{qinv}}\Theta+\lambda^{-1}(L_h-\lambda)\gamma_h[\Theta]\Phi_h^++(L_h-\lambda)\tilde{\Psi}\\[0.2cm]
&=&\Theta+\gamma_h[\Theta]\Phi_h^+-\lambda L_h^{\mathrm{qinv}}\Theta-\gamma_h[\Theta]\Phi_h^++(L_h-\lambda)\tilde{\Psi}.\end{array}\end{equation}
In particular, we must find a solution $\tilde{\Psi}\in X_h$ for the equation
\begin{equation}\begin{array}{lcl}L_h\tilde{\Psi}&=&\lambda\tilde{\Psi}+\lambda L_h^{\mathrm{qinv}}\Theta,\end{array}\end{equation}
which we can rewrite as
\begin{equation}\label{equivalentevergelijking}\begin{array}{lcl}[I-\lambda L_h^{-1}]\tilde{\Psi}&=&\lambda L_h^{-1}L_h^{\mathrm{qinv}}\Theta.\end{array}\end{equation}
Note that $L_h^{-1}:X_h\rightarrow X_h$ is also a bounded operator since $X_h\subset Y_h$. Since
\begin{equation}\begin{array}{lcl}\nrm{L_h^{-1}\Psi}_{\textbf{H}^1}&\leq & C_{\mathrm{unif}}\nrm{\Psi}_{\textbf{L}^2}\\[0.2cm]
&\leq &C_{\mathrm{unif}}\nrm{\Psi}_{\textbf{H}^1},\end{array}\end{equation}
we obtain
\begin{equation}
\begin{array}{lcl}
\nrm{L_h^{-1}}_{\mathcal{B}(X_h,X_h)}&\leq &C_{\mathrm{unif}}.
\end{array}
\end{equation}
We choose $\lambda_0$ in such a way that $0<\lambda_0C_{\mathrm{unif}}<1$.
Then it is well-known that $I-\lambda L_h^{-1}$ is invertible as an operator on $X_h$ for $0<|\lambda|<\lambda_0$. Since $\lambda L_h^{-1}L_h^{\mathrm{qinv}}\Theta\in X_h$, we see that (\ref{equivalentevergelijking}) indeed has a unique solution $\tilde{\Psi}\in X_h$. Hence, the equation $(L_h-\lambda)\Psi=\Theta$ always has a unique solution. Proposition \ref{LhlambdaFredholm} states that $L_h-\lambda$ is Fredholm with index $0$, which now implies that $L_h-\lambda$ is indeed invertible.\qed

\vspace*{4pt}\Large{\textbf{Region $R_2$}.}\vspace*{4pt}
\normalsize

We now show that in an appropriate right half-plane, which can be chosen independently of $h$, the spectrum of $-L_h$ is empty. The proof proceeds via a relatively direct estimate that is strongly inspired by \cite[Lemma 5]{BatesInfRange}.
\begin{proposition}\label{spectrumgroot} Assume that (\asref{aannamespuls}{\text{H}}),(\asref{extraaannamespuls}{\text{H}}), (\asref{aannamesconstanten}{\text{H}}), (\asref{aannames}{\text{H}\alpha1}) and (\asref{extraaannames}{\text{H}\alpha}) are satisfied. There exists a constant $\lambda_1>0$ such that for all $\lambda\in\C$ with $\Re \,\lambda\geq\lambda_1$ and all $0<h< h_{**}$ the operator $L_h+\lambda$ is invertible.\end{proposition}
\vspace*{4pt}\noindent\textit{Proof.} Write
\begin{equation}\begin{array}{lcl}\lambda_1&=&1+g_*+\frac{1}{2}(1-\rho),\end{array}\end{equation}
where $g_*$ is defined in Lemma \ref{lemma6bewijs2}. Pick any $\lambda\in\C$ with $\Re \,\lambda\geq \lambda_1$ and any $0<h< h_{**}$. Let $\Psi=(\phi,\psi)\in\textbf{H}^1$ be arbitrary and set $\Theta=L_h\Psi+\lambda\Psi$. Then we see that
\begin{equation}\begin{array}{lcl}\nrm{\Psi}_{\textbf{L}^2}\nrm{\Theta}_{\textbf{L}^2}&\geq &\Re \,\ip{L_h\Psi+\lambda\Psi,\Psi}\\[0.2cm]
&\geq &\Re \,\ip{-\Delta_h\phi,\phi}-\nrm{g_u(\overline{u}_h)}_{L^\infty}\nrm{\phi}_{L^2}^2\\[0.2cm]
&&\qquad -(1-\rho)|\Re \,\ip{\phi,\psi}|+\gamma\rho\nrm{\psi}_{L^2}^2+\Re \,\lambda\nrm{\Psi}_{\textbf{L}^2}^2\\[0.2cm]
&\geq &-g_*\nrm{\phi}_{L^2}^2
-(1-\rho)|\Re \,\ip{\phi,\psi}|
 +\gamma\rho\nrm{\psi}_{L^2}^2+\Re \,\lambda\nrm{\Psi}_{\textbf{L}^2}^2\\[0.2cm]
&\geq &-g_*\nrm{\phi}_{L^2}^2
-(1-\rho)\nrm{\phi}_{L^2}\nrm{\psi}_{L^2}
+\gamma\rho\nrm{\psi}_{L^2}^2+\Re \,\lambda\nrm{\Psi}_{\textbf{L}^2}^2\\[0.2cm]
&\geq &-(g_*+\frac{1}{2}(1-\rho))\nrm{\Psi}_{\textbf{L}^2}^2
+\Re \,\lambda\nrm{\Psi}_{\textbf{L}^2}^2.\end{array}\end{equation}
Hence, we obtain
\begin{equation}\begin{array}{lcl}\Big(\Re \,\lambda-(g_*+\frac{1}{2}(1-\rho))\Big)\nrm{\Psi}_{\textbf{L}^2}&\leq &\nrm{\Theta}_{\textbf{L}^2}.\end{array}\end{equation}
Since $\Re \,\lambda\geq 1+g_*+\frac{1}{2}(1-\rho)$, we obtain the bound $\nrm{\Psi}_{\textbf{L}^2}\leq \nrm{\Theta}_{\textbf{L}^2}$.\par

In particular, if $\Theta=0$ then we necessarily have $\Psi=0$, which implies that $L_h+\lambda$ is injective. Since also $\mathrm{ind}(L_h+\lambda)=0$ by Proposition \ref{LhlambdaFredholm}, this means that $L_h+\lambda$ is invertible.
\qed

\vspace*{4pt}\Large{\textbf{Region $R_3$}.}\vspace*{4pt}
\normalsize

This region is the most delicate to handle on account of the periodicity of the spectrum.
Indeed, one cannot simply take $\Im \lambda \to \pm \infty$ in a fashion that is uniform in $h$.
We pursue a direct approach here, using the Fourier transform to isolate the
problematic part of $L_h + \lambda$, which has constant coefficients.
The corresponding portion of the resolvent can be estimated in a
controlled way by rescaling the imaginary part of $\lambda$.
We remark that an alternative approach could be to factor out the periodicity in
a more operator-theoretic setting, but we do not pursue such an argument here.\par


Pick $\lambda\in\C$ with $\lambda_0<|\Im \,\lambda|\leq\frac{|c_h|}{h}\pi$ and write
\begin{equation}\begin{array}{lcl}\lambda&=&\lambda_\mathrm{r}+i\lambda_{\mathrm{im}}\end{array}\end{equation}
with $\lambda_\mathrm{r},\lambda_{\mathrm{im}}\in\R$. Introducing the new variable $\tau=\Im \,\lambda\xi$, we can write the eigenvalue problem $(L_h+\lambda)(v,w)=0$ in the form
\begin{equation}\label{systeminothervariable}\begin{array}{lcl}c_hv_\tau(\tau)&=&  \frac{1}{\lambda_\mathrm{im} h^2}\sum\limits_{k>0}\alpha_k\Big[v(\tau+kh\lambda_\mathrm{im})+v(\tau-kh\lambda_\mathrm{im})-2v(\tau)\Big]\\[0.2cm]
&&\qquad+\frac{1}{\lambda_\mathrm{im}}g_u\Big(\overline{u}_h(\tau)\Big)v(\tau)
-iv(\tau)-\frac{1}{\lambda_\mathrm{im}}\lambda_\mathrm{r} v(\tau)-\frac{1}{\lambda_\mathrm{im}}w(\tau),\\[0.2cm]
c_hw_\tau(\tau)&=&\frac{1}{\lambda_\mathrm{im}}\Big(\rho v(\tau)-\rho\gamma w(\tau)+\lambda w(\tau)\Big).\end{array}\end{equation}
Our computations below show that the leading order terms in the appropriate $|\lambda_{\mathrm{im}}|\rightarrow \infty$ limit are encoded by the  `homogeneous operator' $H_{h,\lambda}$ that acts as
\begin{equation}\begin{array}{lcl}H_{h,\lambda}v(\tau)&=&c_hv_\tau(\tau)+iv(\tau)-\frac{1}{\lambda_\mathrm{im} h^2}\sum\limits_{k>0}\alpha_k\Big[v(\tau+kh\lambda)+v(\tau-kh\lambda)-2v(\tau)\Big].\end{array}\end{equation}
Writing $\mathcal{H}_{h,\lambda}$ for the Fourier symbol associated to $H_{h,\lambda}$, we see that
\begin{equation}\begin{array}{lcl}\mathcal{H}_{h,\lambda}(i\omega)&=&c_hi\omega+i-\frac{1}{\lambda_\mathrm{im} h^2}\sum\limits_{k>0}\alpha_k\Big[\exp(ihk\lambda_\mathrm{im} \omega)+\exp(-ihk\lambda_\mathrm{im}\omega)-2\Big]\\[0.2cm]
&=&c_hi\omega+i-\frac{2}{\lambda_\mathrm{im} h^2}\sum\limits_{k>0}\alpha_k\Big[\cos(hk\lambda_\mathrm{im} \omega)-1\Big].\end{array}\end{equation}

\begin{lemma}\label{spectrumhoogbewijs0} Assume that (\asref{aannamespuls}{\text{H}}),(\asref{extraaannamespuls}{\text{H}}), (\asref{aannamesconstanten}{\text{H}}), (\asref{aannames}{\text{H}\alpha1}) and (\asref{extraaannames}{\text{H}\alpha}) are satisfied. There exist small constants $\epsilon>0$, $h_*>0$ and $\omega_0>0$ so that for all $\lambda\in\C\setminus\R$, all $0<h<h_*$ and all $\omega\in\R$, the inequality
\begin{equation}
\begin{array}{lcl}
|\Im \,\mathcal{H}_{h,\lambda}(i\omega)|&<&\epsilon
\end{array}
\end{equation}
can only be satisfied if the inequalities
\begin{equation}\label{imHhongelijkheiden}\begin{array}{lcl}|c_{h}\omega|&\leq & \frac{3}{2}\\[0.2cm]
|\omega|&\geq &\omega_0\end{array}\end{equation}
both hold.\end{lemma}
\noindent\textit{Proof.} Note that
\begin{equation}
\begin{array}{lcl}
|\Im \,\mathcal{H}_{h,\lambda}(i\omega)|&=&|c_{h}\omega+1|.
\end{array}
\end{equation}
In particular, upon choosing $\epsilon=\frac{1}{4}$, we see that
\begin{equation}
\begin{array}{lcl}
|\Im \,\mathcal{H}_{h,\lambda}(i\omega)|&<&\epsilon
\end{array}
\end{equation}
implies
\begin{equation}
\begin{array}{lclcl}
\big||c_h\omega|-1\big|&\leq &|c_h\omega+1|
&<&\epsilon
\end{array}
\end{equation}
and hence
\begin{equation}
\begin{array}{lclclclcl}
\frac{1}{2}&<& 1-\epsilon
&\leq & |c_h\omega|
&\leq & 1+\epsilon
&<&\frac{3}{2}.
\end{array}
\end{equation}
Since $c_h\rightarrow c_0\neq 0$ as $h\downarrow 0$, the desired inequalities (\ref{imHhongelijkheiden}) follow. \qed

\begin{lemma}\label{spectrumhoogbewijs} Assume that (\asref{aannamespuls}{\text{H}}),(\asref{extraaannamespuls}{\text{H}}), (\asref{aannamesconstanten}{\text{H}}), (\asref{aannames}{\text{H}\alpha1}) and (\asref{extraaannames}{\text{H}\alpha}) are satisfied. Then there exists a constant $C>0$ such that for all $\omega\in\R$ and $0<h< h_{**}$ and all $\lambda\in\C$ with $|\lambda|>\lambda_0$ and $|\Im\lambda|\leq\frac{|c_h|}{h}\pi$, we have the inequality
\begin{equation}\begin{array}{lcl}|\mathcal{H}_{h,\lambda}(i\omega)|\geq\frac{1}{C}.\end{array}\end{equation}\end{lemma}
\vspace*{4pt}\noindent\textit{Proof.} We show that $\mathcal{H}_{h,\lambda}(i\omega)$ is bounded away from $0$, uniformly in $h,\lambda$ and $\omega$. To do so, we show that the real part of $\mathcal{H}_{h,\lambda}(i\omega)$ can be bounded away from zero, whenever the imaginary part is small, i.e. when (\ref{imHhongelijkheiden}) holds.\par

Recall the function $A(y)=\sum\limits_{k>0}\alpha_k[1-\cos(ky)]$ defined in Assumption (\asref{aannames}{\text{H}\alpha1}), which satisfies $A(y)>0$ for $y\in(0,2\pi)$. A direct calculation shows that $A'(0)=0$ and
\begin{equation}\label{d0}\begin{array}{lcl}A''(0)&=&\sum\limits_{k>0}\alpha_kk^2\\[0.2cm]
&=&1.\end{array}\end{equation}
Hence, we can pick $d_0>0$ in such a way that
\begin{equation}
\begin{array}{lcl}
\frac{1}{y^2}A(y)&>&d_0
\end{array}
\end{equation}
holds for all $0<|y|\leq \frac{3}{2}\pi$. \par

Writing $\mu=h\lambda_\mathrm{im}\omega$, we see
\begin{equation}
\begin{array}{lcl}
\Re \,\mathcal{H}_{h,\lambda}(i\omega)&=&\frac{2\omega^2\lambda_{\mathrm{im}}}{\mu^2}\sum\limits_{k>0}\alpha_k\Big[1-\cos(k\mu)\Big]\\[0.2cm]
&=&\frac{2\omega^2\lambda_{\mathrm{im}}}{\mu^2}A(\mu).
\end{array}
\end{equation}
Now fix $\omega,h,\lambda$ for which $|\Im \,\mathcal{H}_{h,\lambda}(i\omega)|<\epsilon$. The conditions (\ref{imHhongelijkheiden}) now imply that $|\omega|\geq \omega_0$ and $|\mu|\leq h\frac{|c_h|}{h}\pi|\omega|\leq \frac{3}{2}\pi$. Using (\ref{d0}), we hence see that
\begin{equation}\begin{array}{lcl}|\Re \,\mathcal{H}_{h,\lambda}(i\omega)|&=&|\frac{2\omega^2\lambda_{\mathrm{im}}}{\mu^2}A(\mu)|\\[0.2cm]
&\geq &2|\lambda_\mathrm{im}|\omega^2 d_0\\[0.2cm]
&\geq &2\lambda_0\omega_0^2d_0,\end{array}\end{equation}
which shows that $\mathcal{H}_{h,\lambda}(i\omega)$ can indeed be uniformly bounded away from zero.\qed

\begin{proposition}\label{spectrumhoog} Assume that (\asref{aannamespuls}{\text{H}}),(\asref{extraaannamespuls}{\text{H}}), (\asref{aannamesconstanten}{\text{H}}), (\asref{aannames}{\text{H}\alpha1}) and (\asref{extraaannames}{\text{H}\alpha}) are satisfied. There exist constants $\lambda_2>0$ and $\lambda_3>0$ such that for all $\lambda\in\C$ with $\lambda_2\leq|\Im \, \lambda|\leq \frac{|c_h|}{2h}2\pi$ and $-\lambda_3\leq|\Re \,\lambda|\leq \lambda_1$ and all $0<h< h_{**}$ the operator $L_h+\lambda$ is invertible.\end{proposition}
\noindent\textit{Proof.}
Since Proposition \ref{LhlambdaFredholm} implies that $L_h+\lambda$ is Fredholm with index zero, it suffices to prove that $L_h+\lambda$ is injective.\par

Let $\lambda_3=\min\{\frac{1}{2}\rho\gamma,\lambda_*,\tilde{\lambda}\}$, where $\lambda_*$ is defined in (\asref{extraaannamespuls}{\text{H}}) and $\tilde{\lambda}$ is defined in Proposition \ref{LhlambdaFredholm}. Pick $\lambda\in\C$ with $\lambda_0\leq|\Im \, \lambda|\leq \frac{|c_h|}{2h}2\pi$ and $-\lambda_3\leq|\Re \,\lambda|\leq \lambda_1$. Write $\lambda=\lambda_\mathrm{r}+i\lambda_{\mathrm{im}}$ as before. Suppose $\Psi=(v,w)$ satisfies $(L_h+\lambda)\Psi=0$. \par

Write $\hat{v}$ and $\hat{w}$ for the Fourier transforms of $v$ and $w$ respectively. For $f\in L^2$ with Fourier transform $\hat{f}$, the identity
\begin{equation}
\begin{array}{lcl}
H_{h,\lambda}v&=&f
\end{array}
\end{equation}
implies that
\begin{equation}\begin{array}{lcl}\mathcal{H}_{h,\lambda}(i\omega)\hat{v}(i\omega)&=&\hat{f}(i\omega).\end{array}\end{equation}
In particular, we obtain
\begin{equation}\begin{array}{lcl}\hat{v}(i\omega)&=&\frac{1}{\mathcal{H}_{h,\lambda}(i\omega)}\hat{f}(i\omega),\end{array}\end{equation}
which using Lemma \ref{spectrumhoogbewijs} implies that
\begin{equation}\begin{array}{lcl}\nrm{v}_{L^2}&\leq &C\nrm{f}_{L^2}\end{array}\end{equation}
for some constant $C>0$ that is independent of $h,\lambda$ and $\omega$.\par

Since $\Psi$ is an eigenfunction, (\ref{systeminothervariable}) hence yields
\begin{equation}\begin{array}{lcl}\nrm{v}_{L^2}&\leq &C\frac{1}{|\lambda_\mathrm{im}|}(g_*+|\lambda_\mathrm{r}|)\nrm{v}_{L^2}+C\frac{1}{|\lambda_\mathrm{im}|}\nrm{w}_{L^2}.\end{array}\end{equation}
Furthermore, applying a Fourier Transform to the second line of (\ref{systeminothervariable}), we find
\begin{equation}\begin{array}{lcl}\lambda_\mathrm{im} c_hi\omega\hat{w}(i\omega)&=&\rho \hat{v}(i\omega)-\rho\gamma \hat{w}(i\omega)+\lambda \hat{w}(i\omega).\end{array}\end{equation}
Our choice $\lambda_3\leq \frac{1}{2}\rho\gamma$ implies that $-\rho\gamma+\lambda_\mathrm{r}$ is bounded away from $0$. We may hence write
\begin{equation}\begin{array}{lcl}\hat{w}(i\omega)&=&\frac{1}{\rho\gamma-\lambda_\mathrm{r}+i(\omega\lambda_\mathrm{im} c_h-\lambda_\mathrm{im})}\rho\hat{v}(i\omega),\end{array}\end{equation}
which yields the bound
\begin{equation}\begin{array}{lcl}\nrm{w}_{L^2}&\leq &C'\nrm{v}_{L^2}\end{array}\end{equation}
for some constant $C'>0$. Therefore, we obtain that
\begin{equation}\begin{array}{lcl}\nrm{v}_{L^2}&\leq &C''\frac{1}{|\lambda_\mathrm{im}|}\nrm{v}_{L^2}\end{array}\end{equation}
for some constant $C''$, which is independent of $\lambda,h$ and $v$. Clearly this is impossible for $v\neq 0$ if
\begin{equation}\begin{array}{lclcl}|\lambda_\mathrm{im}|&\geq&\lambda_2&:=&2C''.\end{array}\end{equation}
Furthermore, if $v=0$, then clearly also $w=0$. Therefore, we have $\Psi=0$, allowing us to conclude that $L_h+\lambda$ is invertible. \qed

\vspace*{4pt}\Large{\textbf{Region $R_4$.}}\vspace*{4pt}
\normalsize

We conclude our spectral analysis by considering the remaining region $R_4$. This region is compact and bounded away from the origin, allowing us to directly apply the theory developed in \S \ref{singularoperator}.
\begin{corollary}\label{spectrumintermediate} Assume that (\asref{aannamespuls}{\text{H}}),(\asref{extraaannamespuls}{\text{H}}), (\asref{aannamesconstanten}{\text{H}}), (\asref{aannames}{\text{H}\alpha1}) and (\asref{extraaannames}{\text{H}\alpha}) are satisfied. For all $\lambda\in\C$ with $|\lambda|\geq\lambda_0$, $-\lambda_3\leq|\Re \,\lambda|\leq\lambda_1$ and $|\Im \,\lambda|\leq\lambda_2$ and all $0<h< h_{**}$ the operator $L_h+\lambda$ is invertible.\end{corollary}
\noindent\textit{Proof.} The statement follows by applying Proposition \ref{equivalenttheorem4version2} with the choices $(\tilde{u}_h,\tilde{w}_h)\break=(\overline{u}_h,\overline{w}_h)$, $\tilde{c}_h=c_h$ and $M=R_4$.\qed

\vspace*{4pt}\noindent\textit{Proof of Theorem \ref{totalespectrum}.} The result follows directly from Lemma \ref{lemmaperiodiek}, Proposition \ref{lemmaspectrumklein}, Proposition \ref{spectrumgroot}, Proposition \ref{spectrumhoog} and Corollary \ref{spectrumintermediate}.\qed


\section{Green's functions}\label{sectiongreen}

In order to establish the nonlinear stability of the pulse solution $(\overline{u}_h,\overline{w}_h)$, we need to consider two types of Green's functions. In particular, we first study $G_{\lambda}(\xi,\xi_0)$, which can roughly be seen as a solution of the equation
\begin{equation}\begin{array}{lcl}\Big[(L_h+\lambda)G_{\lambda}(\cdot,\xi_0)\Big](\xi)=\delta(\xi-\xi_0),\end{array}\end{equation}
where $\delta$ is the Dirac delta-distribution. We then use these functions to build a
Green's function $\mathcal{G}$ for the linearisation of the LDE (\ref{ditishetalgemeneprobleem}) around the travelling pulse solution.\par

An important difficulty in comparison to the PDE setting is caused by the discreteness of the
spatial variable $j$. In particular, we cannot use a frame of reference in which the solution
$(\overline{u}_h,\overline{w}_h)$ is constant without changing the structure of the equation
(\ref{ditishetalgemeneprobleem}). The Green's function $\mathcal{G}$
will hence be the solution to a non-autonomous problem that satisfies a shift-periodicity condition.
Nevertheless, one can follow the technique
in \cite{BGV2003} to express $\mathcal{G}$ in terms of a contour integral involving
the functions $G_{\lambda}$.\par

A significant part of our effort here is concerned with the construction of these latter functions.
Indeed, previous approaches in \cite{HJHNLS, HJHSTBFHN} all used exponential dichotomies or
variation-of-constants formula's for MFDEs with finite-range interactions. These tools are no longer
available for use in the present infinite-range setting. In particular, we construct
the functions $G_{\lambda}$ in a direct fashion using only Fredholm properties of the operators
$L_h + \lambda$. This makes it somewhat involved to recover the desired exponential decay rates
and to properly isolate the meromorphic terms of order $O(\lambda^{-1})$.\par


From now on, we will no longer explicitly use the $h$-dependence of our system. To simplify our notation, we fix $0<h<h_{**}$ and write
\begin{equation}
\begin{array}{lclcl}
L&:=&L_h,\\
L_\infty&:=&L_{h;\infty},\\
\overline{U}&=&(\overline{u},\overline{w})
&:=&(\overline{u}_h,\overline{w}_h),\\

\Phi^{\pm}&=&(\phi^{\pm},\psi^{\pm})
&:=&(\phi_h^{\pm},\psi_h^{\pm}),\\
c&:=&c_h.
\end{array}
\end{equation}
We emphasize that from now on all our constants may (and will) depend on $h$.\par

We will loosely follow \S 2 of \cite{HJHSTBFHN},
borrowing a number of results from \cite{BGV2003,HJHPOLLUTION} at appropriate times.
%
In particular, we start by
considering the linearisation of the original LDE (\ref{ditishetalgemeneprobleem}) around the travelling
pulse solution $\overline{U}(t)$ given by (\ref{definitieUh}).
To this end,
we introduce the Hilbert space
\begin{equation}\begin{array}{lcl}\mathbb{L}^2&:=&\{V\in(\text{Mat}_2(\R))^\Z:\sum\limits_{j\in\Z}|V(j)|^2<\infty\},\end{array}\end{equation}
in which $\text{Mat}_2(\R)$ is the space of $2\times 2$-matrices with real coefficients which we equip with the maximum-norm $|\cdot|$.
For any $\mathcal{V}\in\mathbb{L}^2$, we often write $\mathcal{V}=\left(\begin{array}{ll}\mathcal{V}^{(1,1)}&\mathcal{V}^{(1,2)}\\ \mathcal{V}^{(2,1)}&\mathcal{V}^{(2,2)}\end{array}\right)$, when we need to refer to the component sequences $\mathcal{V}^{(i,j)}\in\ell^2(\Z;\R)$. For any $t\in\R$ we now introduce the linear operator $\mathcal{A}(t):\mathbb{L}^2\rightarrow\mathbb{L}^2$  that acts as
\begin{equation}\label{defAt}\begin{array}{lcl}\mathcal{A}(t)\cdot \mathcal{V}&=&\frac{1}{c}\left(\begin{array}{ll}A^{(1,1)}(t)&A^{(1,2)}(t)\\ A^{(2,1)}(t)&A^{(2,2)}(t)\end{array}\right)\left(\begin{array}{ll}\mathcal{V}^{(1,1)}&\mathcal{V}^{(1,2)}\\ \mathcal{V}^{(2,1)}&\mathcal{V}^{(2,2)}\end{array}\right),\end{array}\end{equation}
where
\begin{equation}\begin{array}{lcl}(A^{(1,1)}(t)v)_j&=&\frac{1}{h^2}\sum\limits_{k>0} \alpha_k[v_{j+k}+v_{j-k}-2v_j]+g_u\Big(\overline{u}(hj+ct)\Big)v_{j}\\[0.2cm]
(A^{(1,2)}(t)w)_j&=&-w_j\\[0.2cm]
(A^{(2,1)}(t)v)_j&=& \rho v_j\\[0.2cm]
(A^{(2,2)}(t)v)_j&=&-\rho\gamma w_{j}\end{array}\end{equation}
for $v\in \ell^2(\Z;\R)$ and $w\in\ell^2(\Z;\R)$. With all this notation in hand, we can write the desired linearisation as the ODE
\begin{equation}\label{definitieGhuh}\begin{array}{lcl}\frac{d}{dt}\mathcal{V}(t)&=&\mathcal{A}(t)\cdot \mathcal{V}(t)\end{array}\end{equation}
posed on $\mathbb{L}^2$.\par

Fix $t_0\in\R$ and $j_0\in\Z$. Consider the function
\begin{equation}
\begin{array}{lcl}
\R\ni t\mapsto \mathcal{G}^{j_0}(t,t_0)&=&\{\mathcal{G}_j^{j_0}(t,t_0)\}_{j\in\Z}\in\mathbb{L}^2
\end{array}
\end{equation}
that is uniquely determined by the initial value problem
\begin{equation}\label{3.2equivalent}\begin{cases}\frac{d}{dt}\mathcal{G}^{j_0}(t,t_0)&=\mathcal{A}(t)\cdot \mathcal{G}^{j_0}(t,t_0) \\ \mathcal{G}_j^{j_0}(t_0,t_0)&=\delta_j^{j_0}I.\end{cases}\end{equation}
Here we have introduced
\begin{equation}\begin{array}{lcl}\delta_j^{j_0}&=&\begin{cases} 1\text{ if }j=j_0\\ 0\text{ else,}\end{cases}\end{array}\end{equation}
where $I\in\text{Mat}_2(\R)$ is the identity matrix. We remark that $\mathcal{G}_j^{j_0}(t,t_0)$ is an element of $\text{Mat}_2(\R)$ for each $j\in\Z$.\par

This function $\mathcal{G}$ is called the Green's function for the linearisation around our travelling pulse. Indeed, the general solution of the inhomogeneous equation
\begin{equation}\begin{cases}\frac{dV}{dt}&=\enskip \mathcal{A}(t)\cdot V(t)+F(t)\\
V(0)&=\enskip V^0,\end{cases}\end{equation}
where now $V(t)\in \ell^2(\Z;\R^2)\cong \ell^2(\Z;\R^{2\times 1})$ and $F(t)\in \ell^2(\Z;\R^2)\cong \ell^2(\Z;\R^{2\times 1})$, is given by
\begin{equation}\begin{array}{lcl}V_j(t)&=&\sum\limits_{j_0\in\Z}\mathcal{G}_j^{j_0}(t,0)V_{j_0}^0+\int_0^t \sum\limits_{j_0\in\Z}\mathcal{G}_j^{j_0}(t,t_0)F_{j_0}(t_0)\ dt_0.\end{array}\end{equation}
Introducing the standard convolution operator $\ast$, this can be written in the abbreviated form
\begin{equation}\label{ast}
\begin{array}{lcl}
V&=&\mathcal{G}(t,0)\ast V^0+\int_0^t \mathcal{G}(t,t_0)\ast F(t_0)\ dt_0.
\end{array}
\end{equation}

The main result of this section is the following proposition, which shows that we can decompose the Green's function $\mathcal{G}$ into a part that decays exponentially and a neutral part associated with translation along the family of travelling pulses.
\begin{proposition}\label{cor2.8equivalentsterker} Assume that (\asref{aannamespuls}{\text{H}}),(\asref{extraaannamespuls}{\text{H}}), (\asref{aannamesconstanten}{\text{H}}), (\asref{aannames}{\text{H}\alpha1}) and (\asref{extraaannames}{\text{H}\alpha}) are satisfied. For any pair $t\geq t_0$ and any $j,j_0\in\Z$, we have the representation
\begin{equation}\label{zogaanwehemopdelen}\begin{array}{lcl}\mathcal{G}_j^{j_0}(t,t_0)&=&\mathcal{E}_j^{j_0}(t,t_0)+\tilde{\mathcal{G}}_j^{j_0}(t,t_0),\end{array}\end{equation}
in which
\begin{equation}\begin{array}{lcl}\mathcal{E}_j^{j_0}(t,t_0)&=&\frac{h}{\Omega}\left(\begin{array}{ll}\phi^-(hj_0+ct_0)\phi^+(hj+ct)& \psi^-(hj_0+ct_0)\phi^+(hj+ct)\\ \phi^-(hj_0+ct_0)\psi^+(hj+ct)& \psi^-(hj_0+ct_0)\psi^+(hj+ct)\end{array}\right),\end{array}\end{equation}
while $\tilde{\mathcal{G}}$ satisfies the found
\begin{equation}\label{eq:expdecaymathcalG}\begin{array}{lcl}|\tilde{\mathcal{G}}_j^{j_0}(t,t_0)| &\leq & K e^{-\tilde{\delta}(t-t_0)}e^{-\tilde{\delta}|hj+ct-hj_0-ct_0|}\end{array}\end{equation}
for some $K>0$ and $\tilde{\delta}>0$. The constant $\Omega>0$ is given by
\begin{equation}\begin{array}{lcl}\Omega&=&\ip{\Phi^-,\Phi^+}.\end{array}\end{equation}
Furthermore, for any $t\geq t_0$ we have the representation
\begin{equation}\begin{array}{lcl}\mathcal{G}_j^{j_0}(t,t_0)&=&\sum\limits_{i\in\Z}\Big[\mathcal{E}_{j}^i(t,t_0)\mathcal{E}_{i}^{j_0}(t_0,t_0)+\tilde{\mathcal{G}}_j^i(t,t_0)(\delta_{i}^{j_0}I-\mathcal{E}_{i}^{j_0}(t_0,t_0))\Big],\end{array}\end{equation}
which can be abbreviated as
\begin{equation}
\begin{array}{lcl}
\mathcal{G}(t,t_0)&=& \mathcal{E}(t,t_0)\ast\mathcal{E}(t_0,t_0)+\tilde{\mathcal{G}}(t,t_0)\ast\Big(I-\mathcal{E}(t_0,t_0)\Big).
\end{array}
\end{equation}
\end{proposition}
\noindent\textit{Proof of Theorem \ref{nonlinearstability}.}
The Green's function bounds from Proposition \ref{cor2.8equivalentsterker} are strong
enough to follow the program described in the proof of \cite[Prop 2.1]{HJHSTBFHN}.
In particular, using standard fixed point arguments one can foliate the
$\ell^p$-neighbourhood of the wave $\overline{U}$ with the one-parameter family
formed by the stable manifolds of the translates $\overline{U}(\cdot + \vartheta)$.
\qed

\subsection{Construction of the Green's function}
In this subsection we set out to define the functions $G_{\lambda}$ in a more rigorous fashion.
In addition, we use these Green's functions
to formulate
a powerful representation formula for $\mathcal{G}$, see Proposition \ref{thrGreenfunction} below,
following the approach developed in \cite{BGV2003}.

A key role in our analysis is reserved for the operator $L_{\infty;\lambda}$ and the function $\Delta_{L_{\infty;\lambda}}$ from Lemma \ref{hyperbolic}. We will show that $L_{\infty;\lambda}$ has a Green's function $G_{\infty;\lambda}$ which takes values in the space $\text{Mat}_2(\R)$ and
has some useful properties. To this end, we recall the constant $\tilde{\lambda}$ from Lemma \ref{hyperbolic}. For each $\lambda\in\C$ with $\Re \,\lambda\geq-\frac{\tilde{\lambda}}{2}$,
we may now define $G_{\infty;\lambda}:\R\rightarrow \text{Mat}_2(\R)$ by writing
\begin{equation}\label{defvanghinfty}\begin{array}{lcl}G_{\infty;\lambda}(\xi)&=&\frac{1}{2\pi }\int_{-\infty}^\infty e^{i\eta\xi}(\Delta_{L_{\infty;\lambda}}(i\eta))^{-1}\ d\eta.\end{array}\end{equation}
We also introduce the notation
\begin{equation}\begin{array}{lcl}G_{\infty}&=&G_{\infty;0}.\end{array}\end{equation}
Here (\asref{extraaannames}{\text{H}\alpha}) is essential to ensure that these Green's functions decay exponentially.
\begin{lemma}\label{eigenschappenghinfty} Assume that (\asref{aannamespuls}{\text{H}}), (\asref{extraaannamespuls}{\text{H}}), (\asref{aannamesconstanten}{\text{H}}), (\asref{aannames}{\text{H}\alpha1}) and (\asref{extraaannames}{\text{H}\alpha}) are satisfied. Fix $\lambda\in\C$ with $\Re \,\lambda\geq-\frac{\tilde{\lambda}}{2}$. The function $G_{\infty;\lambda}$ is bounded and continuous on $\R\setminus\{0\}$ and $C^1$-smooth on $\R\setminus h\Z$. Furthermore, $(L_{\infty}+\lambda)G_{\infty;\lambda}(\cdot-\xi_0)$ is constantly zero except at $\xi=\xi_0+h\Z$ and satisfies the identity
\begin{equation}\begin{array}{lcl}\int_{-\infty}^\infty \Big[(L_{\infty}+\lambda) G_{\infty;\lambda} (\cdot-\xi_0)\Big](\xi) f(\xi)\ d\xi
&=&f(\xi_0)\end{array}\end{equation}
for all $\xi\in\R$ and all $f\in \textbf{H}^1$.\par

Finally for each $\chi>0$ there exist constants $K_*>0$ and $\beta_*>0$, which may depend on $\chi$, such that for each $\lambda\in\C$ with $-\frac{\tilde{\lambda}}{2}\leq \Re \,\lambda \leq \chi$ and $|\Im \,\lambda|\leq\frac{\pi |c|}{h}$ we have the bound
\begin{equation}\label{ditiseenbound}
\begin{array}{lcl}
|G_{\infty;\lambda}(\xi-\xi_0)|\leq K_* e^{-\beta_*|\xi-\xi_0|}
\end{array}
\end{equation}
for all $\xi,\xi_0\in\R$.\end{lemma}

Pick $\lambda\in\C\setminus\sigma(-L)$ with $\Re \,\lambda\geq-\frac{\tilde{\lambda}}{2}$. Observe that
\begin{equation}\begin{array}{lcl}L-L_{\infty}&=&\left(\begin{array}{ll}-g_u(\overline{u})+r_0& 0\\ 0&0\end{array}\right).\end{array}\end{equation}
We know that $G_{\infty;\lambda}(\cdot-\xi_0)\in L^2(\R,\text{Mat}_2(\R))$ since it decays exponentially.
This means that we also have the inclusion
\begin{equation}\begin{array}{lcl}[L-L_{\infty}]G_{\infty;\lambda}(\cdot-\xi_0)&\in & L^2(\R,\text{Mat}_2(\C)).\end{array}\end{equation}
Hence, it is possible to define the function $G_{\lambda}$ by writing
\begin{equation}\label{definitievanGhlambda}\begin{array}{lcl}G_{\lambda}(\xi,\xi_0)&=&G_{\infty;\lambda}(\xi-\xi_0)-\Big[(\lambda+L)^{-1}[L-L_{\infty}]G_{\infty;\lambda}(\cdot-\xi_0)\Big](\xi).\end{array}\end{equation}
The next result shows that $G_{\lambda}$ can be interpreted as the Green's function of $L+\lambda$. It is based on \cite[Lemma 2.6]{HJHSTBFHN}.
\begin{lemma}\label{existenceofgreen} Assume that (\asref{aannamespuls}{\text{H}}), (\asref{extraaannamespuls}{\text{H}}), (\asref{aannamesconstanten}{\text{H}}), (\asref{aannames}{\text{H}\alpha1}) and (\asref{extraaannames}{\text{H}\alpha}) are satisfied. For $\lambda\in\C\setminus\sigma(-L)$ with $\Re \,\lambda\geq-\frac{\tilde{\lambda}}{2}$ we have that $G_{\lambda}(\cdot,y)$ is continuous on $\R\setminus\{y\}$ and $C^1$-smooth on $\R\setminus\{y+kh:k\in\Z\}$. Furthermore, it satisfies
\begin{equation}\label{integraalGheigenschap}\begin{array}{lcl}\int_{-\infty}^\infty \Big[(\lambda+L) G_{\lambda} (\cdot,\xi_0)\Big](\xi)f(\xi)\ d\xi&=&f(\xi_0)\end{array}\end{equation}
for all $\xi\in\R$ and all $f\in \textbf{H}^1$.
\end{lemma}

The link between our two types of Green's functions is provided by the following key result. It is based on \cite[Theorem 4.2]{BGV2003}, where it was used to study one-sided spatial discretisation schemes for systems with conservation laws.
\begin{proposition}\label{thrGreenfunction} Assume that (\asref{aannamespuls}{\text{H}}),(\asref{extraaannamespuls}{\text{H}}), (\asref{aannamesconstanten}{\text{H}}), (\asref{aannames}{\text{H}\alpha1}) and (\asref{extraaannames}{\text{H}\alpha}) are satisfied. Let $\chi>\lambda_{\mathrm{unif}}$ be given, where $\lambda_{\mathrm{unif}}$ is as in Lemma \ref{lemmaexistenceGreen}. For all $t\geq t_0$ the Green's function $\mathcal{G}_j^{j_0}(t,t_0)$ of (\ref{3.2equivalent}) is given by
\begin{equation}\label{4.3equivalent}\begin{array}{lcl}\mathcal{G}_j^{j_0}(t,t_0)&=&-\frac{h}{2\pi i}\int\limits_{\chi-\frac{i\pi c}{h}}^{\chi+\frac{i\pi c}{h}}e^{\lambda(t-t_0)}G_{\lambda}(hj+c t,hj_0+c t_0)d\lambda\end{array}\end{equation}
where $G_{\lambda}$ is the Green's function of $\lambda+L$ as defined in (\ref{definitievanGhlambda}).\end{proposition}

Our first task is to collect
several basic facts concerning the operators
$L_h$ and $L_\infty$ that will allow us
to establish Lemma's \ref{eigenschappenghinfty} and \ref{existenceofgreen}. In particular,
we need to isolate and explicitly compute
the part of the Fourier integral
(\ref{defvanghinfty})
that behave as $\abs{\eta}^{-1}$ and $\abs{\eta}^{-2}$
as $\eta \to \pm \infty$, as these lead
to the discontinuities in $G_{\infty;\lambda}$ and
its derivative.

\begin{lemma}\label{deltaremainscont} Assume that (\asref{aannames}{\text{H}\alpha1}) and (\asref{extraaannames}{\text{H}\alpha}) are satisfied. Consider any bounded function $f:\R\rightarrow\R$ which is continuous everywhere except at some $\xi_0\in\R$. Then $\Delta_h f$ is continuous everywhere except at $\{\xi_0+hk:k\in\Z\}$. Moreover, if $f$ is differentiable except at $\xi_0$ and $f'$ is bounded, then $\Delta_h f$ is differentiable everywhere except at $\{\xi_0+hk:k\in\Z\}$ and $[\Delta_h f]'(\xi)=[\Delta_h f'](\xi)$. \end{lemma}
\vspace*{4pt}\noindent\textit{Proof.} For convenience we set $\xi_0=0$. Pick $\xi\in\R$ with $\xi\notin\{kh:k\in\Z\}$. Then $f$ is continuous in each point $\xi+kh$ for $k\in\Z$. Choose $\epsilon>0$. Since $f$ is bounded and $\sum\limits_{j=1}^\infty |\alpha_j|<\infty$, we can pick $K>0$ in such a way that
\begin{equation}\begin{array}{lcl}2\nrm{ f}_\infty\frac{1}{h^2}\sum\limits_{j=K}^\infty|\alpha_j|&<&\frac{\epsilon}{2}.\end{array}\end{equation}
For $j\in\{1,...,K-1\}$ we can pick $\delta_j>0$ in such a way that
\begin{equation}\begin{array}{lcl}\frac{1}{h^2}|\alpha_j|\Big|f(\xi+y+hj)-f(\xi+hj)\Big|&<&\frac{\epsilon}{2^{K+1}}\end{array}\end{equation}
for all $y\in\R$ with $|y|<\delta_j$. Let $\delta=\min\{\delta_j:1\leq j<K\}>0$. Then for $y\in\R$ with $|y|<\delta$ we obtain
\begin{equation}\begin{array}{lcl}|\Delta_h f(\xi+y)
-\Delta_h f(\xi)|&\leq&\frac{1}{h^2}\sum\limits_{j=K}^\infty|\alpha_j|\Big(|f(\xi+y+jh)|+|f(\xi+jh)|\Big)\\[0.2cm]
&&\qquad +\frac{1}{h^2}\sum\limits_{j=1}^{K-1}|\alpha_j|\Big|f(\xi+y+jh)-f(\xi+jh)\Big|\\[0.2cm]
&\leq &\frac{2}{h^2}\sum\limits_{j=K}^\infty|\alpha_j|\nrm{f}_\infty+\sum\limits_{j=1}^{K-1}\frac{\epsilon}{2^{K+1}}\\[0.2cm]
&<&\frac{\epsilon}{2}+\frac{\epsilon}{2}\\[0.2cm]
&=&\epsilon.\end{array}\end{equation}
So $\Delta_h f$ is continuous outside of $\{kh:k\in\Z\}$.\par

Writing
\begin{equation}
\begin{array}{lcl}
f_n(\xi)&=&\frac{1}{h^2}\sum\limits_{j=1}^n\alpha_j\Big[f(\xi+hj)+f(\xi-hj)-2f(\xi)\Big]
\end{array}
\end{equation}
for $n\in\Z_{>0}$, we can compute
\begin{equation}
\begin{array}{lcl}
f_n'(\xi)&=&\frac{1}{h^2}\sum\limits_{j=1}^n\alpha_j\Big[f'(\xi+hj)+f'(\xi-hj)-2f'(\xi)\Big].
\end{array}
\end{equation}
This allows us to estimate
\begin{equation}
\begin{array}{lcl}
|f_n'(\xi)-(\Delta_hf')(\xi)|&\leq &\frac{1}{h^2}\sum\limits_{j=n+1}^\infty|\alpha_j|4\nrm{f'}_\infty.
\end{array}
\end{equation}
In particular, the sequence $\{f_n'\}$ converges uniformly to $\Delta_hf'$ from which it follows that
\begin{equation}
\begin{array}{lcl}
(\Delta_h f)'(\xi)&=&\frac{1}{h^2}\sum\limits_{j=1}^\infty\alpha_m\Big[f'(\xi+hj)+f'(\xi-hj)-2f'(\xi)\Big]\\[0.2cm]
&=& (\Delta_h f')(\xi).
\end{array}
\end{equation}\qed

\vspace*{4pt}\noindent\textit{Proof of Lemma \ref{eigenschappenghinfty}.}
Pick $\chi>0$ and set
\begin{equation}
\begin{array}{lcl}
R&=&\{\lambda\in\C:-\frac{\tilde{\lambda}}{2}\leq \Re \,\lambda \leq \chi\hbox{ and }|\Im \,\lambda|\leq\frac{\pi |c|}{h}\}.
\end{array}
\end{equation}
The proof of Lemma \ref{hyperbolic2} implies that we can choose $\beta_*>0$ and $K_*>0$ in such a way that
\begin{equation}
\begin{array}{lcl}
\nrm{\Delta_{L_{\infty;\lambda}}(z)^{-1}}&\leq &\frac{K_*}{1+|\Im\,z|}
\end{array}
\end{equation}
for all $\lambda\in R$ and all $z\in\C$ with $|\Re\,z|\leq 2\beta_*$. In particular, it follows that
$(y\mapsto\Delta_{L_{\infty;\lambda}}(iy)^{-1})\in L^2(\R)$. By the Plancherel Theorem it follows that
$G_{\infty;\lambda}$ is a well-defined function in $L^2(\R)$. In particular, it is bounded.
Shifting the integration path in (\ref{defvanghinfty}) in the standard fashion described in \cite{HJHCM,MPA},
we obtain the bound
\begin{equation}
\begin{array}{lcl}
|G_{\infty;\lambda}(\xi-\xi_0)|\leq K_* e^{-\beta_*|\xi-\xi_0|}
\end{array}
\end{equation}
for all $\xi,\xi_0\in\R$ and $\lambda\in R$.
\par

We loosely follow the approach of \cite[\S 5.1]{HJHPOLLUTION}, which considers a similar setting for Green's functions for Banach space-valued operators with finite range interactions. Pick $\lambda\in R$. We rewrite the definition of $\Delta_{L_{\infty;\lambda}}$ given in (\ref{definitiedeltalh}) in the more general form
\begin{equation}
\begin{array}{lcl}
\frac{1}{c}\Delta_{L_{\infty;\lambda}}(z)&=&
z-B_{\infty;\lambda}e^{z\cdot},
\end{array}
\end{equation}
For $\alpha\in\R$ close to $0$ we
introduce the expression $\mathcal{R}_{L_{\infty;\lambda};\alpha}$ by
\begin{equation}\label{eq5.23poll}
\begin{array}{lcl}
\mathcal{R}_{L_{\infty;\lambda};\alpha}(z)&=& c \Delta_{L_{\infty;\lambda}}(z)^{-1}-\frac{1}{z-\alpha}
 -\frac{B_{\infty;\lambda}e^{z\cdot}-\alpha}{(z-\alpha)^2}
\end{array}
\end{equation}
for $z\in\C$ unequal to $\alpha$ and $|\Re\,z|\leq 2\beta_*$. Since we can compute
\begin{equation}
\begin{array}{lcl}
c 
\Delta_{L_{\infty;\lambda}}(z)^{-1}&=&\Big[z-\alpha+\big(\alpha-B_{\infty;\lambda}e^{z\cdot}\big)\Big]^{-1}\\[0.2cm]
&=&(z-\alpha)^{-1}\Big[1+(z-\alpha)^{-1}\big(\alpha-B_{\infty;\lambda}e^{z\cdot}\big)\Big]^{-1}\\[0.2cm]
&=&(z-\alpha)^{-1}\Big[1-(z-\alpha)^{-1}\big(\alpha-B_{\infty;\lambda}e^{z\cdot}\big)
  +\mathcal{O}((z-\alpha)^{-2})\Big] ,
\end{array}
\end{equation}
we obtain the estimate 
\begin{equation}
\begin{array}{lcl}
|\mathcal{R}_{L_{\infty;\lambda};\alpha}(iy)|&\leq &\frac{K_*}{1+|y|^3},
\end{array}
\end{equation}
for all $y \in \mathbb{R}$, possibly after increasing $K_*$. \par 

Exploiting the decomposition (\ref{eq5.23poll}), we write
\begin{equation}
\begin{array}{lcl}
G_{\infty;\lambda}&=& \frac{1}{c}
  \mathcal{M}_\alpha
  +\frac{1}{c}
  \mathcal{R}_\alpha,
\end{array}
\end{equation}
where we have introduced
\begin{equation}
\begin{array}{lcl}
\mathcal{M}_\alpha(\xi)&=&\frac{1}{2\pi}\int_{-\infty}^\infty e^{i\eta\xi}\Big(\frac{1}{i\eta-\alpha}
  -\frac{B_{\infty;\lambda}e^{i\eta\cdot}-\alpha}{(i\eta-\alpha)^2}\Big)\ d\eta,\\[0.2cm]
\mathcal{R}_\alpha(\xi)&=&\frac{1}{2\pi}\int_{-\infty}^\infty e^{i\eta\xi}\mathcal{R}_{L_{\infty;\lambda};\alpha}(i\eta)\ d\eta
\end{array}
\end{equation}
for any $\alpha\in\R\setminus\{0\}$ and $\xi\in\R$.
Using \cite[Lemma 5.8]{HJHPOLLUTION} we can explicitly compute
\begin{equation}
\begin{array}{lcl}
\mathcal{M}_\alpha (\xi)&=&-e^{\alpha\xi}H(-\xi)-\Big[B_{\infty;\lambda}-\alpha\Big]\Big( \cdot e^{\alpha\cdot}H(-\cdot)\Big)(\xi),
\end{array}
\end{equation}
where we have introduced the Heaviside function $H$ as
\begin{equation}
\begin{array}{lcl}
H(\xi)&=&\begin{cases} I,\enskip \xi>0\\ \frac{1}{2} I,\enskip \xi=0\\ 0,\enskip \xi<0.\end{cases}
\end{array}
\end{equation}
Since $\xi\mapsto \xi e^{\alpha\xi}H(-\xi)$ is continuous everywhere and differentiable outside of $\xi=0$,
Lemma \ref{deltaremainscont} implies that $\mathcal{M}_\alpha$ is continuous everywhere outside of $\xi=0$
and differentiable outside of $\{hk:k\in\Z\}$. Moreover, we have the jump discontinuity
\begin{equation}\label{discontinuity}
\begin{array}{lcl}
\mathcal{M}_\alpha(0^+)-\mathcal{M}_{\alpha}(0^-)&=& I
\end{array}
\end{equation}
and we can easily compute
\begin{equation}
\begin{array}{lcl}
\mathcal{M}_\alpha'(\xi)&=&\alpha\mathcal{M}_\alpha(\xi)-[B_{\infty;\lambda}-\alpha]\Big[e^{\alpha\cdot}H(-\cdot)\Big](\xi),
\end{array}
\end{equation}
from which it follows that
\begin{equation*}
\begin{array}{lcl}
\frac{1}{c} L_{\infty;\lambda}\mathcal{M}_\alpha(\xi)&=&\mathcal{M}_\alpha'(\xi)-
  B_{\infty;\lambda}\mathcal{M}_\alpha(\xi)\\[0.2cm]
&=&-\alpha e^{\alpha\xi}H(-\xi)-\alpha\Big[B_{\infty;\lambda}-\alpha\Big]
  \Big( \cdot e^{\alpha\cdot}H(-\cdot)\Big)(\xi)\\[0.2cm]
&&\qquad -[B_{\infty;\lambda}-\alpha]
    \Big[e^{\alpha\cdot}H(-\cdot)\Big](\xi)+B_{\infty;\lambda}\Big[e^{\alpha\cdot}H(-\cdot)\Big](\xi)
\end{array}
\end{equation*}\begin{equation}
\begin{array}{lcl}&&\qquad+B_{\infty;\lambda}\Big[[B_{\infty;\lambda}-\alpha]\Big( \cdot e^{\alpha\cdot}H(-\cdot)\Big)(\ast)\Big](\xi)\\[0.2cm]
&=&[B_{\infty;\lambda}-\alpha]\Big[[B_{\infty;\lambda}-\alpha]\Big( \cdot e^{\alpha\cdot}H(-\cdot)\Big)(\ast)\Big](\xi).
\end{array}
\end{equation}
Since $\mathcal{R}_{L_{\infty;\lambda};\alpha}\in L^1(\R)$ we see that $\mathcal{R}_\alpha$ is continuous. Therefore, $G_{\infty;\lambda}$ is continuous outside of $\xi=0$. Similarly to \cite[Eq. (5.79)]{HJHPOLLUTION} we observe that
\begin{equation}
\begin{array}{lcl}
\frac{1}{c} \Delta_{L_{\infty;\lambda}}(z)\mathcal{R}_{L_{\infty;\lambda};\alpha}(z)&=&
\frac{(B_{\infty;\lambda}e^{z\cdot}-\alpha)^2}{(z-\alpha)^2},
\end{array}
\end{equation}
which yields
\begin{equation}
\begin{array}{lcl}
\frac{1}{c} L_{\infty;\lambda}\mathcal{R}_\alpha(\xi)&=&\mathcal{R}_\alpha'(\xi)
 -B_{\infty;\lambda}\mathcal{R}_\alpha(\xi)\\[0.2cm]
&=&\frac{1}{2\pi c}\int_{-\infty}^{\infty}e^{i \xi y}\Delta_{L_{\infty;\lambda}}(iy)
   \mathcal{R}_{L_{\infty;\lambda};\alpha}(iy)dy\\[0.2cm]
&=&\frac{1}{2\pi }\int_{-\infty}^{\infty}e^{i \xi y}\frac{(B_{\infty;\lambda}e^{iy\cdot}-\alpha)^2}{(iy-\alpha)^2}dy\\[0.2cm]
&=&-[B_{\infty;\lambda}-\alpha]\Big[[B_{\infty;\lambda}-\alpha]\Big( \cdot e^{\alpha\cdot}H(-\cdot)\Big)(\ast)\Big](\xi),
\end{array}
\end{equation}
using \cite[Lemma 5.8]{HJHPOLLUTION}. In particular, we see that
\begin{equation}
\begin{array}{lcl}
L_{\infty;\lambda}G_{\infty;\lambda}(\xi)&=&0
\end{array}
\end{equation}
for all $\xi$ outside of $\{hk:k\in\Z\}$. Lemma \ref{deltaremainscont} subsequently
shows that $G_{\infty;\lambda}$ is $C^1$-smooth outside of $\{hk:k\in\Z\}$. \par

Fix $f\in\mathbf{H}^1$. For any $\delta>0$ we may compute
\begin{equation}
\begin{array}{lcl}
0&=& \int_\delta^\infty \Big[L_{\infty;\lambda}G_{\infty;\lambda}(\cdot)\Big](\xi)f(\xi)d\xi\\[0.2cm]
&=&\Big[ c G_{\infty;\lambda}f\Big]_{\delta}^\infty-\int_\delta^\infty  c G_{\infty;\lambda}(\xi)f'(\xi)+[cB_{\infty;\lambda}G_{\infty;\lambda}](\xi)f(\xi),
\end{array}
\end{equation}
together with
\begin{equation}
\begin{array}{lcl}
0&=&\Big[ c G_{\infty;\lambda}f\Big]_{-\infty}^{-\delta}-\int_{-\infty}^{-\delta} c G_{\infty;\lambda}(\xi)f'(\xi)
  +[c B_{\infty;\lambda}G_{\infty;\lambda}](\xi)f(\xi).
\end{array}
\end{equation}
Using (\ref{discontinuity}) we can hence compute
\begin{equation}
\begin{array}{lcl}
\int_{-\infty}^\infty \Big[L_{\infty;\lambda}G_{\infty;\lambda}(\cdot)\Big](\xi)f(\xi)d\xi&=&\lim_{\delta\downarrow 0}\Big[c G_{\infty;\lambda}f\Big]_{\delta}^\infty -\Big[c G_{\infty;\lambda}f\Big]_{-\infty}^{-\delta}\\[0.2cm]
&=&f(0).
\end{array}
\end{equation}
%
\qed

\begin{lemma}\label{lemma4.3bachbos} Assume that (\asref{aannamespuls}{\text{H}}),(\asref{extraaannamespuls}{\text{H}}), (\asref{aannamesconstanten}{\text{H}}), (\asref{aannames}{\text{H}\alpha1}) and (\asref{extraaannames}{\text{H}\alpha}) are satisfied. Fix $\lambda\in\C$ with $\Re \,\lambda\geq-\frac{\tilde{\lambda}}{2}$. Then there exist constants $K>0$ and $\beta>0$ so that
for any $g\in\mathbf{L}^2$ and $f\in\mathbf{H}^1$
that satisfy  $(L+\lambda)f=g$,
the pointwise bound
\begin{equation}
\begin{array}{lcl}
|f(\xi)|&\leq &Ke^{-\alpha|\xi|}\nrm{f}_\infty+K\int_{-\infty}^\infty e^{-|\eta-\xi|}g(\eta)d\eta
\end{array}
\end{equation}
holds for all $\xi\in\R$.\end{lemma}
\vspace*{4pt}\noindent\textit{Proof.} On account of Lemma \ref{eigenschappenghinfty} we can lift the results from
\cite[Prop. 5.2-5.3]{MPA} to our current infinite range setting.
The proof of these results are identical, since the estimate \cite[Eq. (5.4)]{MPA} still holds in our
setting on account of (\asref{extraaannames}{\text{H}\alpha}).
A more detailed description for this
procedure can be found in \cite[Lemma 4.1-Lemma 4.3]{BachBos}. \qed

\vspace*{4pt}\noindent\textit{Proof of Lemma \ref{existenceofgreen}.}
Pick $\lambda\in\C\setminus\sigma(-L)$ and compute
\begin{equation}\begin{array}{lcl}(\lambda+L)G_{\lambda}(\cdot,\xi_0)&=&(\lambda+L)G_{\infty;\lambda}(\cdot-\xi_0)-[L-L_{\infty}]G_{\infty;\lambda}(\cdot -\xi_0)\\[0.2cm]
&=&(\lambda+L_{\infty} )G_{\infty;\lambda}(\cdot-\xi_0).\end{array}\end{equation}
The last statement follows immediately from this identity.

Write
\begin{equation}\begin{array}{lcl}\hat{G}_{\infty;\lambda}(\cdot-\xi_0)&=&[L-L_{\infty}]G_{\infty;\lambda}(\cdot-\xi_0).\end{array}\end{equation}
We have already seen that $\hat{G}_{\infty;\lambda}(\cdot-\xi_0)\in L^2(\R,\text{Mat}_2(\C))$. Hence, it follows that
\begin{equation}\begin{array}{lcl}(\lambda+L)^{-1}\hat{G}_{\infty;\lambda}(\cdot-\xi_0)&\in & H^1(\R,\text{Mat}_2(\C)).\end{array}\end{equation}
In particular, this function is continuous. Together with Lemma \ref{eigenschappenghinfty} we obtain that $G_{\lambda}(\cdot,\xi_0)$ is continuous on $\R\setminus\{\xi_0\}$.\par

Set $H=(\lambda+L)^{-1}\hat{G}_{\infty;\lambda}$ and write $H=\left(\begin{array}{ll}H^{(1,1)}&H^{(1,2)}\\ H^{(2,1)}& H^{(2,2)}\end{array}\right)$. Using the definition of $L$ we see that
\begin{equation}\begin{array}{lcl}c \frac{d}{d\xi} H&=&-\lambda H-\hat{G}_{\infty}-\tilde{H},\end{array}\end{equation}
where
\begin{equation}\begin{array}{lcl}\tilde{H}&=&-\left(\begin{array}{ll}-\Delta_h H^{(1,1)}-g_u(\overline{u})H^{(1,1)}+H^{(2,1)} & \Delta_h H^{(1,2)}-g_u(\overline{u})H^{(1,2)}+H^{(2,2)}\\ -\rho H^{(1,1)}+\gamma\rho H^{(2,1)}&-\rho H^{(1,2)}+\gamma\rho H^{(2,2)} \end{array}\right).\end{array}\end{equation}
Since $\overline{u}'\in H^1$ and hence $\overline{u}'$ is continuous, we must have that $\overline{u}$ is continuous. As argued before $\Delta_h H^{(1,1)}$ and $\Delta_h H^{(1,2)}$ are also continuous. Hence, we see that $c \frac{d}{d\xi} H$ is continuous on $\R\setminus\{\xi_0\}$ and thus that $\frac{d}{d\xi}H$ is continuous on $\R\setminus\{\xi_0\}$. Therefore, we obtain that $G_{\lambda}(\cdot,\xi_0)$ is $C^1$-smooth on $\R\setminus\{\xi_0+kh:k\in\Z\}$.\qed

We now proceed to the verification
of the integral representation (\ref{4.3equivalent}).
As a preparation, we need to show that whenever
$\lambda$ has a  sufficiently large real part, the function $G_{\lambda}$ is bounded uniformly by a constant. This result is based on \cite[Lemma 4.1]{BGV2003}.
\begin{lemma}\label{lemmaexistenceGreen} Assume that (\asref{aannamespuls}{\text{H}}),(\asref{extraaannamespuls}{\text{H}}), (\asref{aannamesconstanten}{\text{H}}), (\asref{aannames}{\text{H}\alpha1}) and (\asref{extraaannames}{\text{H}\alpha}) are satisfied. Then there exist constants $K$ and $\lambda_{\mathrm{unif}}$ so that the Green's function $G_{\lambda}$ enjoys the uniform estimate
\begin{equation}\label{4.6equivalent}\begin{array}{lcl}|G_{\lambda}(\xi,\xi_0)|&\leq& K,\end{array}\end{equation}
for all $\xi,\xi_0\in\R$, whenever $\Re \,\lambda>\lambda_\mathrm{unif}$.
\end{lemma}
\vspace*{4pt}\noindent\textit{Proof.} We write $L=c\frac{d}{d\xi}+B$ with
\begin{equation}\begin{array}{lcl}B&=&\left(\begin{array}{ll}-\Delta_h -g_u(\overline{u})& 1\\ -\rho & \gamma\rho\end{array}\right).\end{array}\end{equation}
We introduce $G_{\lambda}^0$ as the Green's function of $(\lambda+c\frac{d}{d\xi})$ viewed as a map from $\textbf{H}^1$ to $\textbf{L}^2$. Luckily, it is well-known that this Green's function admits the estimate
\begin{equation}\begin{array}{lcl}|G_{\lambda}^0(\xi,\xi_0)|&\leq&\frac{1}{|c|}e^{-\Re \,\lambda|\xi-\xi_0|/|c|}.\end{array}\end{equation}
We can look for the Green's function $G_{\lambda}$ as the solution of the fixed point problem
\begin{equation}\label{4.8equivalent}\begin{array}{lcl}G_{\lambda}(\xi,\xi_0)&=&G_{\lambda}^0(\xi,\xi_0)+\int_\R G_{\lambda} (\xi,z)(BG_{\lambda}^0)(z,\xi_0)dz.\end{array}\end{equation}
Since $\lambda+L$ is invertible by Theorem \ref{totalespectrum}, $G_{\lambda}$ must necessarily satisfy the fixed point problem (\ref{4.8equivalent}).\par

For a matrix $A\in\text{Mat}_2(\C)$ we write $A=\left(\begin{array}{ll}A^{(1,1)}&A^{(1,2)}\\A^{(2,1)}&A^{(2,2)}\end{array}\right)$. We make the decomposition
\begin{equation}\begin{array}{lcl}

B&=&B_0+B_1,\end{array}
\end{equation}
where
\begin{equation}\begin{array}{lcl}
B_0&=&\left(\begin{array}{ll}-\Delta_h & 0\\ 0 & 0\end{array}\right),\\[0.4cm]
B_1&=&\left(\begin{array}{ll} -g_u(\overline{u})& 1\\ -\rho & \gamma\rho\end{array}\right).

\end{array}
\end{equation}
We estimate
\begin{equation}
\begin{array}{lcl}
|(B_0 G_{\lambda}^0)(\xi,\xi_0)|&= &|\Delta_h G_{\lambda}^0(\xi,\xi_0)^{(1,1)}|\\[0.2cm]
&\leq &\sum\limits_{j=1}^\infty\Big[ \frac{1}{h^2}|\alpha_j|\Big(|G_{\lambda}^0(\xi+hj,\xi_0)^{(1,1)}|+|G_{\lambda}^0(\xi-hj,\xi_0)^{(1,1)}|+2|G_{\lambda}^0(\xi,\xi_0)^{(1,1)}|\Big)\Big]\\[0.2cm]
&\leq &\frac{1}{|c|}\sum\limits_{j=1}^\infty\Big[\frac{1}{h^2}|\alpha_j|\Big(e^{-\Re \,\lambda|\xi+hj-\xi_0|/|c|} +e^{-\Re \,\lambda|\xi-hj-\xi_0|/|c|}+2e^{-\Re \,\lambda|\xi-\xi_0|/|c|}\Big)\Big]
\end{array}
\end{equation}
and observe that
\begin{equation}
\begin{array}{lcl}
\int_\R|(B_0 G_{\lambda}^0)(\xi,\xi_0) |\ d\xi&\leq &\frac{1}{|c|}\Bigg(\sum\limits_{j=1}^\infty 4\Big[\frac{1}{h^2}|\alpha_j|\frac{1}{\Re \,\lambda/|c|}\Big]\Bigg)\\[0.2cm]
&=&\frac{4}{h^2\Re \,\lambda}\sum\limits_{j=1}^\infty|\alpha_j|.\end{array}\end{equation}
We now fix $G\in L^\infty(\R^2,\text{Mat}_2(\C))$ and consider the expressions
\begin{equation}
\begin{array}{lcl}
\mathcal{I}_0&=&\int_\R \Big[G(\xi,z)(B_0G_{\lambda}^0)(z,\xi_0)\Big]^{(1,1)}dz,\\[0.2cm]
\mathcal{I}_1&=&\int_\R \Big[G(\xi,z)(B_1G_{\lambda}^0)(z,\xi_0)\Big]^{(1,1)}dz.
\end{array}
\end{equation}
Using Fubini's theorem for positive functions to switch the integral and the sum, we obtain the estimates
\begin{equation}\begin{array}{lcl}|\mathcal{I}_0|
&\leq &\nrm{G}_{L^\infty}\int_\R|(B_0 G_{\lambda}^0)(z,\xi_0) |\ dz\\[0.2cm]
&\leq &\nrm{G}_{L^\infty}\frac{4}{h^2\Re \,\lambda}\sum\limits_{j=1}^\infty|\alpha_j|\end{array}\end{equation}
and
\begin{equation}\begin{array}{lcl}|\mathcal{I}_1|
&\leq &\nrm{G}_{L^\infty}\int_{\R}\Big(|g_u(\overline{u}(z))||G_{\lambda}^0(z,\xi_0)^{(1,1)}|+\rho|G_{\lambda}^0(z,\xi_0)^{(1,1)}|\\[0.2cm]
&&\qquad +(1+\gamma\rho)|G_{\lambda}^0(z,\xi_0)^{(2,1)}|\Big)dz\\[0.2cm]
&\leq &\nrm{G}_{L^\infty}\frac{1}{|c|}\int_{\R}\Big(\Big(|g_u(\overline{u}(z))|+\rho+1+\gamma\rho\Big)e^{-\Re \,\lambda|z-\xi_0|/|c|}\Big)dz\\[0.2cm]
&\leq &\nrm{G}_{L^\infty}\frac{1}{|c|}\Big(\nrm{g_u(\overline{u})}_{L^\infty}+\rho+1+\gamma\rho\Big)\left(\frac{1}{\Re \,\lambda/|c|}\right)\\[0.2cm]
&\leq &\nrm{G}_{L^\infty}\Big(g_*+\rho+1+\gamma\rho\Big)\left(\frac{1}{\Re \,\lambda}\right).\end{array}\end{equation}
Similar estimates hold for the other components of $\int_\R G(\xi,z)(BG_{\lambda}^0)(z,\xi_0)dz$. Therefore, the mapping $G \mapsto \int_\R G(\xi,z)(BG_{\lambda}^0)(z,\xi_0)dz$ is a contraction in $L^\infty(\R^2,\text{Mat}_2(\C))$ for $\Re \,\lambda>\lambda_{\mathrm{unif}}$ for $\lambda_{\mathrm{unif}}$ large enough, with $\lambda_{\mathrm{unif}}$ possibly dependent of $h\in(0,h_{**})$. Hence, we get a unique bounded solution of (\ref{4.8equivalent}), which must be $G_{\lambda}$. The desired bound on $G_{\lambda}$ is now immediate.\qed

\vspace*{4pt}\noindent\textit{Proof of Proposition \ref{thrGreenfunction}.} Fix $j_0\in\Z$ and $t_0\in\R$. Since (\ref{3.2equivalent}) is merely a linear ODE in the Banach space $\mathbb{L}^2$, it follows from the Cauchy-Lipschitz theorem that (\ref{3.2equivalent}) indeed has a unique solution $\mathcal{V}:[t_0,\infty)\rightarrow\mathbb{L}^2$. For any $Z\in C_c^\infty (\R;\mathbb{L}^2)$, an integration by parts yields
\begin{equation}\begin{array}{lcl}-Z_{j_0}(t_0)&=&\int_{t_0}^\infty\sum\limits_{j\in\Z} \Big[\Big(\frac{d\mathcal{V}_j}{dt}(t)-(\mathcal{A}(t)\cdot \mathcal{V}(t))_j\Big)Z_j(t)\Big]dt-\sum_{j\in\Z}\mathcal{V}_j(t_0)Z_j(t_0)\\[0.2cm]
&=&\int_{t_0}^\infty \sum\limits_{j\in\Z}\Big[-\frac{dZ_j}{dt}(t)\mathcal{V}_j(t)-(\mathcal{A}(t)\cdot \mathcal{V})_j(t)Z_j(t)\Big]dt
 .\end{array}\end{equation}

We want to show that the function $V_j(t):=\mathcal{G}_j^{j_0}(t,t_0)$ defined by (\ref{4.3equivalent}) coincides with $\mathcal{V}$ on $[t_0,\infty)$. To accomplish this, we define
\begin{equation}\begin{array}{lcl}I&=&\int_{t_0}^\infty \sum\limits_{j\in\Z}\Big[-\frac{dZ_j}{dt}(t)V_j(t)-(\mathcal{A}(t)\cdot V(t))_jZ_j(t)\Big]dt\end{array}\end{equation}
and show that $V$ is a weak solution to (\ref{3.2equivalent}) in the sense that
\begin{equation}\begin{array}{lcl}I&=& -Z_{j_0}(t_0)\end{array}\end{equation}
holds for all $Z \in C_c^\infty (\R;\mathbb{L}^2)$. Indeed, the uniqueness of weak solutions then implies that $V=\mathcal{V}$.\par

Note first that $V(t)=0$ for $t<t_0$, which can be seen by using (\ref{4.6equivalent}) and taking $\chi\rightarrow \infty$ in (\ref{4.3equivalent}). We write $y=hj_0+ct_0$, $\chi_-=\chi-\frac{i\pi c}{h}$ and $\chi_+=\chi+\frac{i\pi c}{h}$. We see that
\begin{equation}\begin{array}{lcl}I
&=&\int\limits_{-\infty}^\infty \sum\limits_{j\in\Z}\Big[-\frac{dZ_j}{dt}(t)V_j(t)-(\mathcal{A}(t)\cdot V(t))_jZ_j(t)\Big]dt,\end{array}\end{equation}
since $V(t)=0$ for $t<t_0$. Moreover, we write
\begin{equation}
\begin{array}{lcl}
\mathbb{G}_j(t)&=&G_{\lambda}(hj+c t,y).
\end{array}
\end{equation}
Using our definition of $V(t)$, we have
\begin{equation}\label{watisI}\begin{array}{lcl}I&=&-\frac{h}{2\pi i}\int\limits_{\chi_-}^{\chi_+}\sum\limits_{j\in\Z}\bigg[\int\limits_{-\infty}^\infty \mathcal{I}_{j}(t,\lambda)dt \bigg]d\lambda,\end{array}\end{equation}
where
\begin{equation}
\begin{array}{lcl}
\mathcal{I}_{j}(t,\lambda)&=&e^{\lambda(t-t_0)}\Big[-\mathbb{G}_j(t)\frac{dZ_j}{dt}(t)-(\mathcal{A}(t)\cdot \mathbb{G}(t))_jZ_j(t)\Big].\end{array}\end{equation}
The permutation of the summations and integrations is allowed by Lebesgue's theorem because $Z$ and $\frac{dZ}{dt}$ are compactly supported and $G_{\lambda}$ is uniformly bounded by (\ref{4.6equivalent}). Fix $\chi_-\leq\lambda\leq\chi_+$ and $j\in\Z$. Using the change of variable $x=hj+ct$ we obtain
\begin{equation}\begin{array}{lcl}\int\limits_{-\infty}^\infty \mathcal{I}_{j}(t,\lambda)dt&=&\frac{1}{c}\int\limits_{x=-\infty}^{x=\infty}\Big[-c\mathbb{G}_j\Big(\frac{x-hj}{c}\Big)\frac{d\mathcal{Z}_j}{dx}+\lambda\mathbb{G}_j\Big(\frac{x-hj}{c}\Big)\mathcal{Z}_j(x,\lambda)\\[0.2cm]
&&\qquad -\Big(\mathcal{A}\Big(\frac{x-hj}{c}\Big)\cdot \mathbb{G}\Big(\frac{x-hj}{c}\Big)\Big)_j\mathcal{Z}_j(x,\lambda)\Big]dx,\end{array}\end{equation}
where
\begin{equation}\begin{array}{lcl}\mathcal{Z}_j(x,\lambda)&=&e^{\lambda((x-hj)/c-t_0)}Z_j\Big(\frac{x-hj}{c}\Big).\end{array}\end{equation}
Exploiting the fact that $Z_j$ and, therefore, $\mathcal{Z}_j$ is compactly supported, (\ref{integraalGheigenschap}) yields
\begin{equation}\begin{array}{lcl}\int\limits_{-\infty}^\infty \mathcal{I}_{j}(t,\lambda)dt&=&\frac{1}{c}\int_{-\infty}^\infty [(L+\lambda)G_{\lambda}(x,y)\mathcal{Z}_j(x,\lambda)]dx\\[0.2cm]
&=&\frac{1}{c}\mathcal{Z}_j(y).\end{array}\end{equation}
Now since $Z_j$ is compactly supported, we can exchange sums and integrals in equation (\ref{watisI}). This allows us to compute
\begin{equation}\begin{array}{lcl}I&=&-\frac{h}{2\pi i}\frac{1}{c}\int\limits_{\chi_-}^{\chi_+}\sum\limits_{j\in\Z}\int\limits_{-\infty}^\infty \mathcal{I}_{j}(t,\lambda)dt d\lambda\\[0.2cm]
&=&-\frac{h}{2\pi i}\frac{1}{c}\int\limits_{\chi_-}^{\chi_+}\sum\limits_{j\in\Z}\mathcal{Z}_j(y,\lambda)d\lambda\\[0.2cm]
 &=&-\frac{h}{2\pi ic}\int\limits_{\chi_-}^{\chi_+}\sum\limits_{j\in\Z}e^{\lambda\frac{(hj_0-hj)}{c}}Z_j\Big(\frac{(hj_0-hj)}{c}+t_0\Big)d\lambda\\[0.2cm]
 &=&-\frac{h}{2\pi ic}\sum\limits_{j\in\Z}\int\limits_{\chi_-}^{\chi_+}e^{\lambda\frac{(hj_0-hj)}{c}}Z_j\Big(\frac{(hj_0-hj)}{c}+t_0\Big)d\lambda\\[0.2cm]
 &=&-\frac{h}{2\pi ic}\sum\limits_{j\in\Z}2\frac{\pi ic}{h}\delta_j^{j_0}Z_j\Big(\frac{(hj_0-hj)}{c}+t_0\Big)\\[0.2cm]
 &=&-Z_{j_0}(t_0),\end{array}\end{equation}
as desired.
\qed

\subsection{Meromorphic expansion of $G_{\lambda}$}\label{sec:grn:mero}

In this subsection we set out to explicitly isolate the
pole at $\lambda = 0$ in the meromorphic expansion
of $G_{\lambda}$. In addition, we show that
both parts of this decomposition
decay exponentially in a $\lambda$-uniform fashion.
This will allow us to shift the integration path in (\ref{4.3equivalent}) to the left of the imaginary axis.
The decomposition (\ref{zogaanwehemopdelen}) for the Green's function $\mathcal{G}$
together with the exponential decay estimates
(\ref{eq:expdecaymathcalG})
can subsequently be read off from the shifted
contour integral.\par


\begin{lemma}\label{exponentieelverval} Assume that (\asref{aannamespuls}{\text{H}}),(\asref{extraaannamespuls}{\text{H}}), (\asref{aannamesconstanten}{\text{H}}), (\asref{aannames}{\text{H}\alpha1}) and (\asref{extraaannames}{\text{H}\alpha}) are satisfied. There exist constants $K_1>0,K_2>0,\delta>0$ and $\tilde{\delta}>0$ such that
\begin{equation}\begin{array}{lcl}|\Phi^+(\xi)|&\leq &K_1 e^{-\delta |\xi|}\nrm{\Phi^+}_\infty,\\[0.2cm]
|\Phi^-(\xi)|&\leq& K_2 e^{-\tilde{\delta}|\xi|}\nrm{\Phi^-}_\infty
\end{array}\end{equation}
for all $\xi\in\R$.\end{lemma}
\vspace*{4pt}\noindent\textit{Proof.} We obtain from Lemma \ref{lemma4.3bachbos} that there are constants $\delta>0$ and $K_1>0$ for which
\begin{equation}\begin{array}{lcl}|\Psi(\xi)|&\leq &K_1 e^{-\delta|\xi|}\nrm{\Psi}_\infty+K_1\int_{-\infty}^\infty e^{-\delta|\xi-\eta|}|\Theta(\eta)|d\eta\end{array}\end{equation}
holds for each $\Psi\in \textbf{H}^1$, where $\Theta=L\Psi$. Since $L\Phi^+=0$ we conclude that
\begin{equation}\begin{array}{lcl}|\Phi^+(\xi)|&\leq &K_1 e^{-\delta|\xi|}\nrm{\Phi^+}_\infty\end{array}\end{equation}
for all $\xi$. Note that the operator $L^*$ is also asymptotically hyperbolic. Hence, there are $\tilde{\delta}>0$ and $K_2>0$ for which
\begin{equation}\begin{array}{lcl}|\Psi(\xi)|&\leq& K_2 e^{-\tilde{\delta}|\xi|}\nrm{\Psi}_\infty+K_2\int_{-\infty}^\infty e^{-\tilde{\delta}|\xi-\eta|}|\Theta(\eta)|d\eta\end{array}\end{equation}
holds for each $\Psi\in \textbf{H}^1$, where $\Theta=L^*\Psi$. Since $L^*\Phi_h^-=0$ we obtain that
\begin{equation}\begin{array}{lcl}|\Phi^-(\xi)|&\leq& K_2 e^{-\tilde{\delta}|\xi|}\nrm{\Phi^-}_\infty\end{array}\end{equation}
for all $\xi$. \qed

\begin{lemma}\label{exponentieelverval2} Assume that (\asref{aannamespuls}{\text{H}}),(\asref{extraaannamespuls}{\text{H}}), (\asref{aannamesconstanten}{\text{H}}), (\asref{aannames}{\text{H}\alpha1}) and (\asref{extraaannames}{\text{H}\alpha}) are satisfied. Then there exist constants $K_3>0$ and $\delta>0$
such that
\begin{equation}\begin{array}{lcl}|(\Phi^\pm)'(\xi)|&\leq &K_3 e^{-\delta|\xi|}\end{array}\end{equation}
for all $\xi\in\R$.
\end{lemma}
\vspace*{4pt}\noindent\textit{Proof.} Lemma \ref{exponentieelverval} implies that
\begin{equation}\begin{array}{lcl}|\Delta_h\phi^+(\xi)|&\leq&\frac{1}{h^2}K_1\sum\limits_{k>0}|\alpha_k|(e^{-\delta|\xi+hk|}+e^{-\delta|\xi-hk|}+2e^{-\delta|\xi|})\\[0.3cm]
&\leq &K_1 e^{-\delta|\xi|}(\frac{1}{h^2}\sum\limits_{k>0}|\alpha_k|(2e^{\delta  hk}+2)),\end{array}\end{equation}
where the last sum converges by (\asref{extraaannames}{\text{H}\alpha}), possibly after decreasing $\delta>0$. Using the fact that
\begin{equation}\begin{array}{lcl}(\Phi^+)'&=&\frac{1}{c}\left(\begin{array}{ll}\Delta_h \phi^++g_u(\overline{u})\phi^+-\psi^+\\ \rho \phi^+-\rho\gamma \psi^+\end{array}\right)\end{array}\end{equation}
we hence see that there exists a constant $K_3>0$ such that
\begin{equation}\begin{array}{lcl}|(\Phi^+)'(\xi)|&\leq &K_3 e^{-\delta |\xi|}.\end{array}\end{equation}
The proof for the bound on $(\Phi^-)'$ is identical.\qed

We recall the spaces
\begin{equation}\label{oude756}\begin{array}{lclcl}X&:=&X_h&=&\{\Theta\in \textbf{H}^1: \ip{\Phi^-,\Theta}=0\}\\[0.2cm]
Y&:=&Y_h&=&\{\Theta\in \textbf{L}^2: \ip{\Phi^-,\Theta}=0\},
\end{array}\end{equation}
together with the operators $L^{-1}$ in the spaces $\mathcal{B}(X,X)$ and in $\mathcal{B}(Y,X)$ that were defined in Proposition \ref{samenvattingh5}. We also recall the notation $L^{\mathrm{qinv}}\Theta$ that was introduced in Corollary \ref{vhwh} for the unique solution $\Psi$ of the equation
\begin{equation}\begin{array}{lcl}L\Psi&=&\Theta-\frac{\ip{\Phi^-,\Theta}}{\ip{\Phi^-,\Phi^+}}\Phi^+\end{array}\end{equation}
in the space $X$, which is given explicitly by
\begin{equation}\label{oude758}\begin{array}{lcl}L^{\mathrm{qinv}}\Theta&=&L^{-1}\Big[\Theta-\frac{\ip{\Phi^-,\Theta}}{\ip{\Phi^-,\Phi^+}}\Phi^+\Big].\end{array}\end{equation}
We now exploit these operators to decompose the Green's function of $\lambda+L$ into a meromorphic and an analytic part. This result is based on \cite[Lemma 2.7]{HJHSTBFHN}.
\begin{lemma}\label{decompositionat0ofgreen} Assume that (\asref{aannamespuls}{\text{H}}),(\asref{extraaannamespuls}{\text{H}}), (\asref{aannamesconstanten}{\text{H}}), (\asref{aannames}{\text{H}\alpha1}) and (\asref{extraaannames}{\text{H}\alpha}) are satisfied. There exists a constant $0<\overline{\lambda}\leq \lambda_0$ such that for all $0<|\lambda|<\overline{\lambda}$ we have the representation
\begin{equation}\begin{array}{lcl}G_{\lambda}(\xi,\xi_0)&=&E_{\lambda}(\xi,\xi_0)+\tilde{G}_{\lambda}(\xi,\xi_0)\end{array}\end{equation}
Here the meromorphic (in $\lambda$) term can be written as
\begin{equation}\begin{array}{lcl}E_{\lambda}(\xi,\xi_0)&=&-\frac{1}{\lambda \Omega}\left(\begin{array}{ll}\phi^-(\xi_0)\phi^+(\xi)&\psi^-(\xi_0)\phi^+(\xi)\\ \phi^-(\xi_0)\psi^+(\xi)&\psi^-(\xi_0)\psi^+(\xi)\end{array}\right)\end{array}\end{equation}
and the analytic (in $\lambda$) term $\tilde{G}_{\lambda}$ is given by
\begin{equation}\label{analyticterm}\begin{array}{lcl}\tilde{G}_{\lambda}(\xi,\xi_0)&=&G_{\infty;\lambda}(\xi-\xi_0) -\Big[[I+\lambda L^{-1}]^{-1}L^{\mathrm{qinv}}(L-L_{\infty}) G_{\infty;\lambda}(\cdot-\xi_0)\Big](\xi)\\[0.2cm]
&&\qquad -\frac{1}{\Omega}\ip{\Phi^-,G_{\infty;\lambda}(\cdot-\xi_0)}\Phi^+(\xi).\end{array}\end{equation}
Here we recall the notation
\begin{equation}\begin{array}{lcl}\Omega&=&\ip{\Phi^-,\Phi^+}.\end{array}\end{equation}
\end{lemma}
\noindent\textit{Proof.} Pick $\lambda\in\C$ with $0<|\lambda|<\lambda_0$. By the proof of Proposition \ref{lemmaspectrumklein} we see that
\begin{equation}\begin{array}{lcl}(L+\lambda)^{-1}\Theta&=&\lambda^{-1}\frac{\ip{\Phi^-,\Theta}}{\Omega}\Phi^++L^{\mathrm{qinv}}\Theta-[I+\lambda L^{-1}]^{-1}\lambda L^{-1}L^{\mathrm{qinv}}\Theta.\end{array}\end{equation}

We now compute
\begin{equation}\begin{array}{lcl}\ip{\Phi^-,(L-L_{\infty})G_{\infty;\lambda}(\cdot-\xi_0)}&=&\ip{\Phi^-,-L_{\infty} G_{\infty;\lambda}(\cdot-\xi_0)}\\[0.2cm]
&=&-\Phi^-(\xi_0)+\lambda\ip{\Phi^-,G_{\infty;\lambda}(\cdot-\xi_0)}.\end{array}\end{equation}
In particular, writing
\begin{equation}
\begin{array}{lcl}
\hat{L}&=&L-L_{\infty},
\end{array}
\end{equation}
we obtain
\begin{equation}\begin{array}{lcl}(L+\lambda)^{-1}\hat{L}G_{\infty}(\cdot-\xi_0)&=&\frac{1}{\lambda \Omega}\left(\begin{array}{ll}\phi^-(\xi_0)\phi^+&\psi^-(\xi_0)\phi^+\\ \phi^-(\xi_0)\psi^+&\psi^-(\xi_0)\psi^+\end{array}\right)\\[0.4cm]
&&\qquad +\frac{\ip{\Phi^-,G_{\infty;\lambda}(\cdot-\xi_0)}}{\Omega}\Phi^+ +L^{\mathrm{qinv}}\hat{L}G_{\infty;\lambda}(\cdot-\xi_0)\\[0.2cm]
&&\qquad -[I+\lambda L^{-1}]^{-1}\lambda L^{-1}L^{\mathrm{qinv}}\hat{L}G_{\infty;\lambda}(\cdot-\xi_0).\end{array}\end{equation}
We may hence write
\begin{equation}\begin{array}{lcl}G_{\lambda}(\xi,\xi_0)&=&E_{\lambda}(\xi,\xi_0)+\tilde{G}_{\lambda}(\xi,\xi_0)\end{array}\end{equation}
with
\begin{equation}\begin{array}{lcl}E_{\lambda}(\xi,\xi_0)&=&-\frac{1}{\lambda \Omega}\left(\begin{array}{ll}\phi^-(\xi_0)\phi^+(\xi)&\psi^-(\xi_0)\phi^+(\xi)\\ \phi^-(\xi_0)\psi^+(\xi)&\psi^-(\xi_0)\psi^+(\xi)\end{array}\right)\end{array}\end{equation}
and
\begin{equation}\begin{array}{lcl}\tilde{G}_{\lambda}(\cdot,\xi_0)&=&G_{\infty;\lambda}(\cdot-\xi_0)-L^{\mathrm{qinv}}\hat{L}G_{\infty\lambda}(\cdot-\xi_0)\\[0.2cm]
&&\qquad +[I+\lambda L^{-1}]^{-1}\lambda L^{-1}L^{\mathrm{qinv}}\hat{L}G_{\infty;\lambda}(\cdot-\xi_0)\\[0.2cm]
&&\qquad -\frac{1}{\Omega}\ip{\Phi^-,G_{\infty;\lambda}(\cdot-\xi_0)}\Phi^+\\[0.2cm]

&=&G_{\infty;\lambda}(\cdot-\xi_0) -[I+\lambda L^{-1}]^{-1}L^{\mathrm{qinv}}\hat{L} G_{\infty;\lambda}(\cdot-\xi_0)\\[0.2cm]
&&\qquad -\frac{1}{\Omega}\ip{\Phi^-,G_{\infty}(\cdot-\xi_0)}\Phi^+.\end{array}\end{equation}
Clearly $E_{\lambda}$ is meromorphic in $\lambda$, while $\tilde{G}_{\lambda}$ is analytic in $\lambda$ in the region $|\lambda|<\lambda_0$.\qed

We fix $\chi>\lambda_{\mathrm{unif}}$, where $\lambda_{\mathrm{unif}}$ was defined in Lemma \ref{existenceofgreen}, and set
\begin{equation}
\begin{array}{lcl}
R&=&\{\lambda\in\C:-\frac{\tilde{\lambda}}{2}\leq \Re \,\lambda \leq \chi\hbox{ and }|\Im \,\lambda|\leq\frac{\pi |c|}{h}\}.
\end{array}
\end{equation}
We now set out to obtain an estimate on the function $\tilde{G}_{\lambda}$ from Lemma \ref{decompositionat0ofgreen} by exploiting the asymptotic hyperbolicity of $L$. We treat each of the terms in (\ref{analyticterm}) separately in the results below.
\begin{lemma}\label{decompositiebewijs0} Assume that (\asref{aannamespuls}{\text{H}}),(\asref{extraaannamespuls}{\text{H}}), (\asref{aannamesconstanten}{\text{H}}), (\asref{aannames}{\text{H}\alpha1}) and (\asref{extraaannames}{\text{H}\alpha}) are satisfied. There exist constants $K_4>0$ and $\tilde{\chi}>0$ such that for all $\lambda\in R$
\begin{equation}\begin{array}{lcl}\big|\ip{\Phi^-,(L-L_{\infty})G_{\infty;\lambda}(\cdot-\xi_0)}\big|&\leq &K_4 e^{-\tilde{\chi}|\xi_0|}.\end{array}\end{equation}\end{lemma}
\vspace*{4pt}\noindent\textit{Proof.} We reuse the notation $\hat{L}=L-L_{\infty}$ from the previous proof. Lemma \ref{eigenschappenghinfty} implies that we can pick constants $\beta_*>0$ and $K_*>0$ in such a way that
\begin{equation}\begin{array}{lcl}|G_{\infty;\lambda}(\xi-\xi_0)|&\leq &K_* e^{-\beta_*|\xi-\xi_0|}\end{array}\end{equation}
for all values of $\xi,\xi_0$. Recall the constants $K_2,\tilde{\delta}$ from Lemma \ref{exponentieelverval} and set $K_3=K_2\nrm{\Phi^-}_\infty$. Then we obtain
\begin{equation}\begin{array}{lcl}&&|\ip{\Phi^-,\hat{L}G_{\infty;\lambda}(\cdot-\xi_0)}|\\[0.2cm]&\leq &\int_{-\infty}^\infty K_3 e^{-\tilde{\delta}|\xi|}g_* K_* e^{-\beta_*|\xi-\xi_0|}\ d\xi\\[0.2cm]
&=&K_3g_*K_*\Big(\frac{1}{\tilde{\delta}+\beta_*}(e^{-\tilde{\delta}|\xi_0|}+e^{-\beta_*|\xi_0|}) +\frac{1}{\beta_*-\tilde{\delta}}(e^{-\tilde{\delta}|\xi_0|}-e^{-
\beta_*|\xi_0|})\Big)\\[0.2cm]
&\leq &K_3g_* K_*\Big(\frac{1}{\tilde{\delta}+\beta_*}2e^{-\min\{\tilde{\delta},\beta_*\}|\xi_0|} +\frac{1}{|\beta_*-\tilde{\delta}|}2e^{-\min\{\tilde{\delta},\beta_*\}|\xi_0|}\Big)\\[0.2cm]
&=&K_4 e^{-\tilde{\chi}|\xi_0|}\end{array}\end{equation}
for some $K_4>0$ and $\tilde{\chi}>0$.\qed

\begin{lemma}\label{decompositiebewijs} Assume that (\asref{aannamespuls}{\text{H}}),(\asref{extraaannamespuls}{\text{H}}), (\asref{aannamesconstanten}{\text{H}}), (\asref{aannames}{\text{H}\alpha1}) and (\asref{extraaannames}{\text{H}\alpha}) are satisfied. There exist constants $K_{10}>0$ and $\tilde{\gamma}>0$ such that for all $\lambda\in R$
\begin{equation}\begin{array}{lcl}\Big|\Big[L^{\mathrm{qinv}}(L-L_{\infty})G_{\infty;\lambda}(\cdot-\xi_0)\Big](\xi)\Big|&\leq &K_{10} e^{-\tilde{\gamma}|\xi|}e^{-\tilde{\gamma}|\xi_0|}\\[0.2cm]
&\leq &K_{10} e^{-\tilde{\gamma}|\xi-\xi_0|}.\end{array}\end{equation}
\end{lemma}
\vspace*{4pt}\noindent\textit{Proof.} We reuse the notation $\hat{L}=L-L_{\infty}$ from the previous proof. Recall the constants $K_1,\delta$ from Lemma \ref{exponentieelverval}. Writing
\begin{equation}\begin{array}{lcl}H_{\xi_0}(\xi)&=&\Big[L^{\mathrm{qinv}}\hat{L}G_{\infty;\lambda}(\cdot-\xi_0)\Big](\xi),\end{array}\end{equation}
we may use Lemma \ref{lemma4.3bachbos} to estimate
\begin{equation}\label{tildeafschatting0}
\begin{array}{lcl}
|H_{\xi_0}(\xi)|&\leq & K_1 e^{-\delta|\xi|}\nrm{H_{\xi_0}}_\infty +K_1\int_{-\infty}^\infty e^{-\delta|\xi-\eta|}|LH_{\xi_0}(\eta)|d\eta.
\end{array}
\end{equation}

Recalling (\ref{oude756})-(\ref{oude758}), we obtain
\begin{equation}\label{tildeafschatting1}\begin{array}{lcl}\nrm{H_{\xi_0}}_{\infty}&\leq &\nrm{H_{\xi_0}}_{\textbf{H}^1}\\[0.2cm]
&\leq &C_{\mathrm{unif}}||\hat{L}G_{\infty;\lambda}(\cdot-\xi_0)-\frac{\ip{\Phi^-,\hat{L}G_{\infty;\lambda}(\cdot-\xi_0)}}{\Omega}\Phi^+||_{\textbf{L}^2}\\[0.2cm]
&\leq &C_{\mathrm{unif}}\Big(1+\frac{\nrm{\Phi^-}_{\textbf{L}^2}}{\Omega}\nrm{\Phi^+}_{\textbf{L}^2}\Big)  \nrm{\hat{L}G_{\infty;\lambda}(\cdot-\xi_0)}_{\textbf{L}^2}\\[0.2cm]
&\leq & K_5\nrm{\hat{L}G_{\infty;\lambda}(\cdot-\xi_0)}_{\textbf{L}^2}\end{array}\end{equation}
for some constant $K_5>0$.\par

Using Lemma \ref{exponentieelverval} we see that there exists a constant $K_6>0$ for which
\begin{equation}
\begin{array}{lcl}
|\overline{u}(\xi)|&=&|\int_{\infty}^\xi \overline{u}'(\xi')\ d\xi'|\\[0.2cm]
&\leq &\int_{\xi}^\infty K_1\nrm{(\overline{u}',\overline{w}')}_{\textbf{L}^2} e^{-\delta |\xi'|}\ d\xi'\\[0.2cm]
&=&K_6 e^{-\delta |\xi|}
\end{array}
\end{equation}
holds for all $\xi\in\R$. Recall that
\begin{equation}\begin{array}{lcl}\hat{L}&=&\left(\begin{array}{ll}-g_u(\overline{u})+r_0& 0\\ 0&0\end{array}\right).\end{array}\end{equation}
Observe that $-g_u(0)+r_0=0$. Then we obtain that
\begin{equation}\label{probleemvgl}
\begin{array}{lcl}
|-g_u(\overline{u}(\xi))+r_0|&\leq & K_7 e^{-\delta|\xi|}
\end{array}
\end{equation}
for all $\xi\in\R$ and for some constant $K_7>0$. Lemma \ref{eigenschappenghinfty} implies that
\begin{equation}
\begin{array}{lcl}
|G_{\infty;\lambda}(\xi-\xi_0)|&\leq &K_*e^{-\beta_*|\xi-\xi_0|}
\end{array}
\end{equation}
for all $\xi\in\R$. Therefore, we must have
\begin{equation}\label{tildeafschatting2}
\begin{array}{lcl}
\nrm{\hat{L}G_{\infty;\lambda}(\cdot-\xi_0)}_{\textbf{L}^2}^2&\leq &\int_\R K_7^2K_*^2 e^{-2\delta|\xi|}e^{-2\beta_*|\xi-\xi_0|}\ d\xi\\[0.2cm]
&\leq & K_8 e^{-2\tilde{\gamma}|\xi_0|}
\end{array}
\end{equation}
for some constants $K_8>0$, $\tilde{\gamma}>0$ with $\tilde{\gamma}\leq \beta_*$, $\tilde{\gamma}\leq \frac{1}{2}\delta$ and $\tilde{\gamma}\leq\frac{1}{2}\tilde{\chi}$. In particular, we obtain the estimate
\begin{equation}\label{tildeafschattingfinal}
\begin{array}{lcl}
\nrm{H_{\xi_0}}_\infty &\leq & K_5\sqrt{K_8}e^{-\tilde{\gamma}|\xi_0|}.
\end{array}
\end{equation}
In a similar fashion, using Lemma \ref{decompositiebewijs0}, we see that
\begin{equation}
\begin{array}{lcl}
|LH_{\xi_0}(\xi)|&\leq &\Big|\Big[\hat{L}G_{\infty;\lambda}(\cdot-\xi_0)\Big](\xi)-\frac{\ip{\Phi^-,\hat{L}G_{\infty;\lambda}(\cdot-\xi_0)}}{\Omega}\Phi^+(\xi)\Big|\\[0.2cm]
&\leq & K_7 K_* e^{-\delta|\xi|}e^{-\beta_*|\xi-\xi_0|}+\frac{1}{\Omega}K_4 e^{-\tilde{\chi}|\xi_0|} K_1 e^{-\delta|\xi|}\\[0.2cm]
&\leq & K_{9}\Big[e^{-2\tilde{\gamma}|\xi|} e^{-\tilde{\gamma}|\xi-\xi_0|}+e^{-\tilde{\gamma}|\xi_0|}e^{-\tilde{\gamma}|\xi|}\Big]
\end{array}
\end{equation}
for all $\xi\in\R$ and some constant $K_{9}>0$. Combining (\ref{tildeafschatting0}) with (\ref{tildeafschatting1}) and (\ref{tildeafschatting2}), we hence obtain
\begin{equation}
\begin{array}{lcl}
|H_{\xi_0}(\xi)|&\leq & K_1 e^{-\delta|\xi|}\nrm{H_{\xi_0}}_\infty +K_1\int_{-\infty}^\infty e^{-\delta|\xi-\eta|}|LH_{\xi_0}(\eta)|d\eta\\[0.2cm]
&\leq &K_1 e^{-\delta|\xi|} K_5 \sqrt{K_8}e^{-\tilde{\gamma}|\xi_0|}\\[0.2cm]
&&\qquad +K_1 \int_{-\infty}^\infty e^{-\delta|\xi-\eta|}K_{9}\Big[e^{-2\tilde{\gamma}|\eta|} e^{-\tilde{\gamma}|\eta-\xi_0|}+e^{-\tilde{\gamma}|\xi_0|}e^{-\tilde{\gamma}|\eta|}\Big]d\eta\\[0.2cm]
&\leq &K_1 e^{-\delta|\xi|} K_5 \sqrt{K_8}e^{-\tilde{\gamma}|\xi_0|}+K_1 \int_{-\infty}^\infty e^{-\delta|\xi-\eta|}2K_{9} e^{-\tilde{\gamma}|\eta|}e^{-\tilde{\gamma}|\xi_0|}d\eta\\[0.2cm]
&\leq & K_{10} e^{-\tilde{\gamma}|\xi|}e^{-\tilde{\gamma}|\xi_0|}\\[0.2cm]
&\leq & K_{10} e^{-\tilde{\gamma}|\xi-\xi_0|}
\end{array}
\end{equation}
for some constant $K_{10}>0$.
\qed

\begin{remark} In the proof of Lemma \ref{decompositiebewijs}, in particular in (\ref{probleemvgl}), we explicitly used that $\overline{U}$ is a pulse solution, instead of a traveling front solution. If one would want to transfer these results to a more general system where the waves have different limits at $\xi=\pm\infty$, then Lemma \ref{decompositiebewijs} would only hold for $\xi_0\geq 0$. However, the definition (\ref{definitievanGhlambda})
remains valid upon using the reference system
at $\xi = - \infty$ instead of $\xi = +\infty$.
This new formulation allows
the desired estimates for $\xi_0 \le 0$ to
be recovered.
\end{remark}

\begin{lemma}\label{afschattingtildeG} Assume that (\asref{aannamespuls}{\text{H}}),(\asref{extraaannamespuls}{\text{H}}), (\asref{aannamesconstanten}{\text{H}}), (\asref{aannames}{\text{H}\alpha1}) and (\asref{extraaannames}{\text{H}\alpha}) are satisfied. There exist constants $K_{13}>0$ and $\omega>0$ such that the function $\tilde{G}_{\lambda}$ from Lemma \ref{decompositionat0ofgreen} satisfies the bound
\begin{equation}\begin{array}{lcl}|\tilde{G}_{\lambda}(\xi,\xi_0)|&\leq &K_{13} e^{-\omega|\xi-\xi_0|}\end{array}\end{equation}
for all $\xi,\xi_0$ and all $0<|\lambda|<\overline{\lambda}$.\end{lemma}
\vspace*{4pt}\noindent\textit{Proof.} As before, we write
\begin{equation}\begin{array}{lcl}H_{\xi_0}(\xi)&=&L^{\mathrm{qinv}}\hat{L}G_{\infty;\lambda}(\cdot-\xi_0)(\xi).\end{array}\end{equation}
Using Lemma \ref{lemma4.3bachbos} Lemma \ref{decompositiebewijs} and (\ref{tildeafschattingfinal}) and recalling (\ref{oude756})-(\ref{oude758}), we obtain the estimate
\begin{equation}\label{inductiestap1}\begin{array}{lcl} |L^{-1}H_{\xi_0}(\xi)|&\leq & K_1 e^{-\alpha|\xi|}\nrm{L^{-1}H_{\xi_0}}_\infty+K_1\int_{-\infty}^\infty e^{-\alpha|\xi-\eta|}|H_{\xi_0}(\eta)|d\eta\\[0.2cm]
&\leq &K_1 e^{-\alpha|\xi|}C_{\mathrm{unif}}\nrm{H_{\xi_0}}_{\textbf{L}^2}+K_1\int_{-\infty}^\infty e^{-\alpha|\xi-\eta|}K_{10} e^{-\tilde{\gamma}|\eta-\xi_0|}d\eta\\[0.2cm]
&\leq &K_1 e^{-\alpha|\xi|}C_{\mathrm{unif}}K_5\sqrt{K_8}e^{-\tilde{\gamma}|\xi_0|}+K_1\int_{-\infty}^\infty e^{-\alpha|\xi-\eta|}K_{10} e^{-\tilde{\gamma}|\eta-\xi_0|}d\eta\\[0.2cm]
&\leq & K_{10}K_{11}e^{-\tilde{\gamma}|\xi-\xi_0|}\end{array}\end{equation}
for some constants $K_{11}>0$ and $2\tilde{\gamma}\leq \alpha$. Using Proposition \ref{samenvattingh5} and (\ref{tildeafschattingfinal}) we obtain that
\begin{equation}\begin{array}{lcl} \nrm{(L^{-1})^nH_{\xi_0}}_{\textbf{H}^1}&\leq &K_5\sqrt{K_8}(C_{\mathrm{unif}})^n e^{-\tilde{\gamma}|\xi_0|}\end{array}\end{equation}
for all $n\in\Z_{>0}$. Continuing in this fashion, we see that
\begin{equation}\begin{array}{lcl} |(L^{-1})^nH_{\xi_0}(\xi)|&\leq &K_{10}K_{11}^n e^{-\tilde{\gamma}|\xi-\xi_0|}\end{array}\end{equation}
for all $n\in\Z_{>0}$. If we set
\begin{equation}
\begin{array}{lcl}
\overline{\lambda}&=&\min\{\frac{\tilde{\lambda}}{2},\lambda_0,\chi,\frac{1}{C_{\mathrm{unif}}K_5\sqrt{K_8}},\frac{1}{K_{11}}\},
\end{array}
\end{equation}
then for each $n\in\Z_{>0}$ and each $0<|\lambda|<\overline{\lambda}$ we have
\begin{equation}\begin{array}{lcl} \nrm{(-\lambda)^n(L^{-1})^nH_{\xi_0}}_{\textbf{H}^1}&\leq &\frac{1}{2}.\end{array}\end{equation}
In particular, it follows that
\begin{equation}
\begin{array}{lcl}
\sum\limits_{n=0}^N (-\lambda)^n(L^{-1})^nH_{\xi_0}&\rightarrow &[I+\lambda L^{-1}]^{-1}H_{\xi_0}
\end{array}
\end{equation}
in $\textbf{H}^1$ as $N\rightarrow \infty$. Since $\textbf{H}^1$-convergence implies point-wise convergence, we conclude that
\begin{equation}
\begin{array}{lcl}\big|[I+\lambda L^{-1}]^{-1}H_{\xi_0}(\xi)\big|&=&
\Big|\sum\limits_{n=0}^\infty (-\lambda)^n(L^{-1})^nH_{\xi_0}(\xi)\Big|\\[0.2cm]
&\leq& \sum\limits_{n=0}^\infty \overline{\lambda}^n K_{11}K_{12}^ne^{-\tilde{\gamma}|\xi-\xi_0|}\\[0.2cm]
&\leq & \frac{K_{11}}{1-\overline{\lambda} K_{11}}e^{-\tilde{\gamma}|\xi-\xi_0|}\\[0.2cm]
&:=&K_{12}e^{-\tilde{\gamma}|\xi-\xi_0|}
\end{array}
\end{equation}
for all $\xi\in\R$ and for some constant $K_{12}>0$.\par

Combining this estimate with Lemma \ref{exponentieelverval} and Lemma \ref{decompositiebewijs} yields the desired bound
\begin{equation}\begin{array}{lcl}|\tilde{G}_{\lambda}(\xi,\xi_0)|&=&\big|G_{\infty;\lambda}(\xi-\xi_0) -\Big[[I+\lambda L^{-1}]^{-1}L^{\mathrm{qinv}}\hat{L} G_{\infty;\lambda}(\cdot-\xi_0)\Big](\xi)\\[0.2cm]
&&\qquad -\frac{1}{\Omega}\ip{\Phi^-,G_{\infty}(\cdot-\xi_0)}\Phi^+(\xi)\big|\\[0.2cm]
&\leq& K_* e^{-\beta_*|\xi-\xi_0|}+K_{12}e^{-\tilde{\gamma}|\xi-\xi_0|}+K_4 \frac{1}{\Omega}e^{-\tilde{\chi}|\xi_0|}K_1 e^{-\delta|\xi|}\nrm{\Phi^+}_{\infty}\\[0.2cm]
&\leq &K_{13}e^{-\omega|\xi-\xi_0|}
\end{array}\end{equation}
for some constants $K_{13}>0$ and $\omega>0$.\qed


We write
\begin{equation}
\begin{array}{lcl}
S&=&\{-\overline{\lambda}+i\omega:\omega\in[-\frac{\pi |c|}{h},\frac{\pi |c|}{h}]\},
\end{array}
\end{equation}
where $\overline{\lambda}$ is defined in the proof of Lemma \ref{afschattingtildeG}.
\begin{lemma}\label{afschattingghzondertilde} Assume that (\asref{aannamespuls}{\text{H}}),(\asref{extraaannamespuls}{\text{H}}), (\asref{aannamesconstanten}{\text{H}}), (\asref{aannames}{\text{H}\alpha1}) and (\asref{extraaannames}{\text{H}\alpha}) are satisfied. Then there exist constants $K>0$ and $\tilde{\beta}>0$ such that for all $\lambda\in S$ we have the bound \begin{equation}\begin{array}{lcl}|G_{\lambda}(\xi,\xi_0)|
&\leq &K e^{-\tilde{\beta}|\xi-\xi_0|}\end{array}\end{equation}
for all $\xi,\xi_0$.
\end{lemma}
\vspace*{4pt}\noindent\textit{Proof.} Fix $\lambda_0\in S$. For $\lambda\in S$ sufficiently close to $\lambda$ we have
\begin{equation}
\begin{array}{lcl}
\Big[L+\lambda\Big]^{-1}&=&\Big[L+\lambda_0+\lambda-\lambda_0\Big]^{-1}\\[0.2cm]
&=&\Big[\Big(L+\lambda_0\Big)\Big(I+(L+\lambda_0)^{-1}(\lambda-\lambda_0)\Big)\Big]^{-1}\\[0.2cm]
&=&\Big[I+(L+\lambda_0)^{-1}(\lambda-\lambda_0)\Big]^{-1}\Big[L+\lambda_0\Big]^{-1}.
\end{array}
\end{equation}
In particular, upon writing
\begin{equation}
\begin{array}{lcl}
H_{\xi_0}(\xi)&=&\Big[[L+\lambda_0]^{-1}\hat{L}G_{\lambda;\infty}(\cdot-\xi_0)\Big](\xi),
\end{array}
\end{equation}
we see that
\begin{equation}
\begin{array}{lcl}
G_{\lambda}(\xi,\xi_0)-G_{\infty;\lambda}(\xi-\xi_0)&=&\Big[[I+(L+\lambda_0)^{-1}(\lambda-\lambda_0)]^{-1}H_{\xi_0}\Big](\xi).
\end{array}
\end{equation}
Using Lemma \ref{lemma4.3bachbos} we can pick constants $k_{\lambda_0}>0$ and $\alpha_{\lambda_0}>0$ in such a way that
\begin{equation}
\begin{array}{lcl}
|H_{\xi_0}(\xi)|&\leq & k_{\lambda_0}e^{-\alpha_{\lambda_0}|\xi|}\nrm{H_{\xi_0}}_\infty+k_{\lambda_0}\int_{-\infty}^\infty e^{-\alpha_{\lambda_0}|\xi-\eta|}|(L+\lambda_0)H_{\xi_0}(\eta)|d\eta.
\end{array}
\end{equation}
Recall the constant $C_S$ appearing in Proposition \ref{equivalenttheorem4version2}. This allows us to estimate
\begin{equation}
\begin{array}{lcl}
\nrm{H_{\xi_0}}_\infty&\leq &\nrm{H_{\xi_0}}_{\textbf{H}^1}\\[0.2cm]
&\leq & C_S \nrm{\hat{L}G_{\lambda_0}(\xi,\xi_0)}_{\textbf{L}^2}\\[0.2cm]
&\leq & C_S \sqrt{K_8}e^{-\tilde{\gamma}|\xi_0|}.
\end{array}
\end{equation}
This yields the bound
\begin{equation}
\begin{array}{lcl}
|H_{\xi_0}(\xi)|&\leq & k_{\lambda_0}e^{-\alpha_{\lambda_0}|\xi|}C_S \sqrt{K_8}e^{-\tilde{\gamma}|\xi_0|}+k_{\lambda_0}\int_{-\infty}^\infty e^{-\alpha_{\lambda_0}|\xi-\eta|}|\hat{L}G_{\lambda;\infty}(\eta,\xi_0)|d\eta\\[0.2cm]
&\leq & k_{\lambda_0}e^{-\alpha_{\lambda_0}|\xi|}C_S \sqrt{K_8}e^{-\tilde{\gamma}|\xi_0|}\\[0.2cm]
&&\qquad +k_{\lambda_0}\int_{-\infty}^\infty e^{-\alpha_{\lambda_0}|\xi-\eta|}K_7K_*e^{-\delta|\eta|}e^{-2\beta_*|\eta-\xi_0|}d\eta\\[0.2cm]
&\leq & k_{\lambda_0;2}e^{-\alpha_{\lambda_0;2}|\xi-\xi_0|}
\end{array}
\end{equation}
for some constants $k_{\lambda_0;2},\alpha_{\lambda_0;2}$, which may depend on $\lambda_0$, but not on $\lambda$. Arguing as in (\ref{inductiestap1}), we obtain
\begin{equation}
\begin{array}{lcl}
\big|[L+\lambda_0]^{-1}H_{\xi_0}(\xi)\big|&\leq & k_{\lambda_0}e^{-\alpha_{\lambda_0}|\xi|}\nrm{[L+\lambda_0]^{-1}H_{\xi_0}}_\infty\\[0.2cm]
&&\qquad +k_{\lambda_0}\int_{-\infty}^\infty e^{-\alpha_{\lambda_0}|\xi-\eta|}|H_{\xi_0}(\eta)|d\eta\\[0.2cm]
&\leq & k_{\lambda_0;2}k_{\lambda_0;3}e^{-\alpha_{\lambda_0;2}|\xi-\xi_0|}
\end{array}
\end{equation}
for some constant $k_{\lambda_0;3}>0$, which may depend on $\lambda_0$, but not on $\lambda$. Following the same steps as the proof of Lemma \ref{afschattingtildeG} and setting
\begin{equation}
\begin{array}{lcl}
\epsilon_{\lambda_0}&=&\min\{\frac{1}{k_{\lambda_0}C_S \sqrt{K_8}},\frac{1}{k_{\lambda_0;3}}\},
\end{array}
\end{equation}
we conclude that
\begin{equation}
\begin{array}{lcl}
|G_{\lambda}(\xi,\xi_0)-G_{\infty;\lambda}(\xi-\xi_0)|&=&
\big|\Big[[I+[L+\lambda_0]^{-1}(\lambda-\lambda_0)]^{-1}H_{\xi_0}\Big](\xi)\big|\\[0.2cm]
&\leq &k_{\lambda_0;4}e^{-\alpha_{\lambda_0;2}|\xi-\xi_0|}
\end{array}
\end{equation}
holds for each $\lambda\in S$ with $|\lambda-\lambda_0|<\epsilon_{\lambda_0}$, for some constant $k_{\lambda_0;4}>0$, which may depend on $\lambda_0$. In particular, we obtain that
\begin{equation}
\begin{array}{lcl}
|G_{\lambda}(\xi,\xi_0)|&\leq & k_{\lambda_0;4}e^{-\alpha_{\lambda_0;2}|\xi-\xi_0|}+K_*e^{-\beta_*|\xi-\xi_0|}\\[0.2cm]
&\leq &k_{\lambda_0;5}e^{-\alpha_{\lambda_0;2}|\xi-\xi_0|}
\end{array}
\end{equation}
holds for each $\lambda\in S$ with $|\lambda-\lambda_0|<\epsilon_{\lambda_0}$, for some constant $k_{\lambda_0;5}>0$, which may depend on $\lambda_0$.\par

Since $S$ is compact we can find $\lambda_1,...,\lambda_n\in S$ in such a way that
\begin{equation}
\begin{array}{lcl}
S&\subset &\bigcup\limits_{i=1}^n \{\lambda\in \C:|\lambda-\lambda_i|<\epsilon_{\lambda_i}\}.
\end{array}
\end{equation}
Setting
\begin{equation}
\begin{array}{lcl}
K&=&\max\{k_{\lambda_i;5}:i\in\{1,...n\}\},\\[0.2cm]
\tilde{\beta}&=&\min\{\alpha_{\lambda_i;2}:i\in\{1,...,n\}\},
\end{array}
\end{equation}
we conclude that
\begin{equation}
\begin{array}{lcl}
|G_{\lambda}(\xi,\xi_0)|&\leq & Ke^{-\tilde{\beta}|\xi-\xi_0|}
\end{array}
\end{equation}
holds for all $\lambda\in S$ and all $\xi,\xi_0\in\R$.
\qed

\subsection{Decomposition into stable and center modes}

In this final subsection we establish
Proposition \ref{cor2.8equivalentsterker}.
In particular,
the decomposition
(\ref{zogaanwehemopdelen}) and the exponential
bounds (\ref{eq:expdecaymathcalG})
for the Green's function $\mathcal{G}$
can be found by using the
splitting of $G_{\lambda}$ obtained in
\S\ref{sec:grn:mero}.
This is performed in Lemma
\ref{cor2.8equivalent},
which is based on \cite[Corollary 2.8]{HJHSTBFHN}.

We subsequently carefully study the terms appearing
in (\ref{zogaanwehemopdelen}) and show
that they can be interpreted as a spectral decomposition
that splits the flow associated to the linear
system (\ref{definitieGhuh})
into two invariant subspaces.
The stable component decays exponentially in a uniform fashion, while the center component can be described explicitly.


\begin{lemma}\label{cor2.8equivalent} Assume that (\asref{aannamespuls}{\text{H}}),(\asref{extraaannamespuls}{\text{H}}), (\asref{aannamesconstanten}{\text{H}}), (\asref{aannames}{\text{H}\alpha1}) and (\asref{extraaannames}{\text{H}\alpha}) are satisfied. For any pair $t\geq t_0$ and any $j,j_0\in\Z$, we have the representation
\begin{equation}\begin{array}{lcl}\mathcal{G}_j^{j_0}(t,t_0)&=&\mathcal{E}_j^{j_0}(t,t_0)+\tilde{\mathcal{G}}_j^{j_0}(t,t_0)\end{array}\end{equation}
in which
\begin{equation}\begin{array}{lcl}\mathcal{E}_j^{j_0}(t,t_0)&=&\frac{h}{\Omega}\left(\begin{array}{ll}\phi^-(hj_0+ct_0)\phi^+(hj+ct)& \psi^-(hj_0+ct_0)\phi^+(hj+ct)\\ \phi^-(hj_0+ct_0)\psi^+(hj+ct)& \psi^-(hj_0+ct_0)\psi^+(hj+ct)\end{array}\right),\end{array}\end{equation}
while $\tilde{\mathcal{G}}$ satisfies the bound
\begin{equation}\label{preciezeafschattingGtilde}\begin{array}{lcl}|\tilde{\mathcal{G}}_j^{j_0}(t,t_0)| &\leq & K e^{-\tilde{\beta}(t-t_0)}e^{-\tilde{\beta}|hj+ct-hj_0-ct_0|}\end{array}\end{equation}
for some $K>0$ and $\tilde{\beta}>0$. \end{lemma}
\noindent\textit{Proof.} Recall the representation of $\mathcal{G}_j^{j_0}$ from Proposition \ref{thrGreenfunction}. Note that $G_{\lambda}(\xi,\xi_0)$ is meromorphic for $\lambda$ in the strip $\{\lambda\in\C: \Re \,\lambda\geq-\lambda_3,|\Im \,\lambda|\leq\frac{c\pi}{h}\}$ with a simple pole at $\lambda=0$ by Lemma \ref{decompositionat0ofgreen}, Lemma \ref{existenceofgreen} and Theorem \ref{totalespectrum}. Lemma \ref{decompositionat0ofgreen} also implies that the residue of $G_{\lambda}(\xi,\xi_0)$ in $\lambda=0$ is given by
\begin{equation}\begin{array}{lcl}\text{Res}(G_{\lambda}(\xi,\xi_0),0)&=&-\frac{1}{\Omega}\left(\begin{array}{ll}\phi^-(\xi_0)\phi^+(\xi)&\psi^-(\xi_0)\phi^+(\xi)\\ \phi^-(\xi_0)\psi^+(\xi)&\psi^-(\xi_0)\psi^+(\xi)\end{array}\right).\end{array}\end{equation}
We write
\begin{equation}\begin{array}{lcl}H(\cdot,\xi_0)&=&e^{2\pi i \frac{1}{h}k\xi_0}(L+\lambda+2\pi i k\frac{c}{h})e_{-2\pi i \frac{1}{h}k}G_{\lambda}(\cdot,\xi_0).\end{array}\end{equation}
In a similar fashion as in the proof of Lemma \ref{lemmaperiodiek} we see that for $k\in\Z$ we have
\begin{equation}
\begin{array}{lcl}
(L+\lambda+2\pi i k\frac{c}{h})e_{-2\pi i \frac{1}{h}k}&=&e_{-2\pi i \frac{1}{h}k}(L+\lambda).
\end{array}
\end{equation}
Therefore, it follows that
\begin{equation}\begin{array}{lcl}H(\cdot,\xi_0)&=&e^{2\pi i\frac{1}{h}k\xi_0}(L+\lambda+2\pi i k\frac{c}{h})e_{-2\pi i \frac{1}{h}k}G_{\lambda}(\cdot,\xi_0)\\[0.2cm]
&=&e^{2\pi i\frac{1}{h}k\xi_0}e_{-2\pi i \frac{1}{h}k}(L+\lambda)G_{\lambda}(\cdot,\xi_0).\end{array}\end{equation}
For any $f\in \textbf{H}^1$ we may hence compute
\begin{equation}\begin{array}{lcl}\int H(\xi,\xi_0)f(\xi_0)\ d\xi_0&=&\int e^{2\pi i\frac{1}{h}k\xi_0}e^{-2\pi i \frac{1}{h}k\xi}(L +\lambda)G_{\lambda}(\cdot,\xi_0)(\xi)f(\xi_0)\ d\xi_0\\[0.2cm]
&=&e^{-2\pi i \frac{1}{h}k\xi}[e^{2\pi i\frac{1}{h}k\xi}f(\xi)]\\[0.2cm]
&=&f(\xi).\end{array}\end{equation}
Therefore, by the invertibility of $L+\lambda+2\pi ik\frac{c}{h}$, we must have
\begin{equation}\begin{array}{lcl}G_{\lambda+2\pi ik\frac{c}{h}}(\xi,\xi_0)&=&e^{2\pi ik\frac{1}{h}(\xi_0-\xi)}G_{\lambda}(\xi,\xi_0).\end{array}\end{equation}
Now recall the constants $\chi,\chi_+,\chi_-$ from (the proof of) Proposition \ref{thrGreenfunction} and define
\begin{equation}\begin{array}{lcl}\overline{\lambda}^-&=&-\frac{\overline{\lambda}}{2}-i\frac{\pi c}{h}\\[0.2cm]
\overline{\lambda}^+ &=&-\frac{\overline{\lambda}}{2}+i\frac{\pi c}{h}.\end{array}\end{equation}
Writing $x=hj+ct,y=hj_0+ct_0$, we see that
\begin{equation}\begin{array}{lcl}\int_{\overline{\lambda}^-}^{\chi_-}e^{\lambda(t-t_0)}G_{\lambda}(x,y)\ d\lambda &=&\int_{\overline{\lambda}^+}^{\chi_+}e^{(\lambda+2\pi i \frac{c}{h})(t-t_0)}e^{-2\pi i\frac{1}{h}(y-x)}G_{\lambda}(x,y)\ d\lambda\\[0.2cm]
&=&\int_{\overline{\lambda}^+}^{\chi_+}e^{\lambda(t-t_0)}G_{\lambda}(x,y)\ d\lambda.\end{array}\end{equation}
Hence, if we integrate the function $e^{\lambda(t-t_0)}G_{\lambda}(hj+ct,hj_0+ct_0)$ along the rectangle with edges $-\frac{\overline{\lambda}}{2}-i\frac{\pi c}{h},-\frac{\overline{\lambda}}{2}+i\frac{\pi c}{h},\chi-i\frac{\pi c}{h}$ and $\chi+i\frac{\pi c}{h}$, then the integrals from $\chi-i\frac{\pi c}{h}$ to $-\frac{\overline{\lambda}}{2}-i\frac{\pi c}{h}$ and from $-\frac{\overline{\lambda}}{2}+i\frac{\pi c}{h}$ to $\chi+i\frac{\pi c}{h}$ cancel each other out. In particular, again writing $x=hj+ct,y=hj_0+ct_0$, the residue theorem implies
\begin{equation}\begin{array}{lcl}\mathcal{G}_{j}^{j_0}(t,t_0)&=&\frac{-h}{2\pi i}\int_{\chi-i\frac{\pi c}{h}}^{\chi+i\frac{\pi c}{h}}e^{\lambda(t-t_0)}G_{\lambda}(x,y)\ d\lambda\\[0.3cm]
&=&\frac{h}{2\pi i}\int_{-\frac{\overline{\lambda}}{2}-i\frac{\pi c}{h}}^{-\frac{\overline{\lambda}}{2}+i\frac{\pi c}{h}}e^{\lambda(t-t_0)}G_{\lambda}(x,y)\ d\lambda\\[0.4cm]
&&\qquad +\frac{h}{\Omega}\left(\begin{array}{ll}\phi^-(y)\phi^+(x)& \psi^-(y)\phi^+(x)\\ \phi^-(y)\psi^+(x)& \psi^-(y)\psi^+(x)\end{array}\right).\end{array}\end{equation}
Using Lemma \ref{afschattingghzondertilde} we also get the estimate
\begin{equation}\begin{array}{lcl}|\frac{h}{2\pi i}\int_{-\frac{\overline{\lambda}}{2}-i\frac{\pi c}{h}}^{-\frac{\overline{\lambda}}{2}+i\frac{\pi c}{h}}e^{\lambda(t-t_0)}G_{\lambda}(x,y)\ d\lambda| &\leq &\frac{h}{2\pi}\frac{2c\pi}{h}e^{-\overline{\lambda}(t-t_0)}Ke^{-\tilde{\beta}|x-y|},\end{array}\end{equation}
which yields the desired bound (\ref{preciezeafschattingGtilde}).\qed

For any $t\in \R$, we introduce the suggestive notation
\begin{equation}\begin{array}{lcl}\Pi^c(t)&=&\mathcal{E}(t,t)\end{array}\end{equation}
together with
\begin{equation}\begin{array}{lcl}\Pi^s(t)&=&I-\Pi^c(t).\end{array}\end{equation}
Recalling the notation introduced in (\ref{ast}), we set out to show that $\Pi^c(t)\ast \Pi^c(t)=\Pi^c(t)$ and $\Pi^s(t)\ast \Pi^s(t)=\Pi^s(t)$. Later on, we will view these operators as projections that correspond to the center and stable parts of the flow induced by $\mathcal{G}$ respectively.\par

To establish the identity $\Pi^c(t)\ast \Pi^c(t)=\Pi^c(t)$, it suffices to show that
\begin{equation}\label{projectie}\begin{array}{lcl}\left(\begin{array}{ll}\phi^-(x_{j_0})\phi^+(x_{j})& \psi^-(x_{j_0})\phi^+(x_{j})\\ \phi^-(x_{j_0})\psi^+(x_{j})& \psi^-(x_{j_0})\psi^+(x_{j})\end{array}\right) &=&\frac{h}{\Omega}\sum\limits_{i\in\Z}\left(\begin{array}{ll}\phi^-(x_i)\phi^+(x_{j})& \psi^-(x_i)\phi^+(x_{j})\\ \phi^-(x_i)\psi^+(x_{j})& \psi^-(x_i)\psi^+(x_{j})\end{array}\right)\\[0.4cm]
&&\qquad\times \left(\begin{array}{ll}\phi^-(x_{j_0})\phi^+(x_i)& \psi^-(x_{j_0})\phi^+(x_i)\\ \phi^-(x_{j_0})\psi^+(x_i)& \psi^-(x_{j_0})\psi^+(x_i)\end{array}\right),\end{array}\end{equation}
in which $x_i=hi+ct$ for $i\in\Z$. We now write our linear operator in the form
\begin{equation}\begin{array}{lcl} L\Psi(\xi)&=&c\frac{d}{d\xi}\Psi(\xi)+\sum\limits_{j=-\infty}^\infty A_j(\xi)\Psi(\xi+jh),\end{array}\end{equation}
where
\begin{equation}
\begin{array}{lcl}
A_j(\xi)&=&\begin{cases}\left(\begin{array}{ll}
\frac{1}{h^2}\alpha_{|j|}& 0\\ 0& 0
\end{array}\right)&\text{ if }j\neq 0\\[0.4cm]
\left(\begin{array}{ll}-2\frac{1}{h^2}\sum\limits_{k>0}\alpha_k+g_u(\overline{u}(\xi))&1\\ -\rho&\rho\gamma\end{array}\right)&\text{ if }j=0.\end{cases}
\end{array}
\end{equation}
Before we continue, we first prove a small lemma that will help us to relate discrete inner products with their continuous counterparts.
\begin{lemma}\label{ditwasooitpreviouslemma} Assume that (\asref{aannamespuls}{\text{H}}),(\asref{extraaannamespuls}{\text{H}}), (\asref{aannamesconstanten}{\text{H}}), (\asref{aannames}{\text{H}\alpha1}) and (\asref{extraaannames}{\text{H}\alpha}) are satisfied. For all $\xi\in\R$ we have the identity
\begin{equation}\label{vergelijkinginooitpreviouslemma}\begin{array}{lcl} c\left(\begin{array}{l}\phi^-(\xi)\phi^+(\xi)\\ \psi^-(\xi)\psi^+(\xi)

\end{array}\right)&=&
\sum\limits_{j=-\infty}^\infty \int_0^{hj}B(\xi+\theta-hj)A_j(\xi+\theta-hj)\Phi^+(\xi+\theta) d\theta,\end{array}\end{equation}
where
\begin{equation}
\begin{array}{lcl}
B(\xi)&=&\left(\begin{array}{ll}\phi^-(\xi)&0\\ 0& \psi^-(\xi)

\end{array} \right).
\end{array}
\end{equation}
\end{lemma}
\vspace*{4pt}\noindent\textit{Proof.} Our strategy is to differentiate both sides of (\ref{vergelijkinginooitpreviouslemma}) and to show their derivatives are equal. Starting with the first component, we pick $N\in\Z_{>0}\cup\{\infty\}$ and write
\begin{equation}\begin{array}{lcl}D(N)&:=&\frac{d}{d\xi}\sum\limits_{j=-N}^N \int_0^{hj} \phi^-(\xi+\theta-hj)\Big[A_j(\xi+\theta-hj)\Phi^+(\xi+\theta)\Big]^{(1)}d\theta. \end{array}\end{equation}
For finite $N$, we may compute
\begin{equation}\begin{array}{lcl}D(N) &=&\frac{d}{d\xi}\sum\limits_{j=-N}^N \int_{\xi-hj}^\xi \phi^-(\theta)\Big[A_j(\theta)\Phi^+(\theta+hj)\Big]^{(1)}d\theta\\[0.2cm]
&=&\sum\limits_{j=-N}^N \phi^-(\xi)\Big[A_j(\xi)\Phi^+(\xi+hj)\Big]^{(1)}\\[0.2cm] &&\qquad -\sum\limits_{j=-N}^N \phi^-(\xi-hj)\Big[A_j(\xi-hj)\Phi^+(\xi)\Big]^{(1)}.\end{array}\end{equation}
Now for $j>0$ we have $|A_j(\xi)\Phi^+(\xi+hj)|\leq \frac{1}{h^2}|\alpha_j|$, so the partial sums converge uniformly. Hence, it follows that
\begin{equation}\begin{array}{lcl}D(\infty)&=&\sum\limits_{j=-\infty}^\infty\phi^-(\xi)\Big[A_j(\xi)\Phi^+(\xi+hj)\Big]^{(1)}\\[0.2cm] &&\qquad -\sum\limits_{j=-\infty}^\infty \phi^-(\xi-hj)\Big[A_j(\xi-hj)\Phi^+(\xi)\Big]^{(1)}\\[0.2cm]
&=&\phi^-(\xi)c(\phi^+)'(\xi)+c
(\phi^-)'(\xi)\phi^+(\xi)\\[0.2cm]
&=&c(\phi^-\overline{u}')'(\xi),\end{array}\end{equation}
since $\Phi^+\in\ker(L)$ and $\Phi^-\in\ker(L^*)$.\par

We now set out to show that both sides of (\ref{vergelijkinginooitpreviouslemma}) converge to zero as $\xi\rightarrow\infty$. Pick $\epsilon>0$ and let $N\in\Z_{>0}$ be large enough to ensure that
\begin{equation}\begin{array}{lcl}\sum\limits_{j\geq N}\frac{1}{h^2}j|\alpha_j|&\leq &\frac{\epsilon}{4(1+\nrm{\phi^-}_\infty)\nrm{\Phi^+}_\infty
}.\end{array}\end{equation} In addition, let $\Xi$ be large enough to have
\begin{equation}\begin{array}{lcl}|\phi^-(\xi)|&\leq &\frac{\epsilon}{4(1+\sum\limits_{j=-N}^N |j\alpha_{|j|}|)\nrm{\Phi^+}_\infty}\end{array}\end{equation}
for all $\xi\geq \Xi-N$. This $\Xi$ exists since $\phi^-\in H^1$. For such $\xi$ we may estimate
\begin{equation}\begin{array}{lcl}|\sum\limits_{j=-\infty}^\infty \int_0^{hj} \phi^-(\xi+\theta-hj)\Big[A_j(\xi+\theta-hj)\Phi^+(\xi+\theta)\Big]^{(1)}d\theta|&<&\epsilon,\end{array}\end{equation}
which allows us to compute
\begin{equation}\begin{array}{lcl}&&\lim\limits_{\xi\rightarrow\infty}\sum\limits_{j=-\infty}^\infty \int_0^{hj} \phi^-(\xi+\theta-hj)\Big[A_j(\xi+\theta-hj)\Phi^+(\xi+\theta)\Big]^{(1)}d\theta=0\\[0.2cm]
&=&\lim\limits_{\xi\rightarrow\infty}c\phi^-(\xi)\phi^+(\xi).\end{array}\end{equation}
With that we have proved our claim. Furthermore, we can repeat the arguments above to obtain
\begin{equation}\begin{array}{lcl}c\psi^-(\xi)\psi^+(\xi)&=&
\sum\limits_{j=-\infty}^\infty \int_0^{hj} \psi^-(\xi+\theta-hj) \Big[A_j(\xi+\theta-hj)\Phi^+(\xi+\theta)\Big]^{(2)}d\theta.\end{array}\end{equation}\qed

We are now ready to show that $\Pi^c(t)\ast \Pi^c(t)=\Pi^c(t)$ and $\Pi^s(t)\ast \Pi^s(t)=\Pi^s(t)$. This result is based on the first part of \cite[Lemma 2.9]{HJHSTBFHN}.
\begin{lemma}\label{projections} Assume that (\asref{aannamespuls}{\text{H}}),(\asref{extraaannamespuls}{\text{H}}), (\asref{aannamesconstanten}{\text{H}}), (\asref{aannames}{\text{H}\alpha1}) and (\asref{extraaannames}{\text{H}\alpha}) are satisfied. Then $\Pi^c(t)\ast \Pi^c(t)=\Pi^c(t)$ and $\Pi^s(t)\ast \Pi^s(t)=\Pi^s(t)$ for all $t\in\R$.\end{lemma}
\noindent\textit{Proof.} For $k\in\Z$ we write $x_k=hk+ct$. In addition, for notational convenience we write $B_j(\theta)=\Big[A_j(\theta)\Phi^+(\theta+hj)\big]^{(1)}$ for $j\in\Z$ and $\theta\in\R$. Using the results from Lemma \ref{ditwasooitpreviouslemma} we may compute
\begin{equation}\begin{array}{lcl}c\sum\limits_{k=-\infty}^\infty\phi^-(x_k)\phi^+(x_k)&=&\sum\limits_{k=-\infty}^\infty\sum\limits_{j=-\infty}^\infty \int\limits_0^{hj}\phi^-(x_k+\theta-hj) B_j(x_k+\theta-hj)d\theta\\[0.2cm]
&=&\sum\limits_{k=-\infty}^\infty \sum\limits_{j=-\infty}^\infty \int\limits_{x_{k-j}}^{x_k} \phi^-(\theta) B_j(\theta) d\theta\\[0.2cm]
&=&\sum\limits_{j=-\infty}^\infty\int\limits_{-\infty}^\infty  j\phi^-(\theta)B_j(\theta)d\theta,\end{array}\end{equation}
where we were allowed to interchange the two infinite sums because
\begin{equation}\begin{array}{lcl}\big|\sum\limits_{k=-N}^N \int\limits_{x_{k-j}}^{x_k} \phi^-(\theta)B_j(\theta)d\theta\big|&\leq &\big|\int\limits_{-\infty}^\infty  j\phi^-(\theta)B_j(\theta)d\theta\big|\\[0.2cm]
&\leq &\nrm{\phi^-}_1 \frac{1}{h^2}|j\alpha_{|j|}|\nrm{\phi^+}_\infty\end{array}\end{equation}
holds for all $N\in\Z_{>0}$ and $j\in\Z$. This expression is summable over $j$, allowing us to apply Lebesgue's theorem. On the other hand, we have
\begin{equation}\begin{array}{lcl}c\int\limits_{-\infty}^\infty \phi^-(\xi)\phi^+(\xi)\ d\xi&=&\int\limits_{-\infty}^\infty \sum\limits_{j=-\infty}^\infty \int\limits_0^{hj} \phi^-(\xi+\theta-hj)B_j(\xi+\theta-hj)d\theta \ d\xi\\[0.2cm]
&=&\sum\limits_{j=-\infty}^\infty \int\limits_0^{hj}\int\limits_{-\infty}^\infty\phi^-(\xi+\theta-hj)B_j(\xi+\theta-hj)\ d\xi d\theta \\[0.2cm]
&=&\sum\limits_{j=-\infty}^\infty \int\limits_0^{hj}\int\limits_{-\infty}^\infty\phi^-(\xi-hj)B_j(\xi-hj) d\xi d\theta \\[0.2cm]
&=&\sum\limits_{j=-\infty}^\infty hj\int\limits_{-\infty}^\infty\phi^-(\xi-hj)B_j(\xi-hj)d\xi. \end{array}\end{equation}
Interchanging the integral with the sum was allowed since $\phi^-$ and $\phi^+$ decay exponentially, say $|\phi^-(x)|\leq \kappa e^{-\delta|x|}$ and $|\phi^+(x)|\leq \kappa e^{-\delta|x|}$. In particular, for each $N\in\Z_{>0}$ and each $\xi\in\R$ we have
\begin{equation}\begin{array}{lcl}|\sum\limits_{j=-N}^N \int\limits_0^{hj} \phi^-(\xi+\theta-hj)B_j(\xi+\theta-hj)d\theta|&\leq &\sum\limits_{j=-\infty}^\infty h \kappa^2 e^{-\delta|\xi|}|j\alpha_{|j|}|\nrm{\Phi^+}_\infty,\end{array}\end{equation}
which is integrable in $\xi$. Furthermore, the interchanging of the two integrals was allowed, since by the exponential decay of $\phi^-$ we also see that for each $j\in\Z,\xi\in\R$ and $\theta\in(0,hj)$ we have
\begin{equation}\begin{array}{lcl}|\phi^-(\xi+\theta-hj)B_j(\xi+\theta-hj)|&\leq &\kappa e^{-\delta|\xi+\theta-hj|}|\alpha_{|j|}|\nrm{\Phi^+}_{\infty}.\end{array}\end{equation}
This is an integrable function for $(\xi,\theta)\in\R\times(0,hj)$, allowing us to apply Fubini's theorem.\par

In particular, we see that
\begin{equation}\begin{array}{lcl}\int\limits_{-\infty}^\infty \phi^-(\xi)\phi^+(\xi)\ d\xi&=&h\sum\limits_{k=-\infty}^\infty\phi^-(x_k)\phi^+(x_k).\end{array}\end{equation}
In the same way we obtain
\begin{equation}\begin{array}{lcl}\int\limits_{-\infty}^\infty \psi^-(\xi)\psi^+(\xi)\ d\xi&=&h\sum\limits_{k=-\infty}^\infty\psi^-(x_k)\psi^+(x_k).\end{array}\end{equation}
By writing out the sums it now follows that indeed (\ref{projectie}) holds.\qed

\vspace*{4pt}\noindent\textit{Proof of Proposition \ref{cor2.8equivalentsterker}.} The calculations above imply that $\mathcal{E}(t,t_0)=\mathcal{E}(t,t_0)\ast\Pi^c(t_0)$, which means that we must also have $\mathcal{E}(t,t_0)\ast\Pi^s(t_0)=0$.\par

Observe that for any $t_0\in\R$, the function $V_j(t):=\left(\begin{array}{l}\phi^+(hj+ct)\\ \psi^+(hj+ct)\end{array}\right)$ is the unique solution to (\ref{definitieGhuh}) with $V_j(t_0)=\left(\begin{array}{l}\phi^+(hj+ct_0)\\ \psi^+(hj+ct_0)\end{array}\right)$. Hence, by the definition of the Green's function $\mathcal{G}(t,t_0)$ we see that
\begin{equation}\begin{array}{lcl}V(t)&=&\mathcal{G}(t,t_0)\ast V(t_0)\end{array}\end{equation}
for all $t\in\R$. Furthermore, we recall that
\begin{equation}
\begin{array}{lcl}
\mathcal{E}_j^{j_0}(t_0,t_0)&=&V_j(t_0)\Phi^-(hj_0+ct_0).
\end{array}
\end{equation}
For $j,j_0\in\Z$ we may hence compute
\begin{equation}\begin{array}{lcl}\Big[\mathcal{G}(t,t_0)\ast\Pi^c(t_0)\Big]_j^{j_0}&=&\sum\limits_{i\in\Z}\mathcal{G}_j^ i(t,t_0)\ast\mathcal{E}_i^{j_0}(t,t_0)\\[0.2cm]
&=&\frac{h}{\Omega}\sum\limits_{i\in\Z}\mathcal{G}_j^i(t,t_0)_i\ast V_i(t_0)\Phi^-(hj_0+ct_0)\\[0.2cm]
&=&\frac{h}{\Omega}V_j(t)\Phi^-(hj_0+ct_0)\\[0.2cm]
&=&\mathcal{E}_j^{j_0}(t,t_0).\end{array}\end{equation}
In particular, we obtain $\mathcal{G}(t,t_0)\ast\Pi^c(t_0)=\mathcal{E}(t,t_0)$ and thus
\begin{equation}\begin{array}{lcl}\tilde{\mathcal{G}}(t,t_0)\ast\Pi^c(t_0)&=&\mathcal{G}(t,t_0)\ast\Pi^c(t_0)-\mathcal{E}(t,t_0)\ast\Pi^c(t_0)\\[0.2cm]
&=&\mathcal{E}(t,t_0)-\mathcal{E}(t,t_0)\\[0.2cm]
&=&0.\end{array}\end{equation}
Therefore, we must have
\begin{equation}\begin{array}{lcl}\mathcal{G}(t,t_0)&=&\mathcal{E}(t,t_0)+\tilde{\mathcal{G}}(t,t_0)\\[0.2cm]
&=&\mathcal{E}(t,t_0)\ast\Big(\Pi^c(t_0)+\Pi^s(t_0)\Big)+\tilde{\mathcal{G}}(t,t_0)\ast\Big(\Pi^c(t_0)+\Pi^s(t_0)\Big)\\[0.2cm]
&=&\mathcal{E}(t,t_0)\ast\Pi^c(t_0)+\tilde{\mathcal{G}}(t,t_0)\ast\Pi^s(t_0),\end{array}\end{equation}
which completes the proof.\qed

\providecommand{\href}[2]{#2}
\providecommand{\arxiv}[1]{\href{http://arxiv.org/abs/#1}{arXiv:#1}}
\providecommand{\url}[1]{\texttt{#1}}
\providecommand{\urlprefix}{URL }
\bibliographystyle{klunumHJ}

\bibliography{ref}

\end{document}